\documentclass[graybox,envcountresetsect,envcountsame]{svmult}

\usepackage{mathptmx}
\usepackage{helvet}
\usepackage{courier}
\usepackage{type1cm}

\usepackage{makeidx}         
\usepackage{graphicx}        
\usepackage{multicol}        
\usepackage[bottom]{footmisc}

\usepackage{amssymb,amsmath,amsfonts}
\usepackage{enumerate}
\usepackage{bm}

\usepackage[frenchb,english]{babel}
\usepackage{dsfont}
\usepackage{changepage}
\usepackage[T1]{fontenc}
\usepackage[utf8]{inputenc}

\usepackage{amsmath}

\usepackage{amsthm}

\usepackage{bookmark}

\usepackage{fancyhdr}

\fancypagestyle{intro}{
\fancyhf{}
\fancyhead[RE]{Introduction}
\fancyhead[LO]{Introduction}
\fancyhead[LE,RO]{\thepage}
}

\fancypagestyle{partI}{
\fancyhf{}
\fancyhead[RE]{\nouppercase Part I. \nouppercase\leftmark}
\fancyhead[LO]{\nouppercase\rightmark}
\fancyhead[LE,RO]{\thepage}
}

\fancypagestyle{partII}{
\fancyhf{}
\fancyhead[RE]{\nouppercase Part II. \nouppercase\leftmark}
\fancyhead[LO]{\nouppercase\rightmark}
\fancyhead[LE,RO]{\thepage}
}

\fancypagestyle{partIII}{
\fancyhf{}
\fancyhead[RE]{\nouppercase Part III. \nouppercase\leftmark}
\fancyhead[LO]{\nouppercase\rightmark}
\fancyhead[LE,RO]{\thepage}
}

\fancypagestyle{partIV}{
\fancyhf{}
\fancyhead[RE]{\nouppercase Part IV. \nouppercase\leftmark}
\fancyhead[LO]{\nouppercase\rightmark}
\fancyhead[LE,RO]{\thepage}
}

\fancypagestyle{appendixI}{
\fancyhf{}
\fancyhead[RE]{\nouppercase Part I. Appendix}
\fancyhead[LO]{\nouppercase\rightmark}
\fancyhead[LE,RO]{\thepage}
}

\fancypagestyle{appendixIII}{
\fancyhf{}
\fancyhead[RE]{\nouppercase Part III. Appendix}
\fancyhead[LO]{\nouppercase\rightmark}
\fancyhead[LE,RO]{\thepage}
}

\fancypagestyle{appendixIV}{
\fancyhf{}
\fancyhead[RE]{\nouppercase Part IV. Appendix}
\fancyhead[LO]{\nouppercase\rightmark}
\fancyhead[LE,RO]{\thepage}
}

\fancypagestyle{referencesGI}{
\fancyhf{}
\fancyhead[RE]{\nouppercase Introduction. References}
\fancyhead[LO]{\nouppercase Introduction. References}
\fancyhead[LE,RO]{\thepage}
}

\fancypagestyle{referencesI}{
\fancyhf{}
\fancyhead[RE]{\nouppercase Part I. References}
\fancyhead[LO]{\nouppercase Part I. References}
\fancyhead[LE,RO]{\thepage}
}

\fancypagestyle{referencesII}{
\fancyhf{}
\fancyhead[RE]{\nouppercase Part II. References}
\fancyhead[LO]{\nouppercase Part II. References}
\fancyhead[LE,RO]{\thepage}
}

\fancypagestyle{referencesIII}{
\fancyhf{}
\fancyhead[RE]{\nouppercase Part III. References}
\fancyhead[LO]{\nouppercase Part III. References}
\fancyhead[LE,RO]{\thepage}
}

\fancypagestyle{referencesIV}{
\fancyhf{}
\fancyhead[RE]{\nouppercase Part IV. References}
\fancyhead[LO]{\nouppercase Part IV. References}
\fancyhead[LE,RO]{\thepage}
}

\usepackage{numprint}

\usepackage[left]{lineno}

\usepackage{rotating}

\makeindex             

\newcommand{\F} {{ \mathcal F }}

\newcommand{\N} {{{\mathbb N}}}
\newcommand{\m}  {\widetilde{\mathcal{M}}}
\newcommand{\n}  {\widetilde{\mathcal{N}}}

\newcommand{\M} {{ \mathcal M }}
\newcommand{\NN} {{ \mathcal N }}

\newcommand{\V} {{ \mathcal V}}

\newcommand{\BB} {{ \mathcal B}}

\newcommand{\eps}    {\varepsilon}

\newcommand{\J}     {{\mathcal J}}
\renewcommand{\L}   {{\mathcal L}}

\newcommand{\q}  {\overline{q}}

\newcommand{\bpf}{\begin{proof}}
\newcommand{\epf}{\end{proof}}

\renewcommand{\E}{\mathbb{E}}

\newcommand{\Z}{\mathbb{Z}}
\newcommand{\R}{\mathbb{R}}

\renewcommand{\I}{\textsc{I}}

\renewcommand{\P}{\mathbb{P}}

\def\interior#1{\smash{\mathop{#1}\limits^{\lower1pt\hbox{$\scriptscriptstyle\circ$}}}}

\spnewtheorem{assumption}[theorem]{Assumption}{\normalfont}{\normalfont}
\spnewtheorem{notation}[theorem]{Notation}{\itshape}{\rmfamily}

\allowdisplaybreaks

\numberwithin{equation}{section}
\numberwithin{theorem}{section}
\numberwithin{figure}{section}
\numberwithin{section}{chapter}
\numberwithin{table}{section}

\makeatletter
\def\blfootnote{\xdef\@thefnmark{}\@footnotetext}
\makeatother

\begin{document}

\vspace{5cm}

\thispagestyle{empty}

\noindent{\LARGE Tom Britton}
\vspace{2mm}

\noindent{\LARGE Etienne Pardoux}
\vspace{2mm}

\noindent{\LARGE (Eds.)}
\vspace{2mm}

\vspace{4cm}

\noindent{\Huge {Stochastic Epidemic Models with Inference}}

\vspace{4cm}

\noindent{\LARGE With Contributions by:}\\

\noindent{\large Frank Ball, Tom Britton, Catherine Lar\'edo, Etienne Pardoux,\\ \noindent David Sirl and Viet Chi Tran}

\newpage




\frontmatter

\tableofcontents

\mainmatter


\part[Stochastic epidemics in a homogeneous community\\ Tom Britton and Etienne Pardoux]{Stochastic epidemics in a homogeneous community\\ \vspace{1cm} \normalsize{Tom Britton\blfootnote{\hspace{-2mm} Tom Britton\\ \email{tom.britton@math.su.se}} and Etienne Pardoux\blfootnote{\hspace{-2mm} Etienne Pardoux\\ \email{etienne.pardoux@univ-amu.fr}}}}

\setcounter{chapter}{0}

\chapter*{Introduction}
\addcontentsline{toc}{chapter}{Introduction}

\pagestyle{partI}

In this Part I of the lecture notes our focus lies exclusively on stochastic epidemic models for a homogeneously mixing community of individuals being of the same type. The important extensions allowing for different types of individuals and allowing for non-uniform mixing behaviour in the community is left for later  parts  in the Notes.

In Chapter 1, 
we present the stochastic SEIR  epidemic model, derive some important properties of it, in particular for the beginning of an outbreak. Motivated by mathematical tractability rather than realism we then study in Chapter 2
the special situation where the model is Markovian, and derive additional results for this sub-model.

What happens later on in the outbreak will depend on our model assumptions, which in turn depend on the scientific questions. In Chapter 3 
we focus on short-term outbreaks, when it can be assumed that the community is fixed and constant during the outbreak; we call these models closed models. In Chapter 4 
we are more interested in long-term behaviour, and then it is necessary to allow for influx of new individuals and that people die, or to include return to susceptibility. Such so-called open population models are harder to analyse -- for this reason we stick to the simpler class of Markovian models.  In this chapter we consider situations where the deterministic model has a unique stable equilibrium, and use both the central limit theorem and large deviation techniques to predict the time at which the disease goes extinct in the population.

The Notes end with an extensive Appendix, giving some relevant probability theory used in the main part of the Notes and also solutions to most of the exercises being scattered out in the different chapters.

\chapter{Stochastic Epidemic Models}\label{TB-EP_chap_StochMood}

This first chapter introduces some basic facts about stochastic epidemic models. We consider the case of a closed community, i.e.\  without influx of new susceptibles or mortality. In particular, we assume that the size of the population is fixed, and
that the individuals who recover from the illness are immune and do not become susceptible again.
We describe the general class of stochastic epidemic models, and define the basic reproduction number, which allows one to determine whether or not a major epidemic may start from the initial infection of a small number of  individuals.
We then approximate the early stage of an outbreak with the help of a branching process, and from this obtain the distribution of the final size (i.e.\ the total number of individuals who ever get infected) in case of a minor outbreak. Finally we discuss the impact of vaccination.

The important problem of estimating model parameters from (various types of) data is left to Part IV of the current volume (also discussed in Chapter 4 of Part III). Here we assume the model parameters to be known.

\section[The stochastic SEIR epidemic model in a closed homogeneous community]{The stochastic SEIR epidemic model in a closed homogeneous community}\label{Sec_Mod-def}

\subsection{Model definition}

Consider a closed population of $N+1$ individuals ($N$ is the number of initially susceptible). At any point in time each individual is either susceptible, exposed, infectious or recovered. Let $S(t),\ E(t),\ I(t)$ and $R(t)$ denote the numbers of individuals in the different states at time $t$ (so $S(t)+ E(t)+I(t)+R(t)=N+1$ for all $t$). The epidemic starts at $t=0$ in a specified state, often the state with one infectious individual, called the index case and thought of as being externally infected, and the rest being susceptible: $(S(0), E(0), I(0), R(0))=(N,0,1,0)$.

\begin{definition}\label{def:SEIR}
While infectious, an individual has infectious contacts according to a Poisson process with rate $\lambda$. Each contact is with  an individual chosen uniformly at random from the rest of the population, and if the contacted individual is susceptible he/she becomes infected -- otherwise the infectious contact has no effect. Individuals that become infected are first latent (called exposed) for a random duration $L$ with distribution $F_L$, then they become infectious for a duration $I$ with distribution $F_I$, after which they become recovered and immune for the remaining time. All Poisson processes, uniform contact choices, latent periods and infectious periods of all individuals are defined to be mutually independent.
\end{definition}

The epidemic goes on until the first time $\tau$ when there are no exposed or infectious individuals, $E(\tau)+I(\tau)=0$. At this time no further individuals can get infected so the epidemic stops. The final state hence consists of susceptible and recovered individuals, and we let $Z$ denote the \emph{final size}, i.e.\ the number of infected (by then recovered) individuals at the end of the epidemic excluding the index case(s): $Z=R(\tau)-I(0)=N-S(\tau)$. The possible values of $Z$ are hence $0,\dots ,N$.

\subsection{Some remarks, submodels and model generalizations}\label{sec:subrem}

Quite often the rate of ``infectious contacts'' $\lambda$ can be thought of as a product of a rate $c$ at which the infectious individual has contact with others, and the probability $p$ that such a contact results in infection given that the other person is susceptible, so $\lambda=c p$. As regards to the propagation of the disease it is however only the product $\lambda$ that matters and since fewer parameters is preferable we keep only $\lambda$.

The rate of infectious contacts is $\lambda$, so the rate at which one infectious has contact with a specific other individual is $\lambda/N$ since each contact is with a uniformly chosen other individual.

First we will look what happens in a very small community/group, but the main focus of these notes is for a large community, and the asymptotics are hence for $N\to\infty$. The parameters of the model, the infection rate $\lambda$, and the latent and infectious periods $L$ and $I$, are defined independently of $N$, but the epidemic is highly dependent on $N$ so when this needs to be emphasized we equip the corresponding notation with an $N$-index, e.g.\ $S^N(t)$ and $\tau^N$ which hence is not a power.

Some special cases of the model have received special attention in the literature. If both $L$ and $I$ are exponentially distributed (with rates $\nu$ and $\gamma$ say), the model is Markovian which simplifies the mathematical analysis a great deal. This model is called the \emph{Markovian SEIR}. If $L\equiv 0$ and $I\sim \mathrm{Exp} (\gamma)$ then we have the \emph{Markovian SIR} (whenever there is no latency period the model is said to be SIR) which is better known under the unfortunate name the \emph{General stochastic epidemic}. Another special case of the stochastic SEIR model is where the infectious period $I$ is non-random. Also here there is a underlying mathematical reason -- when the duration of the infectious period is non-random and equal to $\iota$ say, then an infectious individual has infectious contacts with each other individual at rate $\lambda/N$ during a non-random time implying that the number of contacts with different individuals are \emph{independent}. Consequently, an infectious individual has infectious contacts with each other individual independently with probability $p=1-e^{-(\lambda/N)\iota}$, so the total number of contacts is Binomially distributed, and in the limit as $N\to\infty $ the number of infectious contacts an individual has is Poisson distributed with mean $\lambda\iota$. If further the  latent period is long in comparison to the infectious period then it is possible to identify the infected individuals in terms of \emph{generations}: the first generation are the index cases, the second generation those who were infected by the index case(s), and so one. When the model is described in this discrete time setting and  individuals infect different individuals independently with probability $p$,  this model is the well-known Reed--Frost model named after its inventors Reed and Frost.

The two most studied special cases are hence when the infectious period is exponentially distributed and when it is nonrandom. For real infectious diseases none of these two extremes apply, for influenza for example, the infectious period is believed to be about 4 days, plus or minus one or two days. If one has to choose between these choices a nonrandom infectious period is probably closer to reality.

\begin{figure*}[h]
\begin{center}
\includegraphics[width=0.85\textwidth, height=0.3\textheight]{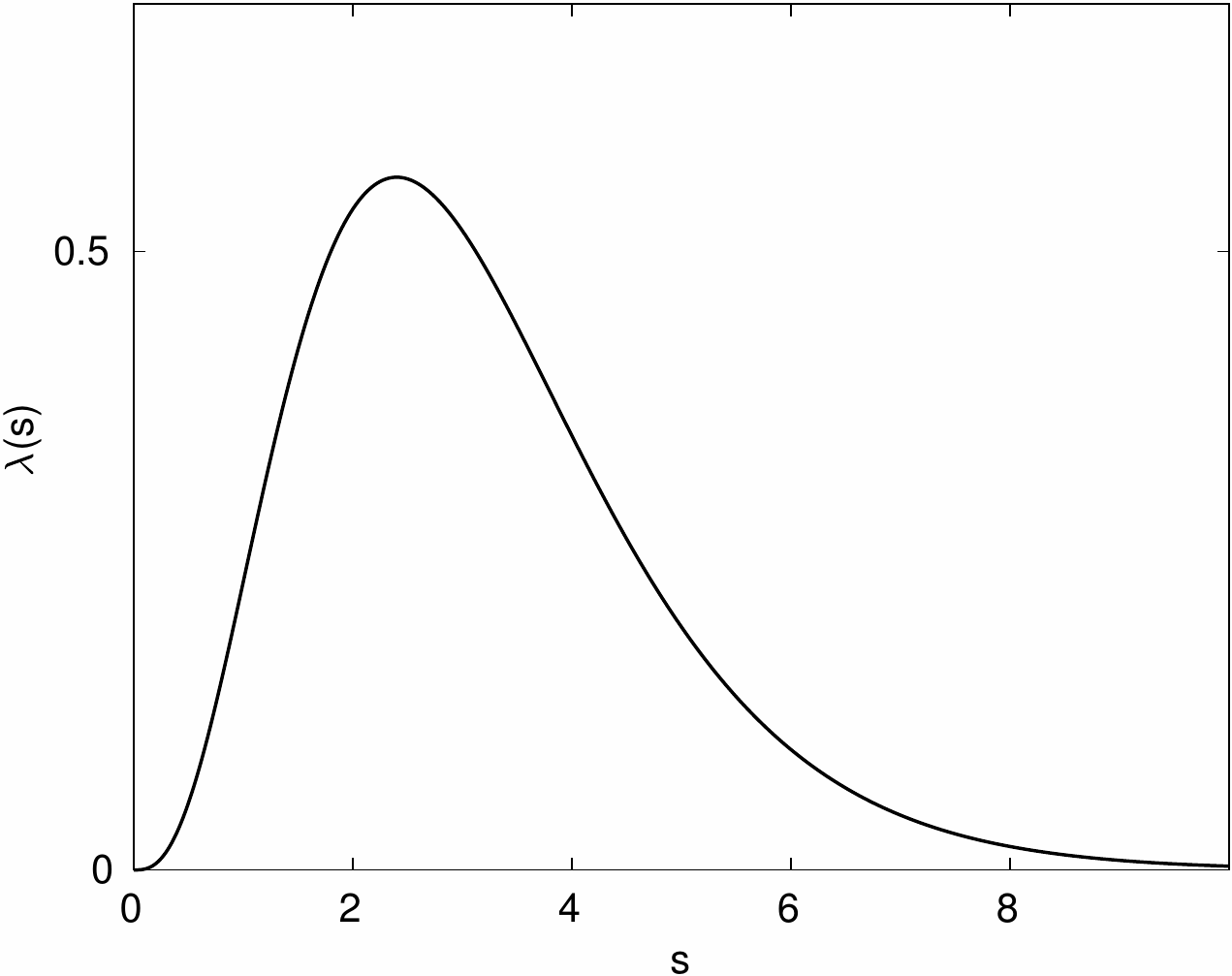}
\end{center}
\caption{Plot of a possible infectivity curve $\lambda (s)$. The time $s$ denotes the time since infection in unit of days.}
\label{fig-infcurve}
\end{figure*}

The stochastic SEIR model in a closed homogeneous community may of course also be generalized towards more realism. Two such extensions have already been mentioned: allowing for individuals to die and new ones to be born, and allowing for some social structures. Some such extensions will be treated in the other articles of the current lecture notes but not here. But even when assuming a closed homogeneously mixing community of homogeneous individuals it is possible to make the model more realistic. The most important such generalization is to let the rate of infectious contact vary with time since infection. The current model assumes there are no infectious contacts during the latent state, and then, suddenly when the latent period ends, the rate of infectious contact becomes $\lambda$ until the infectious period ends when it suddenly drops down to 0 again. In reality, the infectious rate is usually a function $\lambda(s)$ $s$ time units after infection. In most situations $\lambda (s)$ is very small initially (corresponding to the latency period) followed by a gradual increase for some days, and then $\lambda(s)$ starts decaying down towards 0 which it hits when the individual has recovered completely (see Figure \ref{fig-infcurve} for an example where infectivity starts growing after one day and is more or less over after one week).
The function $\lambda (s)$ could be the same for all individuals, or it may be random and hence a stochastic process, i.i.d.\ for different individuals. As regards to the temporal dynamics of the epidemic process, the functional form of $\lambda(s)$ is important, and also its random properties in case it is random. If one is only interested in the final size $\tau$, it is however possible to show that all that affects the final size is the accumulated force of infection, i.e.\ the distribution of $\int_0^\infty \lambda(s)ds$. In particular, if we let $\lambda I$ in the stochastic SEIR model have the same distribution as $\int_0^\infty \lambda(s)ds$ in the more general model, then the two models have the same final size distribution. In that sense, the extended model can be included in the stochastic SEIR model.

\subsection{Two key quantities: $R_0$ and the escape probability}

The most important quantity for this, as well as most other epidemic models, is the \emph{basic reproduction number} (sometimes ``number'' is replaced by ``ratio'') and denoted by $R_0$. In more complicated models its definition and interpretation are sometimes debated, but for the present model it is quite straightforward: $R_0$ denotes the mean number of infectious contacts a typical infected has during the early stage of an outbreak. As the population under consideration is becomes large, this number will coincide with the mean number of infections caused by a typical infected during the early stages of an outbreak. We derive an expression for $R_0$, but before that we should consider its important threshold value of 1. If $R_0>1$ this means that on average an infected infects more than one individual in the beginning of an epidemic. Then the index case on average is replaced by more than one infected, who in turn each are replaced by more than one infected and so on. This clearly suggests that a big community fraction can become infected. If on the other hand $R_0\le 1$, then the same reasoning suggests that there will never be a big community outbreak. Those results hold true which we prove in Section \ref{sec-early-stage} (Corollaries \ref{Th-early-stage-sub} and \ref{Th-early-stage-sup}).

In applications the basic reproduction number $R_0$ is a central quantity of interest. Many studies of disease outbreaks contain estimates of $R_0$ for a specific disease and community, together with modeling conclusions about preventive measures which, if put into place, will reduce the reproduction number $R$ down to below the critical value of 1 when an outbreak is no longer possible (e.g.\ Fraser et al.\ \cite{F2009I}).

Let us now derive an expression for $R_0$. An infected individual has infectious contacts only when infectious, and when in this state the individual has infectious contacts at rate $\lambda$. This means that the expected number of infectious contacts equals
\begin{equation}
R_0=\E(\lambda I)=\lambda \iota .\label{R_0}
\end{equation}
Sometimes the rate $\lambda$ of having infectious contacts is replaced by an over-all rate of contact $c$  multiplied by the probability $p$ of a contact leading to infection, so $\lambda=c\, p$ and $R_0=c\, p\,\iota$ (cf.\ the first lines of the above Subsection \ref{sec:subrem}, Anderson and May \cite{AM91I} and Giesecke \cite{G2017I}).

Another key quantity appearing later several times is the probability for a given susceptible to escape getting infected from a specific infective. The instantaneous infectious force from the infective to this specific susceptible is $\lambda/N$, and the random duration of the infectious period is $I$. Conditional upon $I=x$, the escape probability is hence $e^{-(\lambda/N)x}$, and the unconditional probability to escape infection is therefore
\begin{equation}
 \P (\text{escape infection from an infective})=\E(e^{-\lambda I/N})=\psi_I(-\lambda/N),\label{escape-prob}
\end{equation}
where $\psi_I(b)=\E(e^{bI})$ is the moment generating function of the infectious period (so $\psi(-b)$ is the Laplace transform -- in Part II in this volume the Laplace transform has a separate notation, $\phi$, so $\phi(b)=\psi(-b)$).

\begin{exercise}\label{x2.1}
Consider the Markovian SEIR epidemic in which $\lambda=1.8$, $\nu=2$ and $\gamma=1$ in a village of size $N=100$,  (parameters inspired by Ebola with weeks as time unit). Compute $R_0$ and the escape probability.
\end{exercise}

\begin{exercise}\label{x2.2}
Repeat the previous exercise, but now for the Reed--Frost epidemic with $\lambda=1.8$, $L\equiv 2$ and $I\equiv \iota =1$ in a village of size $N=100$,  (perhaps having more realistic distributions than in the previous exercise).
\end{exercise}

\section{The early stage of an outbreak}\label{sec-early-stage}

We now consider the situation where the community size $N$ is large and study the stochastic SEIR epidemic in the beginning of an outbreak. By ``beginning'' we mean that less than $k=k(N)$ individuals have been infected. Recall from the model definition that infectious individuals have infectious contacts with others independently, each infective at rate $\lambda$. The dependence only appears because individuals can only get infected once, so if an individual has already received an infectious contact, then future infectious contacts with that individual no longer result in someone getting infected. However, in the beginning of an outbreak in a large community it is very unlikely that two infectives happen to have infectious contacts with the same individual. This suggests that during the early phase of an outbreak, infectives infect new individuals more or less independently. This implies that the number of infected can be approximated by a branching process in the beginning of an outbreak, where ``being born'' corresponds to having been infected, and ``giving birth'' corresponds to infecting someone. The current section is devoted to making this approximation rigorous, and thus obtaining asymptotic results for the epidemic in regards to having a minor versus a major outbreak. In the next section  this approximation is exploited in order to determine the distribution of the final size in the case of a minor outbreak. If the epidemic takes off, which happens in the case of a major outbreak, then the approximation that individuals infect others independently breaks down. What happens in this situation is treated in later sections.

First we define the approximating branching process and derive some properties of it. After this we show rigorously that, as $N\to\infty$,  the initial phase of the epidemic process converges to the initial phase of the branching process by using an elegant coupling technique.

The approximating branching process is defined similarly to the epidemic. A newborn individual is first unable to give birth to new individuals for a period with duration $L$ (this period might be denoted childhood in the branching process setting). After this childhood, the individual enters the reproductive stage which last for $I$ units of time. During this period individuals give birth to new individuals at rate $\lambda$ (randomly in time according to a Poisson process with rate $\lambda$). Once the reproductive stage has terminated the individual dies (or at least cannot reproduce and hence plays no further role).

The number of offspring of an individual, $X$, depends on the duration of the reproductive stage $I$. Conditional upon $I=y$, the number of births follow the Poisson distribution  $\text{Poi}(\lambda y)$, so the unconditional distribution of number of offspring is mixed-Poisson, written as $X\sim \text{MixPoi}(\lambda I)$, where $I$ has distribution $F_I$.

If we forget calendar time, and simply study the number of individuals born in each generation, then  our branching process is a Bienaym\'e--Galton--Watson process with offspring distribution being $\text{MixPoi}(\lambda I)$. The mean number of children/offspring equals $m=\E (X)=\E (\E (X|I))=\E (\lambda I)=\lambda\iota$.

\begin{exercise} \label{xEarly.1}
Compute the offspring distribution $\P (X=x)$ explicitly for the two cases: (i) where the infectious period is non-random, $I\equiv \iota$, corresponding to the continuous--time version of the Reed--Frost epidemic; and (ii) for the Markovian SEIR where $I$ is exponential with mean $\iota$.
\end{exercise}

We now show an elegant coupling construction which we will use to show that the epidemic and branching process have similar distributions in the beginning. To this end we define the approximating branching process as well as all epidemics, i.e.\ for each
$N=1,2,\ldots$, on the same probability space.
To this end, let $L_0, L_1,  \dots$ be i.i.d.\ latent periods having distribution $F_L$, and similarly let $I_0, I_1,  \dots$ be i.i.d.\ infectious periods having distribution $F_I$. Further, let $\xi_0(\cdot), \xi_1(\cdot), \dots$ be i.i.d.\ Poisson processes having intensity $\lambda$, and let $U_1,U_2,\dots$ be i.i.d.\ $U(0,1)$ random variables. All random variables and Poisson processes are assumed to be mutually independent. These will be used to construct the branching process as well as the stochastic SEIR epidemic for each $N$ as follows.

\begin{definition}\label{DefApprBrPr}{\rm The approximating branching process}.
At time $t=0$ there is one new born ancestor having label $0$. Let the ancestor have childhood length $L_0$ and reproductive stage for a duration $I_0$ (so the ancestor dies at time $L_0+I_0$), during which the ancestor gives birth at the time points of the Poisson process $\xi_0(\cdot)$. If the jump times of this Poisson process are denoted $T_{0,1}<T_{0,2}< ...$ and $X_0$ denotes the number of jumps prior to $I_0$, then the ancestor gives birth at the time points $L_0+T_{0,1}, \dots ,L_0+T_{0,X_0}$ (the set is empty if $X_0=0$).  The first born individual is given label 1, and having childhood period $L_1$, reproductive period $I_1$ and birth process $\xi_1(\cdot)$. This individual gives birth according to the same rules (starting the latency period at time $L_0+T_{0,1}$), and the next individual born, either to individual 0 or 1, is given label 2 and variables $L_2, I_2$ and birth process $\xi_2(\cdot)$, and so on. This defines the branching process, and we let $L(t), I(t), R(t)$ respectively denote the numbers of individuals in the childhood state, in the reproductive state and dead, respectively, at time $t$. The total number of individuals born up to time $t$, excluding the ancestor/index case, is denoted by $Z(t)=L(t)+I(t)+R(t)-1$ in the branching process, and the ultimate number ever born, excluding the ancestor, is denoted by $Z$ which may be finite or infinite.
\end{definition}

We now define the epidemic for any fixed $N$ (in the epidemic childhood corresponds to latent and reproductive stage to being infectious). This is done similarly to the branching process with the exception that we now keep track of which individuals who get infected using the uniform random variables $U_1,U_2,\dots$.

\begin{definition}{\rm The stochastic SEIR epidemic with $N$ initial susceptibles}.
We label the $N+1$ individuals $0,1, \dots , N$, with the index case having label 0 and the others being labelled arbitrarily. As for the branching process, the index case is given latency period $L_0$, infectious period $I_0$ and contact process $\xi_0(\cdot)$ and the epidemic is started at time $t=0$. The infectious contacts of the index case occur at the time points $L_0+T_{0,1}, \dots ,L_0+T_{0,X_0}$. The first infectious contact is with individual $[U_1N]+1$, the integer part of $NU_1$ plus 1 (this picks an individual uniformly among $1,\dots ,N$). This individual, $k$ say, then becomes infected (and latent) and is given latent period, infectious period and contact process $L_1, I_1$ and $\xi_1(\cdot)$. The next infectious contact (from either the index case or individual $k$) will be with individual $[U_2N]+1$. If the contacted person is individual $k$ then nothing happens, but otherwise this new individual gets infected (and latent), and so on. Infectious contacts only result in infection if the contacted individual is still susceptible. When a contact is with an already infected individual the branching process has a birth whereas there is no infection in the epidemic -- we say a ``ghost'' was infected when comparing with the branching process. Descendants of all ghosts are also ignored in the epidemic.  The epidemic goes on until there are no latent or infectious individuals. This will happen within a finite time (bounded by $\sum_{j=0}^N(L_j+I_j)$). The final number of infected individuals excluding the index case is as before denoted $Z^N\in [0,\dots ,N]$. Similar to before we let $L^N(t), I^N(t), R^N(t)$ denote the numbers of latent, infectious and recovered individuals at time $t$, and now we can also define the number of susceptibles $S^N(t)=N+1- L^N(t) - I^N(t) -R^N(t)$.
\end{definition}

In our model the index case cannot be contacted. This is of course unrealistic but simplifies notation. In the limit as $N$ gets large this assumption has no effect. We now state two important results for these constructions of the branching process and epidemics.

\begin{theorem}
The definition above agrees with the earlier definition of the Stochastic SEIR epidemic in a homogeneous community.
\end{theorem}

\bpf
The latent and infectious periods have the desired distributions, and an infective has infectious contacts with others at overall rate $\lambda$, and each time such a contact is with a uniformly selected individual as desired.
\epf

We now prove that the branching process and the epidemic process (with population size $N$) are identical up to a time point which tends to infinity in probability as $N\to\infty$. To this end, we let $M^N$ denote the number of infections prior to the first ghost (i.e.\ how many uniformly selected individuals $[U_kN]$ there were before someone was reselected. If this never happens we set $M^N=\infty$. Let $T^N$ denote the time at which the first ghost appears (and if this never happens we also set $T^N=\infty$).

\begin{theorem}\label{Th-br-epid-coupling}
The branching process and $N$-epidemic agree up until $T^N$:\linebreak $(L^N(t),I^N(t),R^N(t))=(L(t),I(t),R(t))$ for all $t\in [0,T^N )$. Secondly, $T_N\to\infty $ and $M^N\to\infty$ in probability  as $N\to\infty$.
\end{theorem}

\bpf
The first statement of the proof is obvious. The only difference between the epidemic and the branching process in our construction is that specific individuals are contacted in the epidemic, and up until the first time when some individual is contacted again, each infectious contact results in infection just as in the branching process.

As for the second part of the theorem we first compute the probability that $M^N$ will tend to infinity, and then that the time $T^N$ until the first ghost appears also tends to infinity. It is easy to compute $\P (M^N>k)$ since this will happen if and only if all the first $k$ contacts are with distinct individuals:
\begin{align*}
\P (M^N>k)=1\times \frac{N-1}{N}\times \dots \times \frac{N-k}{N}=\prod_{j=0}^k\left( 1-\frac{j}{N}\right) .
\end{align*}
(This formula is identical to the celebrated (...) birthday problem if $N+1=365$ and $k$ is the size of the class.) For fixed $k$ we see that this probability tends to 1 as $N\to\infty$. We can in fact say more. We have the following lower bound (which is easily proved by recurrence):
\begin{align*}
\P (M^N>k)=\prod_{j=0}^k\left( 1-\frac{j}{N}\right) \ge1-\sum_{j=1}^k\frac{j}{N}  =1-\frac{(k+1)k}{2N}.
\end{align*}
As a consequence, we see that $\P (M^N>k(N))\to 1$ as long as $k=k(N)=o(\sqrt{N})$. In particular $M^N\to\infty$ in probability as $N\to\infty$.  In what follows we write w.l.p.\ for ``with large probability'', meaning with a probability tending to 1 as $N\to\infty$. The consequence hence implies that all infectious contacts up to $k(N)$ will w.l.p.\ be with distinct individuals and thus will result in infections. So, up until $k(N)$ individuals have been infected, the epidemic can be approximated by a branching process for any $k(N)=o(\sqrt{N})$. Let $Z(t)$ denote the number of individuals born before $t$ in the branching process (excluding the ancestor) and $Z^N(t)=N-S^N(t)$ the number of individuals that have been infected before $t$ (excluding the index case) in the $N$-epidemic. Since the epidemic and branching process agree up until $T^N$ it follows that $Z(t)=Z^N(t)$ for $t<T^N$. But, since $k(N)<M^N$ w.l.p.\ it follows that $\inf \{t;Z(t)=k(N)\}\le T^N$ w.l.p. If the branching process is (sub)critical, then $Z(t)$ remains bounded as $t\to\infty$, so $T^N=+\infty$ w.l.p. Consider now the supercritical case. From Section \ref{Sec_Cont_Br_Pr} (Proposition \ref{contBGW_large_t}) we know that $Z(t)=O_p (e^{rt})$ where the Malthusian parameter $r$ solves the equation
\begin{equation}\label{Malthus-eq}
\int_0^\infty e^{-rs} \lambda (s)ds=1.
\end{equation}
The function $\lambda (s)$ is the rate at which an individual gives birth $s$ time units after being born, so $\lambda (s)=\lambda \P(\text{infectious at $s$})$ and hence $\lambda (s)=\lambda \P (L<s<L+I)$ for our model. We thus have that $k(N)\le ce^{rT^N}$ w.l.p., which implies that $T^N\ge \log k(N)/r - \log c$. So if for example $k(N)=N^{1/3}$, which clearly satisfies $k(N)=o(\sqrt{N})$, it follows that $T^N\to\infty$ in probability.
\epf

Theorem \ref{Th-br-epid-coupling} shows that the epidemic behaves like the branching process up to a time point tending to infinity as $N\to\infty$, and that the number of infections/births by then also tends to infinity. This implies that we can use theory for branching processes to obtain results for the early part of the epidemic. We state these important results in the following corollaries; the first corollary is for the subcritical and critical cases and the second corollary is for the supercritical case. Recall that $R_0=\lambda \E (I)$, the basic reproduction number in the epidemic and the mean offspring number in the branching process.

\begin{corollary}\label{Th-early-stage-sub}
If $R_0\le 1$, then $(L^N(t),I^N(t),R^N(t))=(L(t),I(t),R(t))$ for all $t\in [0,\infty )$ w.l.p. As a consequence, $\P( Z^N=k)\to \P (Z=k)$ as $N\to\infty$, and in particular $Z^N$ is bounded in probability.
\end{corollary}

\bpf
In Theorem \ref{Th-br-epid-coupling} it was shown that the epidemic and branching process agree up until there has been $M^N$ births, where $M^N>N^{1/3}$ w.l.p.\   for example. But from branching process theory (Proposition \ref{BGW}) we know that this will happen with a probability tending to 0 with $N$ when $R_0\le 1$, implying that $T^N=\infty$ w.l.p.
\epf

\begin{corollary}\label{Th-early-stage-sup}
If $R_0> 1$, then for finite $k$: $\P (Z^N=k)\to \P (Z=k)$ as $N\to\infty$. Further, $\{Z^N\to\infty\}$ with the same probability as $\{ Z=\infty\}$, 
which is the complement to the extinction probability, the latter being the smallest solution to the equation $z=g(z)$ described in Proposition \ref{BGW}.
\end{corollary}
\bpf
Also this corollary is a direct consequence of Theorem \ref{Th-br-epid-coupling} and properties of branching processes. If only $k$ births occur, then there will be no ghost w.l.p., implying that the epidemic and the branching process agree forever w.l.p. On the other hand, the coupling construction showed that $M^N\to\infty$ on the other part of the sample space, and $Z\ge Z^N\ge M^N$ which completes the proof.
\epf

The two corollaries state that the epidemic and branching process coincide forever as long as the branching process stays finite. If the branching process grows beyond all limits (only possible when $R_0>1$) then the epidemic and branching process will not remain identical even though also the epidemic tends to infinity with $N$. For any fixed $N$ we have $0\le Z^N\le N$ which clearly is different from $Z=\infty$ in that case. The distribution of $Z^N$ on the part of the sample space where $Z^N\to\infty$ is treated below in Section \ref{TB-EP_sec_LLNCLT_final_size}.

The two corollaries show that the final number infected $Z^N$ will be small with a probability equal to the extinction probability of the approximating branching process, and it will tend to infinity with the remaining (explosion) probability. In Section \ref{TB-EP_sec_LLNCLT_final_size} we study the distribution of $Z^N$ (properly normed) and then see that the distribution is clearly bimodal with one part close to 0 and the other part being $O(N)$. These two parts are referred to as \emph{minor outbreak} and \emph{major outbreak} respectively.

What happens during the early stage of an outbreak is particularly important when considering so-called \emph{emerging epidemic outbreaks}. Then statistical inference based on this type of branching process approximation is often used. For example, in \cite{WHO2014I} a branching process approximation that is very similar to the SEIR branching process of Definition \ref{DefApprBrPr} is used for modelling the spread of Ebola during the early stage of the outbreak in West Africa in 2014.

\begin{exercise}\label{xEarly.2}
Use the branching process approximation of the current section to compute the probability of a major outbreak of the SEIR epidemic assuming that $I\equiv \iota$ (the continuous time Reed--Frost case), and $I\sim \mathrm{Exp}(\gamma)$ (the Markovian SIR) with $\gamma=1/\iota$. Only one of them will be explicit. Compute things numerically for $R_0=1.5$ and $\iota=1$.
\end{exercise}

\begin{exercise}\label{xEarly.3}
Use the branching process approximation of the current section to compute the exponential growth rate $r$ for the following two cases: $L\equiv 0$ and $I\equiv \iota$ (the continuous time Reed--Frost), and $L\equiv 0$ and $I\sim \mathrm{Exp}(\gamma=1/\iota)$ (the Markovian SIR). Compute $r$ numerically for the two cases when $R_0=1.5$ and $\iota =\gamma=1$.
\end{exercise}

\section[The final size of the epidemic in case of no major outbreak]{The final size of the epidemic in case of no major outbreak}\label{FinalSizeMinor}

Let $Z^N$ denote the final size of the epidemic (i.e.\ the total number of individuals that get infected during the outbreak) but now also including the initially infected individual.
In the case of no major outbreak, if the total population size $N$ is large enough, $Z^N$ is well approximated by the total number of individuals in a branching process, as we saw in the previous section.
Hence we consider $Z$ as the total number of individuals ever born in a branching process (including the ancestor), where the number of offspring of the $k$-th individual is $X_k$. Let $X_1,X_2,\ldots$ be i.i.d.\ $\N$-valued random variables.
We start by establishing
an identity which is an instance of Kemperman's formula, see e.g.\ Pitman \cite{PitI} page 123.

\begin{figure}[!ht]
\begin{center}
\begin{tabular}{c}
\includegraphics[height=6cm,width=9cm]{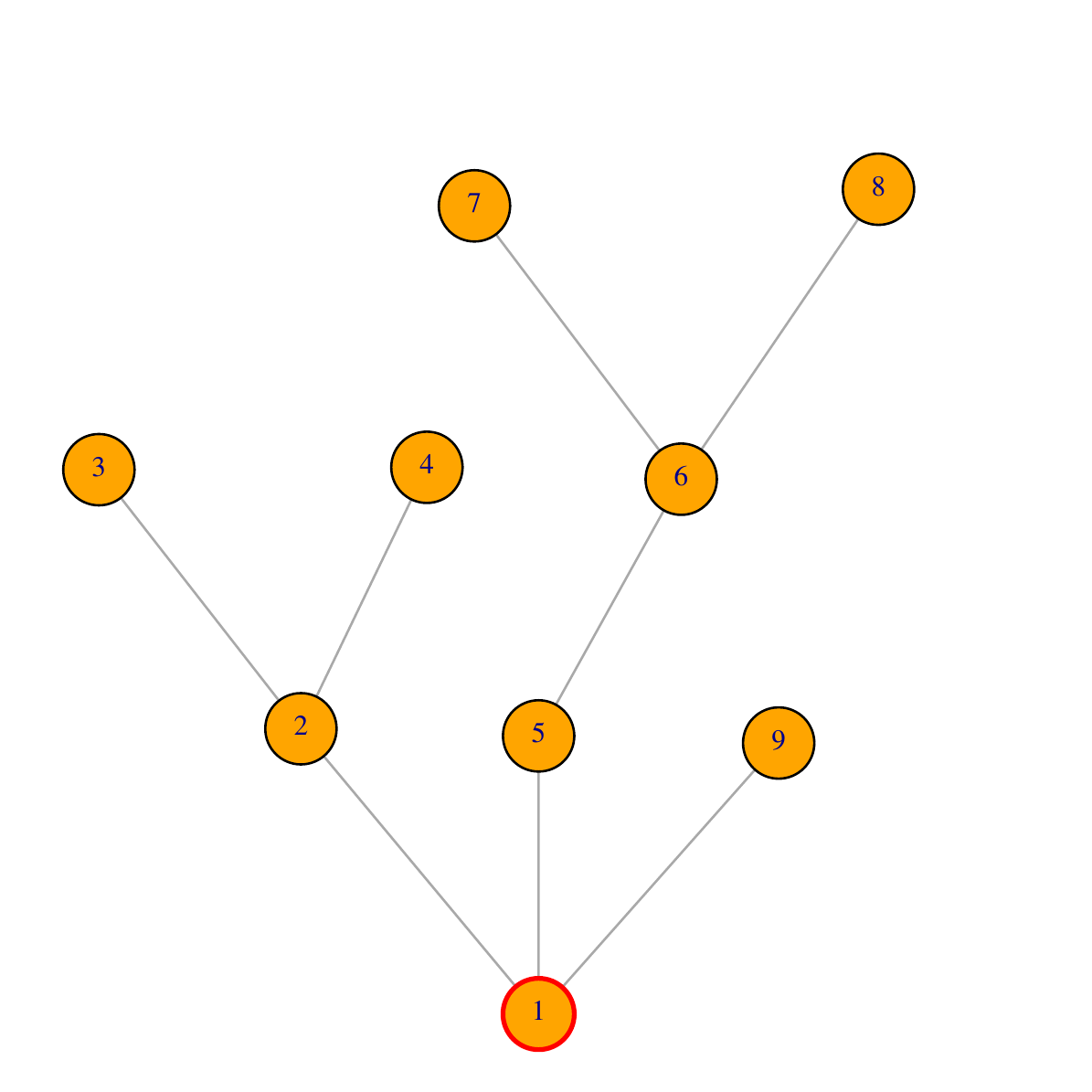}\\
\includegraphics[height=6cm,width=9cm]{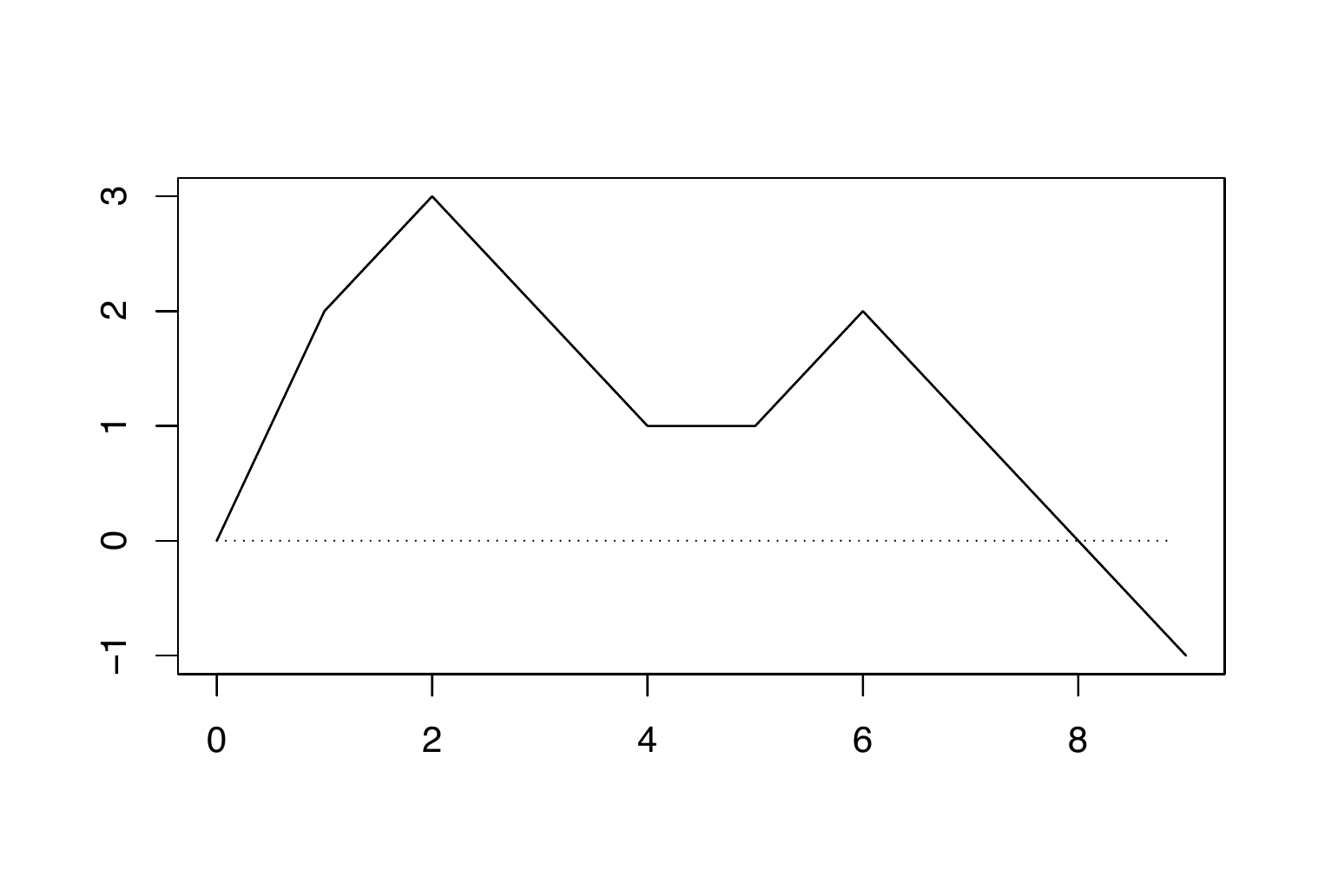}
\end{tabular}
\caption{Top: the tree. Bottom: the random walk $S_k$.
Here $X_1=3$, $X_2=2$, $X_3=0$, $X_4=0$, $X_5=1$, $X_6=2$, $X_7=X_8=X_9=0$,
$Y_1=2$, $Y_2=1$, $Y_3=Y_4=-1$, $Y_5=0$, $Y_6=1$, $Y_7=Y_8=Y_9=-1$.}
\label{fig:dfs}
\end{center}
\end{figure}

\begin{proposition}\label{pro:Kemperman}
For all $k\ge1$,
\[ \P(Z=k)=\frac{1}{k}\P(X_1+X_2+\cdots+X_k=k-1).\]
\end{proposition}
\bpf  Consider the process of depth--first search of the genealogical tree of the infected individuals.
This procedure can be defined as follows.
The tree is explored starting from the root. Suppose we have visited k vertices. The next visit will be to the
leftmost still unexplored son of this individual, if any; otherwise to the leftmost unexplored son of the most recently visited node among those having not yet visited son(s), see Figure \ref{fig:dfs}. $X_1$ is the number of sons of the root, who is the first visited individual.
$X_k$ is the number of sons of the $k$-th visited individual.
This exploration of the tree ends at step $k$ if and only if
$X_1\ge1$, $X_1+X_2\ge2$, $X_1+X_2+X_3\ge3$, ... $X_1+X_2+\cdots X_{k-1}\ge k-1$, and
$X_1+X_2+\cdots+X_k=k-1$. Let us rewrite those conditions.
Define
\begin{align*}
Y_i&=X_i-1, \ i\ge1,\\
S_k&=Y_1+Y_2+\cdots+Y_k.
\end{align*}

\noindent A trajectory $\{Y_i,\ 1\le i\le k\}$ explores a tree of size $k$ if and only if the following conditions are satisfied
\[(C_k)\quad S_0=0, S_1\ge0, S_2\ge0,\ldots, S_{k-1}\ge0, S_k=-1.\]
Indeed, it is easy to convince oneself that it is the case if there is only one generation: if the ancestor has $k-1$ children, then $Y_1=k-2$,
and $Y_2=\cdots=Y_{k}=-1$, hence $(C_k)$ holds. If one attaches one generation trees to some of the leaves of the previous tree, then one replaces
a unique $-1$ step by an excursion upwards which finishes at the same level as the replaced step. Iterating this procedure, we see that the exploration of a general tree with $k$ nodes satisfies $(C_k)$.

The statement of the proposition is equivalent to
\[ \P(Z=k)=\frac{1}{k}\P(Y_1+Y_2+\cdots+Y_k=-1).\]
Denote by $V_k$ the set of sequences of $k$ integers $\ge-1$ which satisfy conditions $(C_k)$, and
$U_k$ the set of sequences of $k$ integers $\ge-1$ which satisfy the unique condition $S_k=-1$.
We use circular permutations operating on the $Y_i$'s. For $1\le i, \ell\le k$, let
\[  (i+\ell)_k=\begin{cases} i+\ell,&\text{if $i+\ell\le k$};\\
i+\ell-k,&\text{if $i+\ell>k$.}\end{cases} \]
For each $1\le\ell\le k$, let $Z^\ell_i=Y_{(i+\ell)_k}$, $S^\ell_j=\sum_{i=1}^j Z^\ell_i$ for $1\le i\le k$.
Clearly $S^\ell_k=-1$ for all $\ell$ as soon as $(C_k)$ is satisfied. On the other hand $S^k\equiv S$
is the only trajectory which satisfies conditions $(C_k)$. The other $S^\ell$ hit the value $-1$ before rank $k$, see Figure \ref{fig:dfs}.
The $Z^\ell$'s are sequences of integers
 $\ge-1$ of length $k$, whose sum equals $-1$.  Finally to each element of $V_k$ we have associated
$k$ distinct elements of $U_k$, all having the same probability.

Reciprocally, to one element $S$ of $U_k\backslash V_k$, choosing $\ell=\underset{1\le i\le k}{\operatorname{argmin}}S_i$
and using the above transformation, we deduce that $S^\ell\in V_k$.

Finally, to each trajectory of $V_k$, we associate $k$ trajectories of $U_k$, who all have the same
 probability, and which are such that the inverse transformation gives back the same trajectory of $V_k$. The result is proved.  \epf

Note that from branching process theory (Proposition \ref{BGW}), we have clearly
\[\sum_{k\ge1}\P(Z=k)\begin{cases} = 1,&\text{ if $\E R_0\le1$};\\
< 1,&\text{ if $\E R_0>1$},\end{cases}\]
which is not so obvious from the proposition.

We now deduce the exact law of $Z$ from Proposition \ref{pro:Kemperman} in two cases which are probably the two most interesting cases for epidemics models. First we consider the case where the $X_i$s are Poisson, which is the situation of the continuous time Reed--Frost model, where the infectious period is non-random. Second we consider the case where the $X_i$s are geometric, which is the case in the Markovian model.

\begin{example}
Suppose that the joint law of the $X_i$s is Poi$(\mu)$, with $0<\mu<1$.
Then $X_1+\cdots+X_k\sim\text{Poi}(k\mu)$, and consequently
\begin{align*}
\P(Z=k)&=\frac{1}{k}\P(X_1+\cdots+X_k=k-1)\\
&=e^{-\mu k}\frac{(\mu k)^{k-1}}{k!}.
\end{align*}
This law of $Z$ is called the Borel distribution with parameter $\mu$.
Note that
\begin{align*}
\E Z&=1+\mu+\mu^2+\cdots\\
&=\frac{1}{1-\mu}.
\end{align*}
\end{example}

\begin{example}
Consider now the case where $X_i\sim\mathcal{G}(p)$, where we mean here that $\P(X_i=k)=(1-p)^k p$, $k=0,1,\ldots$. The law of
$X_i+1$ is the geometric distribution with parameter $p$ whose support is $\N$,
in other words $\P(X_i+1>k)=(1-p)^k$. Then $k+X_1+\cdots+X_k$ follows the negative binomial distribution
with parameters $(k,p)$. Hence
\begin{align*}
\P(Z=k)&=\frac{1}{k}\P(k+X_1+\cdots+X_k=2k-1)\\
&=\frac{1}{k}\begin{pmatrix}2k-2\\ k-1\end{pmatrix}p^k(1-p)^{k-1}\\
&=\frac{(2k-2)!}{k! (k-1)!}p^k(1-p)^{k-1}.
\end{align*}
In the case $p>1/2$, $\E Z=(2p-1)^{-1}p$.
\end{example}

\section{Vaccination}\label{Sec_Vacc}
One important reason for modelling the spread of infectious diseases is to better understand effects of different preventive measures, such as for example vaccination, isolation and school closure. When a new outbreak occurs, epidemiologists (together with mathematicians and statisticians) estimate model parameters and then use these to predict effects of various preventive measures, and based on these predictions, health authorities decide upon which preventive measures to put in place, cf.\ \cite{WHO2014I}.

We refer the reader to Part IV in this volume for estimation methods, but in the current section we touch upon the area of modeling prevention. Our focus is on vaccination, and we consider only vaccination prior to the arrival of an outbreak; the situation where vaccination (or other preventive measures) are put into place \emph{during} the outbreak is not considered. ``Vaccination'' can be interpreted in a wider sense. From a mathematical and spreading point of view, the important feature is that the individual cannot spread the disease further, which could also be achieved by e.g.\ isolation or medication. Modelling effects of vaccination is also considered in Part II, Section 2.4, and in Part III, Section 2.6, in the current volume.

Suppose that a fraction $v$ of the community is vaccinated prior to the arrival of the disease. We assume that the vaccine is perfect in the sense that it gives 100\% protection from being infected and hence of spreading the disease (but see the exercise below). This implies that only a fraction $1-v$ are initially susceptible, and the remaining fraction $v$ are immunized (as discussed briefly in Section \ref{Sec_Det-mod}). Hence we can neglect the latter fraction and consider only the initial susceptible part of the community of size $N'=N(1-v)$. However, it is not only the number of initially susceptibles that changes, the rate of having contact with initial susceptibles has also changed to $\lambda'=\lambda (1-v)$, since a fraction $v$ of all contacts are ``wasted'' on vaccinated people. The spread of disease in a partly-vaccinated community can therefore be modelled using exactly the same SEIR stochastic model with the only difference being that we have a different population size $N'$ and a different contact rate parameter $ \lambda '$.

From this we conclude the new reproduction number, which we denote $R_v$ to show the dependence on $v$, satisfies
$$
R_v=\lambda' \E (I)=\lambda (1-v)\E (I)=(1-v)R_0.
$$
As a consequence, a major outbreak in the community is not possible if $R_v\le 1$, which (when $R_0>1$) is equivalent to $v\ge 1-1/R_0$. This limit, called the \emph{critical vaccination coverage} and denoted
\begin{equation}\label{crit-vacc}
 v_c=1-\frac{1}{R_0},
\end{equation}
is hence a very important quantity: if more than this fraction is vaccinated before an outbreak, then the whole community is protected from a major outbreak and not only the vaccinated, a situation called \emph{herd immunity}. Equation (\ref{crit-vacc}) is well known among infectious disease epidemiologists (e.g.\ Giesecke \cite{G2017I}) and is used by public health authorities all over the world to determine the minimal yearly vaccination coverage in vaccination programs of childhood diseases.

If $v<v_c$ there is still a possibility of a major outbreak. The probability for such an outbreak is obtained using earlier results with $\lambda$ replaced by $\lambda'=\lambda (1-v)$: the probability of a minor outbreak is the solution $s_v$ to the equation $s=g_v(s)$, where $g_v(\cdot)$ is the probability generating function of $X_v\sim \text{MixPoi}(\lambda(1-v)I)$, the number of offspring (= new infections) in the case that a fraction $v$ are immunized by vaccination.

In the case when there is a major outbreak, the relative size $z_v$ of the outbreak (among the initially susceptible!) is given by the unique positive solution to the equation
\begin{equation}
1-z=e^{-R_vz}\text{, or equivalently }1-z=e^{-(1-v)R_0z},
\end{equation}
this result is shown in later sections, cf.\  Equation (\ref{final-size-vacc}).
The community fraction getting infected is hence $(1-v)z_v$.

We summarize our result in the following theorem where we let $Z^N_v$ denote the final number infected when a fraction $v$ are vaccinated prior to the outbreak.
\begin{theorem}\label{Th-Vacc}
If $v\ge v_c=1-1/R_0$, then $Z^N_v/N\to 0$ in probability. If $v< v_c=1-1/R_0$, then $Z^N_v/N\Rightarrow Z^\infty_v$ which has a two-point distribution: $\P (Z^\infty_v =0)=s_v$ and $\P (Z^\infty_v = (1-v)z_v)=1-s_v$, where $s_v$ and $z_v$ have been defined above.
\end{theorem}

\begin{exercise} \label{xvacc-numeric}
Consider the Markovian SEIR epidemic with $\lambda=2$, $L\sim \mathrm{Exp}(2)$ and $I\sim \mathrm{Exp}(1)$. Compute the critical vaccination coverage $v_c$. Compute also numerically the probability of a major outbreak, and the community-fraction that will get infected in the case of a major outbreak when $v=0.333$.
\end{exercise}

\begin{exercise} \label{xnon-perf-vacc}
Suppose that the vaccine gives only partial protection to catching and spreading the disease. Suppose that the vaccine has the effect the risk of getting infected by a contact is only 20\% of the risk of getting infected when not vaccinated, but that the vaccine has no effect on infectivity if the person gets infected (such a vaccine is said to be a ``leaky vaccine'' having 80\% efficacy on susceptibility and 0\% efficacy on infectivity). Compute the reproduction number $R_v$ in the case that a fraction $v$ is vaccinated with such a vaccine. (Another vaccine response model is ``all-or-nothing'' where a fraction is assumed to receive 100\% effect and the remaining fraction receive no effect from vaccination, for example due to the cold chain being broken for a live vaccine.)
\end{exercise}

\chapter{Markov Models}\label{chap_MarkovMod}

This chapter describes the important class of Markov models. It starts with a presentation of the deterministic ODE models. We then formulate precisely the random Markov epidemic model as a Poisson process driven stochastic differential equation, and establish the law of large numbers (later referred to as LLN), whose limit is precisely the already described ODE model. The next section studies the fluctuations around this LLN limit, which is described by the central limit theorem. Finally we give a diffusion approximation result, i.e.\ a diffusion process (solution of a Brownian motion driven stochastic differential equation) which, again in the case of a large population, is a good approximation of our Poisson process driven model.
One of the earliest references for those three approximation theorems is Kurtz \cite{Ku78I}. See also chapter 11 of Ethier and Kurtz \cite{EKI}.

\section{The deterministic SEIR epidemic model}\label{Sec_Det-mod}
Before analysing the stochastic SEIR model assuming $N\to\infty$ in greater detail in the following subsections, we first derive heuristically  a deterministic counterpart for the Markovian version and study some of its properties, which are relevant also for the asymptotic case of the stochastic model.

Consider the Markovian stochastic SEIR model. There are three types of events: a susceptible gets infected and becomes exposed, an exposed becomes infectious when the latent period terminates, and an infectious individual recovers and becomes immune. Since the model is Markovian all these events happen at rates depending only on the current state, and these rates are respectively given by: $\lambda S(t)I(t)/N$, $\nu E(t)$ and $\gamma I(t)$. When an infection occurs, the number of susceptibles decreases by 1 and the number of exposed increases by 1; when a latency period ends, the number of exposed decreases by 1 and the number of infectives increases by 1; and finally when there is a recovery, the number of infectives decreases by 1 and the number of recovered increases by 1. If we instead look at ``proportions'' (to simplify notation we divide by $N$ rather than the more appropriate choice $N+1$), the corresponding changes are $-1/N$ and $+1/N$. This reasoning justifies a deterministic model for proportions where one should think of an infinite population size allowing the proportions to be continuous. The deterministic SEIR epidemic $(s(t),e(t),i(t),r(t))$ is given by
\begin{align*}
s'(t)&=-\lambda s(t)i(t),
\\
e'(t)&= \lambda s(t)i(t) -\nu e(t),
\\
i'(t)&=\nu e(t)-\gamma i(t),
\\
r'(t)&=\gamma i(t).
\end{align*}
We start with all fractions being non-negative and summing to unity, which implies that $s(t)+e(t)+i(t)+r(t)=1$ and all being nonnegative for all $t$.
It is important to stress that this system of differential equations only approximates the \emph{Markovian} SEIR model. If for example the latent and infectious stages are non-random, then a set of differential-delay equations would be the appropriate approximation. If these durations are random but not exponential one possible pragmatic assumption is to use a gamma distribution where the shape parameter is an integer (so it can be seen as a sum of i.i.d.\ exponentials). Then the deterministic approximation would be a set of differential equations where the state space has been expanded. Just like for the stochastic SEIR model, the deterministic model has to start with a positive fraction of exposed and/or infectives for anything to happen. Most often it is assumed that there is a very small fraction $\epsilon$ of latent and/or infectives.

The case where there is no latent period meaning that $\nu\to\infty$, the deterministic SIR epidemic (or deterministic general epidemic), sometimes called the Kermack--McKendrick equations, has perhaps received more attention in the literature:
\begin{equation}\label{SIRdet}
\left\{
\begin{aligned}
s'(t)&=-\lambda s(t)i(t),
\\
i'(t)&= \lambda s(t)i(t) -\gamma i(t),
\\
r'(t)&=\gamma i(t).
\end{aligned}
\right.
\end{equation}
This system of differential equations (and the SEIR system on the previous page) are undoubtedly the most commonly analysed epidemic models (e.g.\ Anderson and May \cite{AM91I}), and numerous related extended models, capturing various heterogeneous aspects of disease spreading, are published every year in mathematical biology journals.

\begin{figure*}[h]
\begin{center}
\includegraphics[width=0.95\textwidth, height=0.35\textheight]{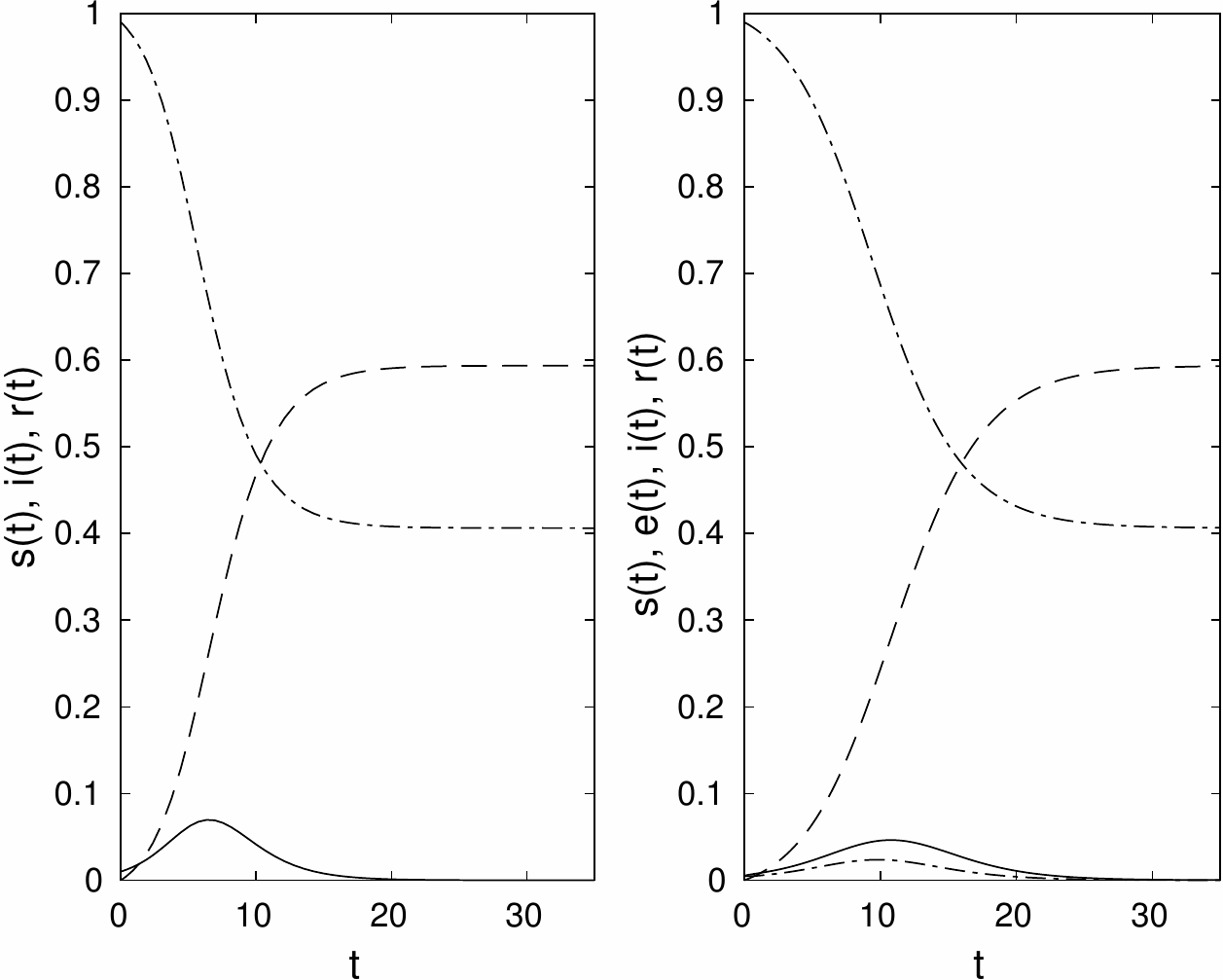}
\end{center}
\caption{Plot of the deterministic SIR (left) and SEIR (right) systems for $\lambda=1.5$ and $\gamma=1$, and with $\nu=1$ in the SEIR model. The dash-dotted curve is the fraction of susceptibles, the solid curve the fraction of infectives, the dashed curve the fraction of recovered, and the lowest curve in the right figure is the fraction of exposed (latent).}
\label{fig-det-syst}
\end{figure*}

The deterministic SEIR and SIR share the two most important properties in that they have the same basic reproduction number $R_0$ and give the same final size (assuming the initial number of infectives/exposed are positive but negligible in both cases), which we now show. In Figure \ref{fig-det-syst} both the SEIR and SIR systems are plotted for the same values of $\lambda=1.5$ and $\gamma=1$ (so $R_0=1.5)$, and with $\nu=1$ in the SEIR system.

From the differential equations we see that $s(t)$ is monotonically decreasing and $r(t)$ monotonically increasing. The differential for $i(t)$ in the SIR model can be written $i'(t)=\gamma i(t)\left( \frac{\lambda}{\gamma}s(t)-1\right)$. The initial value is $i(0)=\epsilon\approx 0$ and $s(0)=1-\epsilon\approx 1$. From this we see that for having $i'(0)>0$ we need that $\lambda/\gamma >1$. If this holds, $i(t)$ grows up until $s(t)<\gamma/\lambda$ after which $i(t)$ decays down to 0. If on the other hand $\lambda/\gamma \le 1$, then $i(t)$ is decreasing from the start and since its initial value is $\epsilon\approx 0$, nothing much will happen so $s(\infty)\approx s(0)\approx 1$ and $r(\infty )\approx r(0)=0$. We hence see that also in the deterministic model, $R_0=\lambda/\gamma$ plays an important role in that whether or not $R_0$ exceeds 1 determines whether there will be a substantial or a negligible fraction getting infected during the outbreak. Note that this is the same $R_0$ as for the Markovian SEIR epidemic. There the infectious period is exponentially distributed with parameter $\gamma$, so $\iota :=E(I)=1/\gamma$.

An important difference between deterministic and stochastic epidemic models lies in the initial values. Stochastic models usually start with a small \emph{number} of infectious individuals (in the model of the current Notes we assumed one initial infective: $I(0)=1$). This implies that the initial \emph{fraction} of infectives tend to 0 as $N\to\infty$. In the deterministic setting we however have to assume a fixed and strictly positive fraction $\epsilon$ of initially infectives (if we start with a fraction 0 of infectives nothing happens in the deterministic model). This implicitly implies that the deterministic model starts to approximate the stochastic counterpart only when the number of infectives in the stochastic model has grown up to a \emph{fraction} $\epsilon$, so a number $N\epsilon$. The earlier part of the stochastic model cannot be approximated by this deterministic model, and as we have seen it might in fact never reach this level (if there is only a minor outbreak).

In order to derive an expression for the ultimate fraction getting infected we use the differential for $s(t)$ (and below also the one for $r(t)$). Dividing by $s$ and multiplying by $dt$ gives the following differential: $ds/s=-\lambda i dt$. Integrating both sides
and recalling that $R_0=\lambda/\gamma$, we obtain
\begin{align*}
\log s(t)-\log s(0)&=-\lambda \int_0^ti(t)dt\\ &=-R_0\int_0^t r'(s)ds\\&=-R_0(r(t)-r(0))= -R_0r(t).
\end{align*}
And since $s(0)=1-\epsilon\approx 1$ and $r(\infty)=1-s(\infty)$ we obtain the following equation for the final size $z=r(\infty)=1-s(\infty)$:
\begin{equation}\label{final-size}
1-z=e^{-R_0z} .
\end{equation}

In Section \ref{subsec_LLN} we show that this final size equation coincides with that of the LLN limit of the final fraction getting infected in the stochastic model (cf.\ Equation (\ref{eq:LLNfinal-size}), which is identical to \eqref{final-size}).

\begin{figure*}[h]
\begin{center}
\includegraphics[width=.8\textwidth, height=0.35\textheight]{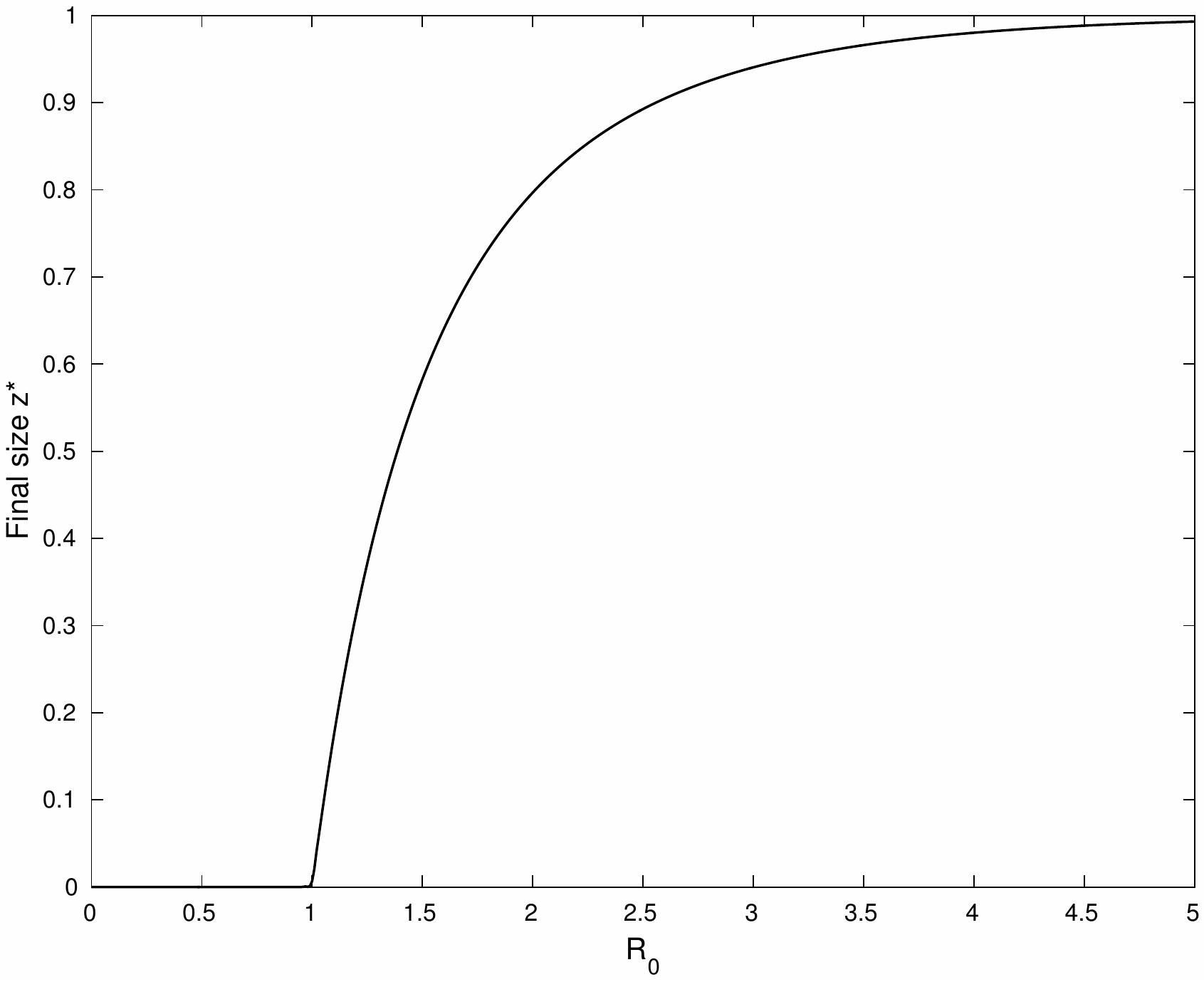}\label{fig-final-size}
\end{center}
\caption{Plot of the final size solution $z^*$ to Equation (\ref{final-size}) as a function of $R_0$. }
\label{fig-R0}
\end{figure*}

The equation always has a root at $z=0$ corresponding to no (or minor) outbreak. It can be shown (cf.\ Exercise \ref{x4.1}) that if and only if $R_0>1$ there is a second solution to (\ref{final-size}), corresponding to the size of a major outbreak, and this solution $z^*$ is strictly positive and smaller than 1. For a given value of $R_0>1$ the solution $z^*$ has to be computed numerically. In Figure \ref{fig-R0} the solution is plotted as a function of $R_0$.

It is important to point out that the final size equation (\ref{final-size}) assumes that, at $t=0$, all individuals (except the very few initially latent and infectives) are susceptible. If a fraction $v$ is initially immune (perhaps due to natural immunity, or vaccination as described in Section \ref{Sec_Vacc}) then $r(0)=v$ and $s(0)=1-v$, resulting in the equation
\begin{equation}
1-z=e^{-R_0z(1-v)}, \label{final-size-vacc}
\end{equation}
where its solution $z_v$ now is interpreted as the fraction among the initially susceptible that get infected. The overall fraction getting infected is hence $z_v(1-v)$. Using the same argument as for the final size without immunity, we conclude that $z=0$ is the only solution if $R_0(1-v)\le 1$. This is equivalent to $v\ge 1-1/R_0$. If immunity was caused by vaccination, this hence suggests that a fraction exceeding $v_c=1-1/R_0$ should be vaccinated; then there will be no outbreak! For this reason, the quantity $v_c=1-1/R_0$ is often called the \emph{critical vaccination coverage}, and if this coverage is reached, so-called \emph{herd immunity} is achieved. Herd immunity implies that not only the vaccinated are protected, but so are also the unvaccinated, since the community is protected from epidemic outbreaks.

\begin{exercise} \label{x4.1}
Show that $z=0$ is the only solution to (\ref{final-size}) when $R_0\le 1$ and that there is a unique positive solution if $R_0>1$. (\emph{Hint}: Study suitable properties of the function $f(z)= e^{-R_0z} +z-1$.)
\end{exercise}

\begin{exercise} \label{x4.2}
Compute the final size numerically for $R_0=1.5$ (e.g.\ influenza), $R_0=3$ (e.g.\ rubella) and $R_0=15$ (e.g.\ measles).
\end{exercise}

\section{Law of Large Numbers}\label{TB-EP_sec_LLN}

Consider a general compartmental model, which takes the form
\[ \mathcal{Z}^N_t=z_N+\sum_{j=1}^kh_jP_j\left(\int_0^t\beta_{N,j}(s,\mathcal{Z}^N_s)ds\right),\]
where the $P_j$s are mutually independent standard (i.e.\  unit rate) Poisson processes, and
$\beta_{N,j}(t,\mathcal{Z}^N_t)$ is the rate of jumps in the direction $h_j$ at time $t$, $h_j$ being a $d$-dimensional vector. $\mathcal{Z}^N_t$ takes values in
$\Z_+^d$. The $i$-th component of $\mathcal{Z}^N_t$ is the number of individuals in the $i$-th compartment
at time $t$. $N$ is a scale parameter. In the case of models with fixed total population size, $N=\sum_{i=1}^dZ^{N,i}_t$
is the total population at any time $t$.
Note that the above formula for $ \mathcal{Z}^N_t$ can be rewritten equivalently, following the comments at the end of Section \ref{sec-Po-proc} in the Appendix below, as
\[ \mathcal{Z}^N_t=z_N+\sum_{j=1}^kh_j\int_0^t\int_0^{\beta_{N,j}(s,\mathcal{Z}^N_s)}Q_j(ds,du),\]
where $Q_1,\ldots,Q_k$ are mutually independent Poisson random measures on $\R^2_+$, with mean measure
$ds\, du$.

We now define
\[ Z^N_t=N^{-1}\mathcal{Z}^N_t\]
the vector of rescaled numbers of individuals in the various compartments. In the case of a constant population size equal to $N$, the components of the vector $Z^N_t$ are the proportions of the total population in the various compartments at time $t$. The equation for $Z^N_t$ reads, with $x_N=N^{-1}z_N$,
\begin{equation*}
 Z^N_t=x_N+\sum_{j=1}^k\frac{h_j}{N}P_j\left(\int_0^t\beta_{N,j}(s,NZ^N_s)ds\right).
 \end{equation*}

\begin{example}\label{example1}{The SIR model.}

One important example is that of the $SIR$ model with constant population size.
Suppose there is no latency period  and that the duration of infection satisfies $I\sim\mathrm{Exp}(\gamma)$.
 In that case, let $S(t)$, $I(t)$, and $R(t)$) denote respectively the number of susceptibles, infectives and recovered at time $t$.

%
In this model, two types of events happen:
\begin{enumerate}
\item infection of a susceptible (such an event decreases $S(t)$ by one, and increases $I(t)$ by one, so $h_1=(-1,1,0)$);
these events happen at rate
\[\beta_{N,1}(t,\mathcal{Z}_t)= \frac{\lambda}{N}S(t)I(t),\quad\text{where }\lambda=cp;\]
\item recovery of an infective (such an event decreases $I(t)$ by one, and increases $R(t)$ by one, so $h_2=(0,-1,1)$); these
events happen at rate
\[\beta_{N,2}(t,\mathcal{Z}_t)=\gamma I(t).\]
\end{enumerate}
Hence we have the following equations, with $P_1(t)$ and $P_2(t)$ two standard mutually independent Poisson
processes:
 \begin{align*}
S(t)&=S(0)-P_1\left( \frac{\lambda}{N}\int_0^t S(r)I(r)dr\right),\\
I(t)&=I(0)+P_1\left( \frac{\lambda}{N}\int_0^t S(r)I(r)dr\right)-P_2\left(\gamma\int_0^t I(r)dr\right),\\
R(t)&=R(0)+P_2\left(\gamma\int_0^t I(r)dr\right).
\end{align*}
We can clearly forget about the third equation, since
$R(t)=N-S(t)-I(t)$.

We now define $(S^N(t),I^N(t))=(N^{-1} S(t),N^{-1} I(t))$. We have
\begin{align*}
S^N(t)&=S^N(0)-\frac{1}{N}P_1\left(N\lambda\int_0^tS^N(r)I^N(r)dr\right),\\
I^N(t)&=I^N(0)+\frac{1}{N}P_1\left(N\lambda\int_0^t S^N(r)I^N(r)dr\right)-\frac{1}{N}P_2\left(N\gamma\int_0^tI^N(r)dr\right).
\end{align*}
\end{example}

The above model assumes that $\lambda$ and $\gamma$ are constant, but in applications at least $\lambda$ may depend upon $t$.

\begin{example}\label{ex:SEIRS}{The SEIRS model with demography.}

We now describe one rather general example. We add to the preceding example
 the state $E$ and the fact that removed individuals lose their immunity at a certain rate, which gives the SEIRS model. In addition, we add demography. There is an influx of susceptible individuals at
 rate $\mu N$, and each individual, irrespective of its type, dies at rate $\mu$.
 This gives the following stochastic differential equation

 \begin{align*}
 S(t)&=S(0)-P_{se}\left(\frac{\lambda}{N}\int_0^tS(r)I(r)dr\right)+P_{rs}\left(\rho\int_0^tR(r)dr\right)\\&\quad
 +P_b(\mu Nt)-P_{ds}\left(\mu\int_0^tS(r)dr\right),\\
 E(t)&=E(0)+P_{se}\left(\frac{\lambda}{N}\int_0^tS(r)I(r)dr\right)-P_{ei}\left(\nu\int_0^tE(r)dr\right)
\\&\quad -P_{de}\left(\mu\int_0^tE(r)dr\right),\\
 I(t)&=I(0)+P_{ei}\left(\nu\!\int_0^t\! E(r)dr\right)-P_{ir}\left(\gamma\!\int_0^t\! I(r)dr\right)
 -P_{di}\left(\mu\!\int_0^t\! I(r)dr\right),\\
 R(t)&=R(0)+P_{ir}\left(\gamma\!\int_0^t\! I(r)dr\right)-P_{rs}\left(\rho\!\int_0^t\!R(r)dr\right)
 -P_{dr}\left(\mu\!\int_0^t\! R(r)dr\right).
 \end{align*}
 In this system, the various Poisson processes are standard and mutually independent. The indices should be self-explanatory.
 Note that the rate of births is $\mu\times N$ rather than $\mu\times$ the actual number of individuals
 in the population, in order to avoid the pitfalls of branching processes (either exponential growth or extinction).
 Also, the probability $S(t)/N(t)$ that an infective meets a susceptible (where $N(t)$ denotes the total population at time $t$) is approximated by $S(t)/N$ for the sake of mathematical simplicity. Note however that $\frac{N(t)}{N}\to1$ a.s. as
 $N\to\infty$, see Exercise \ref{exer4.1.1} below. The equations for the proportions in the various compartments read

 \begin{align*}
 S^N(t)&=S^N(0)-\frac{1}{N}P_{se}\left(N\lambda\int_0^tS^N(r)I^N(r)dr\right)
 +\frac{1}{N}P_{rs}\left(N\rho\int_0^tR^N(r)dr\right)\\
 &\quad +\frac{1}{N}P_b(\mu Nt)-\frac{1}{N}P_{ds}\left(\mu N\int_0^tS^N(r)dr\right),\\
 E^N(t)&=E^N(0)+\frac{1}{N}P_{se}\left(N\lambda\int_0^t S^N(r)I^N(r)dr\right)\\&\quad
 -\frac{1}{N}P_{ei}\left(\nu N\int_0^tE^N(r)dr\right)
 -\frac{1}{N}P_{de}\left(\mu N\int_0^tE^N(r)dr\right),
 \end{align*}

\begin{align*}
I^N(t)&=I^N(0)+\frac{1}{N}P_{ei}\left(\nu N\int_0^tE^N(r)dr\right)-\frac{1}{N}P_{ir}\left(N\gamma\int_0^tI^N(r)dr\right)\\&\quad-\frac{1}{N}P_{di}\left(\mu N\int_0^tI^N(r)dr\right),\\
R^N(t)&=R^N(0)+\frac{1}{N}P_{ir}\left(N\gamma\int_0^tI^N(r)dr\right)-
\frac{1}{N}P_{rs}\left(N\rho\int_0^tR^N(r)dr\right)
\\&\quad-\frac{1}{N}P_{dr}\left(\mu N\int_0^tR^N(r)dr\right).
 \end{align*}
\end{example}

\begin{example}\label{ex:SEIRSb} {A variant of the SEIRS model with demography.}

In the preceding example, we decided to replace the true proportion of susceptibles by its approximation
$S(t)/N$, in order to avoid complications. There is another option, which is to force the population to remain constant. The most natural way to achieve this is to assume that each death event coincides with a birth event. Every susceptible, exposed, infected, removed individual dies at rate $\mu$. Each death is compensated by the birth of a susceptible. The equation for the evolution of $(S(t),E(t),I(t),R(t))$ reads

 \begin{align*}
 S(t)&=S(0)-P_{se}\left(\frac{\lambda}{N}\int_0^tS(r)I(r)dr\right)+P_{rs}\left(\rho\int_0^tR(r)dr\right)\\&\quad
 +P_{ds}\left(\mu\int_0^tS(r)dr\right)+P_{de}\left(\mu\int_0^tE(r)dr\right)+P_{di}\left(\mu\int_0^tI(r)dr\right)\\
 &\quad+P_{dr}\left(\mu\int_0^tR(r)dr\right)-P_{ds}\left(\mu\int_0^tS(r)dr\right),\\
 E(t)&=E(0)+P_{se}\left(\frac{\lambda}{N}\int_0^tS(r)I(r)dr\right)-P_{ei}\left(\nu\int_0^tE(r)dr\right)
 \\&\quad-P_{de}\left(\mu\int_0^tE(r)dr\right),\\
 I(t)&=I(0)+P_{ei}\left(\nu\int_0^tE(r)dr\right)-P_{ir}\left(\gamma\int_0^tI(r)dr\right)
 -P_{di}\left(\mu\int_0^tI(r)dr\right),\\
 R(t)&=R(0)+P_{ir}\left(\gamma\int_0^tI(r)dr\right)-P_{rs}\left(\rho\int_0^tR(r)dr\right)
 -P_{dr}\left(\mu\int_0^tR(r)dr\right).
 \end{align*}
  The equations for the proportions in the various compartments read
   \begin{align*}
 S^N(t)&=S^N(0)-\frac{1}{N}P_{se}\left(N\lambda\int_0^tS^N(r)I^N(r)dr\right)\\&
 \quad+\frac{1}{N}P_{rs}\left(N\rho\int_0^tR^N(r)dr\right)
 +\frac{1}{N}P_{ds}\left(N\mu\int_0^tS^N(r)dr\right)\\&\quad+\frac{1}{N}P_{de}\left(N\mu\int_0^tE^N(r)dr\right)+\frac{1}{N}P_{di}\left(N\mu\int_0^tI^N(r)dr\right)\\
 &\quad+\frac{1}{N}P_{dr}\left(N\mu\int_0^tR^N(r)dr\right)-\frac{1}{N}P_{ds}\left(N\mu\int_0^tS^N(r)dr\right),\\
  E^N(t)&=E^N(0)\!+\!\frac{1}{N}P_{se}\!\left(\!N\lambda\!\int_0^t\!S^N(r)I^N(r)dr\!\right)\!-\!\frac{1}{N}P_{ei}\!\left(\!N\nu\!\int_0^t\!E^N(r)dr\!\right)
 \\&\quad
 -\frac{1}{N}P_{de}\left(N\mu\int_0^tE^N(r)dr\right),\\
 I^N(t)&=I^N(0)+\frac{1}{N}P_{ei}\left(\nu\int_0^tE^N(r)dr\right)-\frac{1}{N}P_{ir}\left(\gamma\int_0^tI^N(r)dr\right)\\&\quad
 -\frac{1}{N}P_{di}\left(N\mu\int_0^tI^N(r)dr\right),\\
 R^N(t)&=R^N(0)+\frac{1}{N}P_{ir}\left(N\gamma\int_0^tI^N(r)dr\right)-\frac{1}{N}P_{rs}\left(N\rho\int_0^tR^N(r)dr\right)
 \\&\quad-\frac{1}{N}P_{dr}\left(N\mu\int_0^tR^N(r)dr\right).
 \end{align*}
\end{example}

In the three above examples, for each $j$,  $\beta_{N,j}(t,Nz)=N\beta_j(t,z)$, for some
$\beta_j(t,z)$ which does not depend upon $N$. We shall assume from now on that this is the case
in our general model, namely that
\[ \beta_{N,j}(t,Nz)=N\beta_j(t,z),\quad \text{for all } 1\le j\le k,\ N\ge1,\ z\in\R^d_+.\]

\begin{remark} We could assume more generally that
\[ \beta_{N,j}(t,Nz)=N\widetilde{\beta}_{N,j}(t,z),\ \text{ where }\widetilde{\beta}_{N,j}(t,z)\to\beta_j(t,z),\]
locally uniformly as $N\to\infty$.
\end{remark}

 Finally our model reads
 \begin{equation}\label{e:Nmodel}
  Z^N_t=x_N+\sum_{j=1}^k\frac{h_j}{N}P_j\left(\int_0^tN\beta_{j}(s,Z^N_s)ds\right).
  \end{equation}

We note that in the first example above, $0\le Z^N_j(t)\le1$ for all $1\le j\le k$, $t\ge0$, $N\ge1$.
In the second example however, such a simple upper bound does not hold, but a much weaker assumption will suffice.
%
%

We assume that all $\beta_j$ are locally bounded, which is clearly satisfied in all examples we can think of, so that
for any $K>0$,
\begin{equation}\label{beta_bound}
 C(T,K):=\sup_{1\le j\le k}\sup_{0\le t\le T}\sup_{|z|\le K}\beta_j(t,z)<\infty.
 \end{equation}
We first prove the Law of Large Numbers for Poisson processes.
\begin{proposition}
Let $\{P(t),\ t\ge0\}$ be a rate $\lambda$ Poisson process.
Then
\[t^{-1} P(t)\to\lambda\ \text{ a.s. as  }t\to\infty.\]
\end{proposition}
\bpf  Consider first for $n\in\Z_+$
\begin{align*}
n^{-1}P(n)&=n^{-1}\sum_{i=1}^n[P(i)-P(i-1)]\\
&\to\lambda\ \text{ a.s. as  }n\to\infty
\end{align*}
from the standard strong Law of Large Numbers, since the random variables $P(i)-P(i-1)$, $1\le i\le n$
are i.i.d.\ Poisson with parameter $\lambda$.
Now
\begin{align*}
t^{-1}P(t)&=\frac{[t]}{t}[t]^{-1}P([t])+t^{-1}\{P(t)-P([t])\},\\
\text{so }\ \left|t^{-1}P(t)-\lambda\right|&\le\left|\frac{[t]}{t}[t]^{-1}P([t])-\lambda\right|+t^{-1}\{P([t]+1)-P([t])\}.
\end{align*}
But
\begin{align*}
t^{-1}\{P([t]+1)-P([t])\}=t^{-1}P([t]+1)-t^{-1}P([t])
\end{align*}
is the difference of two sequences which converge a.s. towards the same limit, hence it converges to 0 a.s. \epf

Define the continuous time martingales (see Section \ref{sec:contmart} in Appendix A) $M_j(t)=P_j(t)-t$, $1\le j\le k$.
We have
\begin{align*}
Z^N_t=x_N+\int_0^tb(s,Z^N_s)ds+\sum_{j=1}^k\frac{h_j}{N}M_j\left(\int_0^tN\beta_{j}(s,Z^N_s)ds\right),
\end{align*}
where
\[ b(t,x)=\sum_{j=1}^k h_j\beta_j(t,x).\]

Consider the $k$-dimensional process $\mathcal{M}^N(t)$
whose $j$-th component is defined as
\[\mathcal{M}^N_j(t):=\frac{1}{N}M_j\left(N\int_0^t \beta_j(r,Z^N_r)dr\right).\]

From the above, we readily deduce the following.
\begin{proposition}\label{lgn_input}
For any $K>0$, let $\tau_K:=\inf\{t>0,\, |Z^N_t|\ge K\}$.
As $N\to\infty$, for all $T>0$, provided \eqref{beta_bound} holds,
\[\sup_{0\le t\le T\wedge\tau_K}|\mathcal{M}^N(t)|\to0\ \text{ a.s.}\]
\end{proposition}
\bpf
In order to simplify the notation we treat the case $d=1$.
It follows from \eqref{beta_bound} that, if $M(t)=P(t)-t$ and $N$ is large enough,
\begin{align*}
\sup_{0\le t\le T\wedge \tau_K}|\mathcal{M}^N(t)|\le \frac{1}{N}\sup_{0\le r\le N T C(T,K)}|M(r)|.
\end{align*}
From the previous proposition, for all $t>0$,
\[\frac{P(Nt)}{N}\to t\quad\text{a.s. as }N\to\infty.\]
Note that we have pointwise convergence of a sequence of increasing functions
towards a continuous (and of course increasing) function.
Consequently from the second Dini Theorem (see e.g.\ pages 81 and 270 in Polya and Szeg\"o \cite{PoSzI}), this convergence is uniform on any compact
interval, hence for all $T>0$,
\[ \frac{1}{N}\sup_{0\le r\le NTC(T,K)}|M(r)|\to0\ \text{ a.s.}\]
\epf

Concerning the initial condition, we assume that for some $x\in[0,1]^d$, $x_N=[Nx]/N$, where
$[Nx]$ is of course a vector of integers.
We can now prove the following theorem.
\begin{theorem}\label{th:LLN}
 {\bf Law of Large Numbers} Assume that the initial condition is given as above,
that $b(t,x)=\sum_{j=1}^k\beta_j(t,x)h_j$ is locally Lipschitz as a function of $x$, locally uniformly
in $t$, that \eqref{beta_bound} holds and that the unique solution of the ODE
\[  \frac{dz_t}{dt}=b(t,z_t),\quad z_0=x\]
does not explode in finite time.
Let $Z^N_t$ denote the solution of the SDE \eqref{e:Nmodel}. Then $Z^N_t\to z_t$ a.s. locally uniformly in $t$, where $\{z_t,\, t\ge0\}$ is the unique solution of the above ODE.
\end{theorem}
Needless to say, our theorem applies to the general model \eqref{e:Nmodel}. We shall describe below three specific models to which we can apply it. Note that if the initial fraction of infected is zero, then the fraction of infected is zero for all $t\ge0$.
\smallskip

\bpf
We have
\begin{align*}
Z^N_t=x_N+\int_0^tb(s,Z^N_s)ds+\sum_{j=1}^k h_j\mathcal{M}^N_j(t).
\end{align*}
Let us fix an arbitrary $T>0$. We want to show uniform convergence on $[0,T]$. Let $K:=\sup_{0\le t\le T}|z_t|+C$, where $C>0$ is arbitrary, and let $\tau_K=\inf\{t>0, |Z^N_t|\ge K\}$.
 Since $b(t,\cdot)$ is locally Lipschitz,
\[ c_{T,K}:=\sup_{0\le t\le T,\, x\not=x', |x|,|x'|\le K}\frac{|b(t,x)-b(t,x')|}{|x-x'|}<\infty.\]
For any $0\le t\le T$, if we define $Y^N_t=\sum_{j=1}^k h_j\mathcal{M}^N_j(t)$, we have
\begin{align*}
|Z^N_{t\wedge\tau_K}-z_{t\wedge\tau_K}|&\le|x_N-x|+c_{T,K}\int_0^{t\wedge\tau_K}|Z^N_s-z_s|ds+|Y^N_{t\wedge\tau_K}|\\
&\le\eps_N\exp{(c_{T,K} t)},
\end{align*}
where $\eps_N:=|x_N-x|+\sup_{0\le t\le T\wedge\tau_K}|Y^N_t|$ and we have used Gronwall's Lemma
\ref{le:Gron} below. It follows from our assumption on $x_N$ and Proposition \ref{lgn_input}
that $\eps_N\to0$ as $N\to\infty$. The result follows, since as soon as $\eps_N\exp{(c_{T,K} T)}\le C$, $\tau_K\ge T$. \epf

\begin{remark}
Showing that a stochastic epidemic model (for population proportions) converges to a particular deterministic process is important also for applications. This motivates the use of deterministic models, which are easier to analyse, in the case of large populations.
\end{remark}

\begin{lemma}\label{le:Gron} {\bf Gronwall}
Let $a,b\ge0$ and $\varphi:[0,T]\to\R$ be such that for all\linebreak $0\le t\le T$,
\[ \varphi(t)\le a + b\int_0^t\varphi(r)dr.\]
Then $\varphi(t)\le a e^{bt}$.
\end{lemma}
\bpf
We deduce from the assumption that
\[ e^{-bt}\varphi(t)-be^{-bt}\int_0^t\varphi(r)dr\le a e^{-bt},\]
or in other words
\[  \frac{d}{dt}\left(e^{-bt}\int_0^t\varphi(r)dr\right)\le a e^{-bt}.\]
Integrating this inequality, we deduce
\[  e^{-bt}\int_0^t\varphi(r)dr\le a\frac{1-e^{-bt}}{b}.\]
Multiplying by $be^{bt}$ and exploiting again the assumption yields the result.
\epf

\begin{example}{ The SIR model.}\label{ex:LLN-SIR}
It is clear that Theorem \ref{th:LLN} applies to Example \ref{example1}. The limit of $(S^N(t),I^N(t))$ is the solution $(s(t),i(t))$ of the ODE
\begin{align*}
s'(t)&=-\lambda s(t)i(t),\\
i'(t)&=\lambda s(t)i(t)-\gamma i(t).
\end{align*}
\end{example}

\begin{example}\label{Ex_SEIRS_3} {The SEIRS model with demography (continued).}
Again Theorem \ref{th:LLN} applies to Example \ref{ex:SEIRS}.  The limit of $(S^N(t),E^N(t),I^N(t),$ $R^N(t))$ is the solution $(s(t),e(t),i(t),r(t))$ of the ODE
\begin{align*}
s'(t)&=\mu(1-s(t))-\lambda s(t)i(t)+\rho r(t),\\
e'(t)&=\lambda s(t)i(t)-(\nu+\mu)e(t),\\
i'(t)&=\nu e(t)-(\gamma+\mu) i(t),\\
r'(t)&= \gamma i(t) -(\rho+\mu) r(t).
\end{align*}
Note that of we define the total renormalized population as $n(t)=s(t)+e(t)+i(t)+r(t)$, then
it is easy to deduce from the above ODE that $n'(t)=\mu(1-n(t))$, consequently $n(t)=1+e^{-\mu t}(n(0)-1)$.
If $n(0)=1$, then $n(t)\equiv1$, and we can reduce the above model to a three-dimensional model
(and to a two-dimensional model as in the previous example if we are treating the SIR or the SIRS model with demography).
\end{example}

We note that this ``Law of Large Numbers''  approximation is only valid when $s, i>0$, i.e.\
when significant fractions of the
population are infective and are susceptible, in particular at time $0$. The ODE is of course of no help to compute the probability
 that the introduction of a single infective results in a major epidemic.

The vast majority of the literature on mathematical models  in epidemiology considers ODEs  of the type
of equations which we have just obtained.
The probabilistic point of view is more recent.

\begin{exercise}\label{exoRossLGN}
Let us consider Ross's model of  malaria,
which we write in a stochastic form.
Denote by $H(t)$ the number of humans (hosts) who are infected by malaria, and by $V(t)$ the number of mosquitos (vectors)
who are infected by malaria at time $t$. Let $N_H$ denote the total number of humans, and $N_V$ denote the total number of mosquitos, which are assumed to be constant in time. The humans (resp.\ the mosquitos) which are not infected are all supposed to be susceptibles. Let
 $m=N_V/N_H$ and denote by $a$ the mean number of bites of humans
by one  mosquito per time unit, $p_{VH}$ the probability that the bite of a  susceptible human by an infected mosquito
infects the human, and by $p_{HV}$ the probability that a susceptible mosquito gets infected while biting an infected human. We assume that the infected humans
(resp.\ mosquitos)  recover at rate $\gamma$ (resp.\ at rate $\mu$).
\begin{enumerate}
\item
What is the mean number of bites that a human suffers per time unit?
\item Given  4 mutually independent standard Poisson processes $P_1(t)$, $P_2(t)$, $P_3(t)$ and $P_4(t)$,
 justify the following as a stochastic model of the propagation of malaria.
\begin{align*}
H(t)&=H(0)+P_1\left(\!ap_{VH}\!\int_0^t\!V(s)\frac{N_H-H(s)}{N_H}ds\right)-P_2\left(\!\gamma\!\int_0^t\!H(s)ds\right)\\
V(t)&=V(0)+P_3\left(\!a m p_{HV}\! \int_0^t\! H(s)\frac{N_V-V(s)}{N_V}ds\right)-P_4\left(\!\mu\!\int_0^t\! V(s)ds\right).
\end{align*}
\item Define now (with $N_H=N$, $N_V=m N$)
\[  h_N(t)=\frac{H(t)}{N_H},\quad v_N(t)=\frac{V(t)}{N_V}.\]
Write the equation for the pair $(h_N(t),v_N(t))$. Show that as $N\to\infty$, with $m$ constant, $(h_N(t),v_N(t))\to (h(t),v(t))$, the solution of Ross's ODE:
\begin{align*}
\frac{dh}{dt}(t)&=ap_{VH} m v(t)(1-h(t))-\gamma h(t),\\
\frac{dv}{dt}(t)&=ap_{HV} h(t)(1-v(t))-\mu v(t).
\end{align*}
\end{enumerate}
\end{exercise}

\section{Central Limit Theorem}\label{TB-EP_sec_CLT}

In the previous section we have shown that the stochastic process describing the evolution of the proportions of the total population in the various compartments converges, in the asymptotic of large population, to the deterministic solution of a system of ODEs.  In the current section we look at fluctuations of the difference between the stochastic epidemic process and its deterministic limit.

We now introduce the rescaled difference between $Z^N_t$ and $z_t$, namely
\[ U^N_t=\sqrt{N}(Z^N_t-z_t).\]
We wish to show that $U^N_t$ converges in law to a Gaussian process. It is clear that
\[ U^N_t=\sqrt{N}(x_N-x)+\sqrt{N}\int_0^t [b(s,Z^N_s)-b(s,z_s)]ds+\sum_{j=1}^k
h_j\m^N_j(t),\]
where for $1\le j\le k$,
\[ \m^N_j(t)=\frac{1}{\sqrt{N}}M_j\left(N\int_0^t\beta_j(r,Z^N_r)dr\right).\]
We certainly need to find the limit in law of the $k$ dimensional process
$\m^N_t$, whose $j$-th coordinate is $\m^N_j(t)$. We prove the following proposition below.
\begin{proposition}\label{cvmart}
As $N\to\infty$, \[\{\m^N_t,\, t\ge0\}\Rightarrow\{\m_t,\, t\ge0\}\]
meaning weak convergence
for the topology of locally uniform convergence,
where for $1\le j\le k$, $\m_j(t)=\int_0^t\sqrt{\beta_j(s,z_s)}dB_j(s)$ and the processes
$B_1(t),$ $\ldots,B_k(t)$ are mutually independent standard Brownian motions.
\end{proposition}
Let us first show that the main result of this section is indeed a consequence of this proposition.

\begin{theorem}\label{th:CLT} {\bf Central Limit Theorem}
In addition to the assumptions of Theorem \ref{th:LLN}, we assume that $x\to b(t,x)$ is of class $C^1$, locally uniformly
in $t$. Then,
as $N\to\infty$, $\{U^N_t,\, t\ge0\}\Rightarrow\{U_t,\, t\ge0\}$, where
\begin{equation}\label{eq:OU}
U_t=\int_0^t\nabla_x b(s,z_s)U_sds+\sum_{j=1}^k h_j\int_0^t\sqrt{\beta_j(s,z_s)}dB_j(s),\ t\ge0.
\end{equation}
\end{theorem}
\bpf
We shall fix an arbitrary $T>0$ throughout the proof.
Let $V^N(s):=\sqrt{N}[b(s,Z^N_s)-b(s,z_s)]$ and $\n^N_t:=\sum_{j=1}^kh_j\m^N_j(t)$. We have
\[ U^N_t=U^N_0+\int_0^tV^N(s)ds+\n^N_t.\]
Let us admit for the moment the following lemma.
\begin{lemma}\label{mean-value}
For each $N\ge1$, $0\le t\le T$ there exists a random $d\times d$ matrix $A^N_t$ such that
\[ V^N_t=\nabla b(t,z_t)U^N_t+A^N_t U^N_t.\]
Moreover, $\sup_{0\le t\le T}\| A^N_t\|\to0$, a.s., as $N\to\infty$.
\end{lemma}
We clearly have
\[ U^N_t=U^N_0+\int_0^t[\nabla b(s,z_s)+A^N_s ]U^N_sds+\n^N_t.\]
It then follows from Gronwall's Lemma that
\[\sup_{0\le t\le T}|U^N_t|\le\left(|U^N_0|+\sup_{0\le t\le T}|\n^N_t|\right)
\exp\left(\sup_{0\le t\le T}\|\nabla b(t,z_t)+A^N_t\| \, T\right).\]
The right-hand side of this inequality is tight,\footnote{A sequence $\xi_n$ of $\R_+$-valued random variables is tight if for any $\eps>0$, there exists an $M_\eps$ such that $\P(\xi_n>M_\eps)\le\eps$, for all $n\ge1$, see Section \ref{tight} in the Appendix.} hence the same is true for the
left-hand side. From
this and Lemma \ref{mean-value} it follows that $R^N_t:=A^N_t U^N_t$ tends to $0$ in probability
as $N\to\infty$, uniformly for $0\le t\le T$.
 Consequently
\begin{align*}
 U^N_t&=\int_0^t\nabla_xb(s,z_s)U^N_sds+W^N_t,\ \text{where}\\
W^N_t&=U^N_0+\int_0^tR^N_sds+\n^N_t.
\end{align*}
The following two hold
\begin{enumerate}
\item
$\sup_{0\le t\le T}|U^N_0+\int_0^tR^N_sds|\to0$ in probability, and from
Proposition \ref{cvmart} $\n^N_t\Rightarrow
\n_t$, hence $W^N_t\Rightarrow\n_t$ for the topology of uniform convergence on $[0,T]$.
\item
The mapping $y\mapsto\Phi(y)$, which to $y\in C([0,T];\R^d)$ associates $x\in C([0,T];\R^d)$, the solution of the ODE
\[ x(t)=\int_0^t\nabla b(s,z_s)x(s)ds+y(t),\]
is continuous.

Indeed, we can construct this mapping by first solving the ODE
\[\dot{z}(t)=\nabla b(t,z_t)[z(t)+y(t)],\quad z(0)=0,\]
and then defining $x(t)=z(t)+y(t)$.
\end{enumerate}
Since
\[U^N=\Phi(W^N),\]
 the result follows from 1. and 2., and the fact that $T$ is arbitrary.
 \epf

\begin{proof}[Proof of Lemma \ref{mean-value}]
For $1\le i\le d$, $0\le t\le T$, define the random function $\rho_{i,t}(\theta)=b_i(t,z_t+\theta(Z^N_t-z_t))$, $0\le \theta\le1$.
The mean value theorem applied to the function $\rho_{i,t}$ implies that for all $0\le i\le d$, $0\le t\le T$, there exists
a random $0<\bar{\theta}_{i,t}<1$ such that \[b_i(t,Z^N_t)-b_i(t,z_t)=\langle\nabla b_i(t,z_t+\bar{\theta}_{i,t}(Z^N_t-z_t)),Z^N_t-z_t\rangle.\]
Applying the same argument for all $1\le i\le d$ yields the first part of the Lemma. Theorem
\ref{th:LLN} and the continuity in $z$ of $\nabla b(t,z)$ uniformly in $t$ imply that\linebreak
$\nabla b_i(t,z_t+\bar{\theta}_{i,t}(Z^N_t-z_t))-\nabla b_i(t,z_t)\to0$ a.s., uniformly in $t$, as $N\to\infty$.
\epf

It remains to prove Proposition \ref{cvmart}. Let us first establish a central limit theorem for standard Poisson processes. Let $\{P_j(t),\ t\ge0\}_{1\le j\le k}$ be $k$ mutually independent standard Poisson processes
 and
$M(t)$ denote the $k$-dimensional process whose $j$-th component is $P_j(t)-t$.
\begin{lemma}\label{TCL-PoisProc}
As $N\to\infty$,
\[ \frac{M(Nt)}{\sqrt{N}}\Rightarrow B(t),\]
where $B(t)$ is a $k$-dimensional standard Brownian motion (in particular $B(t)\sim\mathcal{N}(0,tI)$,
with $I$ the $d\times d$ identity matrix) and the convergence is in the sense of convergence in law
in $D([0,+\infty);\R^k)$.
\end{lemma}
For a definition of the space $D([0,+\infty);\R^k)$ of the $\R^k$-valued c\`al\`ag functions of $t\in[0,\infty)$ and its topology, see section A.5 in the Appendix.
\smallskip

\bpf
It suffices to consider each component separately, since they are independent. So we do as if $k=1$.
We first note that our process is a martingale, whose associated predictable increasing process
is given by $\langle N^{-1/2}M(N\cdot),N^{-1/2}M(N\cdot)\rangle_t=t$. Hence it is tight.

Let us now compute the characteristic function of the random variable\linebreak $N^{-1/2}M(Nt)$.
We obtain
\[\E\left(\exp\left[iuN^{-1/2}M(Nt)\right]\right)=\exp\left(Nt\left[e^{i\frac{u}{\sqrt{N}}}-1-i\frac{u}{\sqrt{N}}\right]\right)\to\exp\left(-t\frac{u^2}{2}\right),\]
as $N\to\infty$. This shows that $N^{-1/2}M(Nt)$ converges in law to an $\mathcal{N}(0,t)$ r.v.

Now let $n\ge1$ and $0<t_1<\dots<t_n$. The random variables $N^{-1/2}M(Nt_1)$, $N^{-1/2}M(Nt_2)-N^{-1/2}M(Nt_1),
\ldots,N^{-1/2}M(Nt_n)-N^{-1/2}M(Nt_{n-1})$ are mutually independent and, if $B(t)$ denotes a standard
one dimensional Brownian motion, the previous argument shows that, with $M(0)=B(0)=0$, for any $1\le k\le n$,
$N^{-1/2}(M(Nt_k)-M(Nt_{k-1}))\Rightarrow B(t_k)-B(t_{k-1})$.
 Thus, since the random variables
$B(t_1), B(t_2)-B(t_1), \ldots,B(t_n)-B(t_{n-1})$ are mutually independent, we have shown that
\begin{align*}
&\left(\frac{M(Nt_1)}{\sqrt{N}},\frac{M(Nt_2)-M(Nt_1)}{\sqrt{N}},\ldots,\frac{M(Nt_n)-M(Nt_{n-1})}{\sqrt{N}}\right)\\
&\qquad\Rightarrow(B(t_1), B(t_2)-B(t_1), \ldots,B(t_n)-B(t_{n-1}))
\end{align*}
 as $N\to\infty$. This proves that the finite dimensional distributions of the process $N^{-1/2}M(Nt)$ converge to those of $B(t)$. Together
with tightness, this shows the lemma. \epf

\begin{proof}[Proof of Proposition \ref{cvmart}]
With the notation of the previous lemma,
\[ \M^N_j(t)=N^{-1/2}M_j\left(N\int_0^t \beta_j(s,Z^N_s)ds\right).\]
We write
\[\M^N_j(t)=N^{-1/2}M_j\left(N\int_0^t\beta_j(s,z_s)ds\right)+\widetilde{\m}^N_j(t),\]
where
\[\widetilde{\m}^N_j(t)=N^{-1/2}M_j\left(N\int_0^t\beta_j(s,Z^N_s)ds\right)-N^{-1/2}M_j\left(N\int_0^t\beta_j(s,z_s)ds\right).\]
For  $C>0$, let $\tau_{N,C}=\inf\{t>0,\, |Z^N_t|> C\}$.  We assume for a moment the identity
\begin{equation}\label{ScdMoment}
 \E\!\left(\!\left|\widetilde{\m}^N_j(t\wedge\tau_{N,C})\right|^2\right)\!\!=\!\E\left(\int_0^{t\wedge\tau_{N,C}}\!\!\left|\beta_j(s,Z^N_s)ds-\int_0^t\beta_j(s,z_s)ds\right|\right).\!\!\!\!
\end{equation}
The above right-hand side is easily shown to converge to $0$ as $N\to\infty$. Jointly with Doob's inequality from Proposition
\ref{Doob2c} in the Appendix, this shows that for all $T>0$, $\eps>0$,
\begin{align*}
\P\left(\sup_{0\le t\le T}\left|\widetilde{\m}^N_j(t)\right|>\eps\right)&\le \P(\tau_{N,C}<T)+\P\left(\sup_{0\le t\le T\wedge\tau_{N,C}}
\left|\widetilde{\m}^N_j(t)\right|>\eps\right)\\
&\le  \P(\tau_{N,C}<T)\\
&\qquad+\frac{4}{\eps^2}\E\left(\int_0^{t\wedge\tau_{N,C}}\left|\beta_j(s,Z^N_s)ds-\int_0^t\beta_j(s,z_s)ds\right|\right).
\end{align*}
It follows from Theorem \ref{th:LLN} that for $C>0$ large enough, both terms on the right tend to $0$, as $N\to\infty$. Consequently
$\sup_{0\le t\le T}\left|\widetilde{\m}^N_j(t)\right|\to0$ in probability as $N\to\infty$.

It remains to note that an immediate consequence of Lemma \ref{TCL-PoisProc} is that
\[N^{-1/2}M_j\left(N\int_0^t\beta_j(s,z_s)ds\right)\Rightarrow B_j\left(\int_0^t\beta_j(s,z_s)ds\right)\]
in the sense of weak convergence in the space $D((0,+\infty);\R)$, and the coordinates are mutually independent. However the two processes
$B_j\!\left(\!\int_0^t\!\beta_j(s,z_s)ds\!\right)$ and $\int_0^t\sqrt{\beta_j(s,z_s)}dB_j(s)$ are two centered Gaussian processes which have the same covariance functions. Hence they have the same law.

We finally need to establish \eqref{ScdMoment}. Following the development in Section \ref{sec-Po-proc} in the Appendix, we can rewrite the local martingale $\widetilde{\m}^N_j(t)$ as follows, forgetting the index $j$, and the time parameter of $\beta$ for the sake of simplifying notations
\begin{align*}
\widetilde{\m}^N(t)&=N^{-1/2}\int_0^t\int_0^\infty{\bf1}_{\{N\beta(z_{s})\le u\le N\beta(Z^N_{s-})\}}\overline{Q}(ds,du)\\
&\quad-N^{-1/2}\int_0^t\int_0^\infty{\bf1}_{\{N\beta(Z^N_{s-})\le u\le N\beta(z_{s})\}}\overline{Q}(ds,du),
\end{align*}
where $\overline{Q}(ds,du)=Q(ds,du)-dsdu$ and $Q$ is a standard Poisson point measure on $\R^2_+$. Noting that the square of each jump of the above martingale equals $N^{-1}$, we deduce from Proposition \ref{SecondMomentMart} in the Appendix that
\begin{align*}
\E\left(\left|\widetilde{\m}^N(t\wedge\tau_{N,C})\right|^2\right)&=N^{-1}\E\int_0^{t\wedge\tau_{N,C}}\int_{N[\beta(Z^N_{s-}\wedge\beta(z_s)]}^{N[\beta(Z^N_{s-}\vee\beta(z_s)]}Q(ds,du)\\
&=N^{-1}\E\int_0^{t\wedge\tau_{N,C}}\int_{N[\beta(Z^N_{s}\wedge\beta(z_s)]}^{N[\beta(Z^N_{s}\vee\beta(z_s)]}dsdu,
\end{align*}
which yields \eqref{ScdMoment}.
\epf

\begin{example} The $SIR$ model. It is clear that Theorem \ref{th:CLT} applies to Example \ref{example1}.
If we define $\begin{pmatrix}U_t \\ V_t\end{pmatrix}=
\lim_{N\to\infty}\sqrt{N}\begin{pmatrix}S^N(t)-s(t)\\ I^N(t)-i(t)\end{pmatrix}$, we have
\begin{align*}
U_t&=-\lambda\int_0^t\left[i(r)U_r+s(r)V_r)\right]dr-\int_0^t\sqrt{\lambda s(r)i(r)}dB_1(r),\\
V_t&=\int_0^t\left[\lambda(i(r)U_r+s(r)V_r)-\gamma V_r\right]dr+\int_0^t\sqrt{\lambda s(r)i(r)}dB_1(r)\\
&\qquad\qquad\qquad\qquad-\int_0^t\sqrt{\gamma i(r)}dB_2(r).
\end{align*}
\end{example}

\begin{remark}
Consider now the SIR model, started with a fixed small number of infectious individuals, all others being susceptible, so that $(S^N(0),I^N(0))$ $\to(1,0)$, as $N\to\infty$. The solution of the ODE from Example \ref{ex:LLN-SIR} starting from $(s(0),i(0))=(1,0)$ is the constant
$(s(t),i(t))\equiv(1,0)$. So in that case the coefficients of the noise in the last example are identically $0$, and, the initial condition of the stochastic model being deterministic, it is natural to assume that $(U_0,V_0)=(0,0)$. Then $(U_t,V_t)\equiv(0,0)$. Consequently Theorem \ref{th:CLT} tells us that, as $N\to\infty$,
for any $T>0$,
\[ \sqrt{N}\binom{S^N(t)-1}{I^N(t)-0}\to0,\ \text{ in probability, uniformly w.r.t. }t\in[0,T].\]
In the case $R_0>1$, i.e.\  $\lambda>\gamma$, with positive probability the epidemic gets off. However, as we shall see in Section \ref{sec:duration-epid} below, this take time of the order of $\log(N)$, and there is no contradiction with the present result.
\end{remark}

We close this section by a discussion of some of the properties of solutions of linear SDEs of the above type, following some of the developments in section 5.6 of Karatzas and Shreve \cite{KSI}. Suppose that $\{A(t),\, t\ge0\}$ and $\{C(t),\, t\ge0\}$ are $d\times d$ matrix-valued measurable and locally bounded deterministic functions of $t$. With $\{B(t),\, t\ge0\}$ being a $d$-dimensional Brownian motion, we consider the SDE
\[ dX_t=A(t)X_tdt+C(t)dB_t,\, t\ge0,\]
$X_0$ being a given $d$-dimensional Gaussian random vector independent of the Brownian motion
$\{B(t)\}$. The solution to this SDE is the $\R^d$-valued process given by the explicit formula
\[ X(t)=\Gamma(t,0)X_0+\int_0^t\Gamma(t,s)C(s)dB_s,\]
where the $d\times d$ matrix $\Gamma(t,s)$ is defined for all $0\le s\le t$ as follows. For each fixed
$s\ge0$, $\{\Gamma(t,s),\, t\ge s\}$ solves the linear ODE
\[ \frac{d\Gamma(t,s)}{dt}=A(t)\Gamma(t,s),\quad \Gamma(s,s)=I,\]
where $I$ denotes the $d\times d$ identity matrix.  It follows that $\{X_t,\, t\ge0\}$ is a Gaussian
process, and for each $t>0$, the mean and the covariance matrix of the Gaussian random vector $X_t$ are given by (denoting by $C^\ast$ the transpose of the matrix $C$)
\begin{align*}
 \E(X_t)&=\Gamma(t,0)\E(X_0),\\ \text{Cov}(X_t)&=\Gamma(t,0)\text{Cov}(X_0)\Gamma^\ast(t,0)
+\int_0^t\Gamma(t,s)C(s)C^\ast(s)\Gamma^\ast(t,s)ds.
\end{align*}
Assume now that $A(t)\equiv A$ and $C(t)\equiv C$ are constant matrices. Then
$\Gamma(t,s)=\exp((t-s)A)$. If we define $V(t):=\text{Cov}(X_t)$, we have that
\[ V(t)=e^{tA}\left[V(0)+\int_0^te^{-sA}CC^\ast e^{-sA^\ast}ds\right]e^{tA^\ast}.\]
If we assume moreover that all the eigenvalues of $A$ have negative real parts, then it is not hard to show that
as $t\to\infty$,
\[V(t)\to V:=\int_0^\infty e^{sA}CC^\ast e^{sA^\ast}ds.\]  In that case the Gaussian law with mean zero and covariance matrix $V$ is an invariant distribution of Gauss--Markov process $X_t$. This means in particular that if $X_0$ has that distribution, then the same is true for $X_t$ for all $t>0$.
We now show the following result, which is often useful for computing the covariance matrix $V$ in particular cases.
\begin{lemma}\label{le:cov_inv}
Under the above assumptions on the matrix $A$,
$V$ is the unique $d\times d$ positive semidefinite symmetric matrix which satisfies
\[ AV+VA^\ast+CC^\ast=0.\]
\end{lemma}
\bpf
Uniqueness follows from the fact that the difference $\bar{V}$ of two solutions satisfies
$A\bar{V}+\bar{V}A^\ast=0$. This implies that for all $x\in\R^d$, $\langle A\bar{V} x,x\rangle=0$.
Since none of the eigenvalues of $A^\ast$ is zero, this implies that $\langle \bar{V} x,x\rangle=0$ for
all eigenvectors $x$ of $A^\ast$, hence for all $x\in \R^d$. Since $\bar{V}$ is symmetric, this implies
that $\bar{V}=0$.

To show that $V$ satisfies the wished identity, assume that the law of $X_0$ is Gaussian  with mean $0$ and Covariance matrix $V$. Then $V$ is also the covariance matrix of $X_t$. Consequently
\[ V=e^{tA}Ve^{tA^\ast}+\int_0^te^{(t-s)A}CC^\ast e^{(t-s)A^\ast}ds.\]
Differentiating with respect to $t$, and letting $t=0$ yields the result.
\epf

We leave the last result as an exercise for the reader.
\begin{exercise}\label{exerOU}
Consider again the case of time varying matrices $A(t)$ and $C(t)$. We assume that $A(t)\to A$ and
$C(t)\to C$ as $t\to\infty$, and moreover that the real parts of all the eigenvalues of $A$ are negative.
Conclude that the law of $X_t$ converges to the Gaussian law with mean $0$ and covariance matrix
$V$ defined as above.
\end{exercise}

\section{Diffusion Approximation}\label{TB-EP_sec_DiffusApprox}

We consider again the vector of proportions in our model as
\begin{equation}\label{model}
Z^N(t)=x+\frac{1}{N}\sum_{j=1}^kh_jP_j\left(\int_0^tN\beta_j(Z^N(s))ds\right).
\end{equation}
From the strong law of large numbers, $\sup_{0\le t\le T}\|Z^N(t)-z_t\|\to0$ almost surely
as $N\to\infty$, for all $T>0$, where $z_t$ solves the ODE
\[\dot{z}_t=b(z_t),\ z_0=x;\ \text{ where } b(x)=\sum_{j=1}^kh_j\beta_j(x).\]
 We now consider a diffusion approximation $X^N_t$ of the above model, which
solves the SDE
\[X^N_t=x+\int_0^tb(X^N_s)ds+\sum_{j=1}^k\frac{h_j}{\sqrt{N}}\int_0^t\sqrt{\beta_j(X^N_s)}dB^j_s,\]
where $B^1,\ldots,B^k$ are mutually independent standard Brownian motions.
Let us define the Wasserstein-$1$ distance on the interval $[0,T]$
between two $\R^d$-valued processes $U_t$ and $V_t$ as
\[W_{1,T}(U,V)=\inf\E\left(\|U-V\|_T\right),\]
where, if $x:[0,T]\to\R^d$, $\|x\|_T=\sup_{0\le t\le T}\|x(t)\|$, and the above infimum is over all couplings of the two
processes  $U(t)$ and $V(t)$, i.e.\ over all ways of defining jointly the two processes, while respecting the two
marginal laws of $U$ and $V$. We shall use the two following well--known facts about the Wasserstein distance: it is a distance (and satisfies the triangle inequality);
if $U_n$ is a sequence of random elements of $D([0,T];\R^d)$ which converges in law to a continuous process $U$, and is such that the sequence of random variables $\| U_n\|_T$ is uniformly integrable, then $W_{1,T}(U_n,U)\to0$ as $n\to\infty$.

The aim of this section is to establish the following theorem.
\begin{theorem}\label{th:diffapprox}
For all $T>0$, as $N\to\infty$,
\[ \sqrt{N}W_{1,T}(Z^N,X^N)\to0,\]
or in other words, $W_{1,T}(Z^N,X^N)=\text{o}(N^{-1/2})$.
\end{theorem}
\bpf
We have proved in Theorem \ref{th:LLN} that $\sup_{0\le t\le T}\|Z^N_t-z_t\|_T\to0$ almost surely, as $N\to\infty$, and moreover
$\sqrt{N}(Z^N-z)\Rightarrow U$ as $N\to\infty$, where the above convergence holds for the
topology of uniform convergence on the interval $[0,T]$, and $U$ is the Gaussian process solution of the SDE
\[U_t=\int_0^t\nabla b(z_s)U_sds+\sum_{j=1}^k h_j\int_0^t\sqrt{\beta_j(z_s)}dB^j_s.\]
It is not hard to prove the following.
\begin{exercise} As $N\to\infty$, $\sup_{0\le t\le T}\|X^N_t-z_t\|_T\to0$ almost surely, and moreover
$\sqrt{N}(X^N-z)\Rightarrow U$.
\end{exercise}
We first note that from the triangle inequality
\begin{align*}
W_{1,T}(\sqrt{N}(Z^N - z),\sqrt{N}(X^N - z))&\le W_{1,T}(\sqrt{N}(Z^N-z),U)\\
&\qquad+W_{1,T}(\sqrt{N}(X^N-z),U)\\
&\to0,
\end{align*}
as $N\to\infty$.
Moreover
\begin{align*}
W_{1,T}(\sqrt{N}(Z^N\!\!-\!z),\sqrt{N}(X^N\!\!-\!z))&=\!\!\!\inf_{\text{couplings}}\!\!\!\E\|\sqrt{N}(Z^N\!\!-\!z)\!-\!\sqrt{N}(X^N\!\!-\!z)\|_T\\
&=\!\!\!\inf_{\text{couplings}}\!\!\!\sqrt{N}\E\|(Z^N-z)-(X^N-z)\|_T\\
&=\!\!\!\inf_{\text{couplings}}\!\!\!\sqrt{N}\E\|Z^N-X^N\|_T\\
&=\sqrt{N}W_{1,T}(Z^N,X^N).
\end{align*}
Theorem \ref{th:diffapprox} follows from the two last computations. \epf

\begin{remark}
If we combine the law of large numbers and the central limit theorem which have been established in the previous two sections, we conclude that $Z^N_t-z_t-N^{-1/2}U_t = \circ(N^{-1/2})$. In other words, if we replace $Z^N_t$ by the Gaussian process $z_t+N^{-1/2}U_t$, the error we make, at least on any given
finite time interval, is small compared to $N^{-1/2}$.  The same is true for the diffusion approximation
$X^N_t$.
\end{remark}

\chapter{General Closed Models}  

In this chapter we go back to the general model, i.e.\ not assuming exponential latent and infectious periods implying that the epidemic process is Markovian. We consider models which are closed in the sense that there is no influx of new susceptibles during the epidemic.
No birth, no immigration, and the removed individual are either dead or recovered, with an immunity which they do not lose in the considered time frame.

In this context, the epidemic will stop sooner or later. The questions of main interest are: the evaluation of the duration of the epidemic, and the total number of individuals which are ever infected. The first section gives exact results concerning the second issue in small communities. The rest of the chapter is concerned with large communities. We present the Sellke construction, and then use it to give a law of large number and a central limit theorem for the number of infected individuals. Finally we study the duration of the epidemic.

\section[Exact results for the final size in small communities]{Exact results for the final size in small communities}\label{Sec_Exact}

In earlier sections it is often assumed that the population size $N$ is large. In other situations this is not the case, for example in planned infectious disease experiments in veterinary science the number of studied animals is of the order 5--20 (e.g.\ Quenee et al.\ \cite{Q2011I}), and in such cases law of large numbers and central limit theorems have not yet kicked in, which motivates the current section about exact results in small populations.

It turns out that it is quite complicated to derive expressions for the distribution of the final size, even when $N$ is quite small. The underlying reason for this is that there are many ways in which an outbreak can result in exactly $k$ initially susceptible individuals getting infected. We illustrate this by computing the final size distribution $\{ p^{(N)}_k\}$ for the Reed--Frost model for $N=1$, 2 and 3. We then derive a recursive formula for the final outcome of the full model valid for general $N$ and $k$ (but numerically unstable for $N$ larger than, say, 40).

Consider the Reed--Frost epidemic where the probability to infect a given susceptible equals $p$ ($=1-e^{-\lambda \iota/N}$). And let $N=1$, one susceptible and one infectious individual to start with. The possible values of $Z$ are then 0 and 1, and obviously we have $p^{(1)}_0=\P (Z=0|N=1)=1-p$ and $p^{(1)}_1=p$. For $N=2$ things are slightly more complicated. No one getting infected is easy: $p^{(2)}_0=(1-p)^2$, since both individuals have to escape infection from the index case. For $Z=1$ to occur, the index case must infect exactly one of the two remaining, but further, this individual must not infect the third person: $p^{(2)}_1=\binom{2}{1}p(1-p)*(1-p)$. Finally, the probability of $Z=2$ is of course the complimentary probability, but it can also be obtained by considering the two possibilities for this to happen: either the index case infects both, or else the index case infects exactly one of the two, and that individual in turn infects the remaining individual: $p^{(2)}_2=p^2+2p(1-p)*p$.

For $N=3$ initial susceptibles the situation becomes even more complicated. It is best to write down the different epidemic generation chains at which individuals get infected. We always have one index case. The chain in which the index case infects two individuals who in turn together infect the last individual,  is denoted $1\to 2\to 1\to 0$. The probability for such a chain can be computed sequentially for each generation keeping in mind: how many susceptibles there are at risk, how many that get infected and what is the risk of getting infected (the complimentary probability of escaping infection). The probability for the chain just mentioned is given by
$$
\P (1\to 2\to 1\to 0|N=3)= \binom{3}{2}p^2(1-p)^1*(1-(1-p)^2).
$$
The last factor comes from the final individual getting infected when there were two infected individuals in the previous generation (so the escape probability equals $(1-p)^2$). We hence see that the probability of a chain is the product of (different) binomial probabilities. The final size probabilities are then obtained by writing down the different possible chains giving the desired final outcome:
\begin{align*}
p^{(3)}_0 & =\P (1\to 0)=(1-p)^3
\\
p^{(3)}_1 & =\P (1\to 1\to 0)=\binom{3}{1}p(1-p)^2 * (1-p)^2
\\
p^{(3)}_2 & = \P (1\to 2\to 0)+\P (1\to 1\to 1\to 0)
\\
     & =\binom{3}{2}p^2(1-p)*((1-p)^2) + \binom{3}{1}p(1-p)^2 *\binom{2}{1}p(1-p)*(1-p)
\\
p^{(3)}_3 & = \P (1\to 3\to 0)+\P (1\to 2\to 1\to 0) + \P (1\to 1\to 2\to 0)
\\ & \hskip2cm + \P (1\to 1\to 1\to 1\to 0)
\\
     &= \dots
\end{align*}

\begin{exercise} \label{x3.1}
Compute $p_3^{(3)} $ explicitly by computing the probabilities of the different chains. Check that $\sum_{k=0}^3p_3^{(k)}=1$ for any $p\in [0,1]$.
\end{exercise}

For general $N$ it is possible to write down the outcome probability for a specific chain as follows. If we denote the number of susceptibles and infectives in generation $k$ by $(S_k,I_k)$, then the epidemic starts with $(S_0, I_0)=(s_0,i_0)=(N,1)$. From a chain $1\to i_1\to \dots i_j\to 0$ (so $i_{j+1}=0$) the number of susceptibles in generation $k$ is also known from the relation $s_k=s_0-\sum_{j=1}^ki_k$. We use this when we compute the binomial probabilities of a given generation of the chain, these binomial probabilities depend on: how many were at risk, how many infectives there were in the previous generation, and how many to be infected in the current. Finally, the probability of a chain is the product of the different binomial probabilities of the different generations. From this we obtain the following so called chain-binomial probabilities
$$
\P (1\to i_1\to \dots i_j\to 0)= \prod_{k=1}^{j+1} \binom{s_{k-1}}{i_k}\left( 1-(1-p)^{i_{k-1}}\right)^{i_k}\left( (1-p)^{i_{k-1}}\right)^{s_{k-1}-i_k}.
$$
As seen, these expression are quite long albeit explicit. However, computing the final outcome probabilities $p_N(k),\ k=0,\dots , N$, is still tedious since there are many different possible chains resulting in exactly $k$ getting infected at the end of the epidemic. Further, things become even more complicated when considering different distributions of the infectious period than a constant infectious period as is assumed for the Reed--Frost epidemic model.

However, it is possible to derive a recursive formula for the final number infected $p_N(k)$, see e.g.\ Ball \cite{B86I}, which we now show. The derivation of the recursion of the final size uses two
main ideas: a Wald's identity for the final size and the total
infection pressure, and the exchangeability of individuals making
it possible to express the probability of having $k$
infections among the initially $N$ susceptibles in terms of the probability of getting all $k$ infected
in the subgroup containing those $k$ individuals and the index case, and the probability that the remaining $N-1-k$ individuals escape infection from that group.

Let us start with the latter. Fix $N$ and write $\bar \lambda=\lambda/N$. As before we let $Z^N$ denote the total number infected excluding the index case(s), explicitly showing the dependence on the number of initially susceptible $N$. Since individuals are exchangeable we can label the individuals according to the order in which they get infected. The index case is labelled $0$, the individuals who get infected during the outbreak are labelled: $1,\dots , Z^N$, and those who avoid infection according to any order $Z^N+1, \dots ,N$. With this labelling we define the total infection pressure $A^{N}$ by
\begin{equation}
A^{N}=\bar \lambda\sum_{i=0}^{Z^N} I_i \label{tot-inf-pres}
\end{equation}
i.e.\ the infection pressure, exerted on any individual, during the complete outbreak (sometimes referred to as the ``total cost" or the ``severity'' of the epidemic).

As earlier we let $p_i^{(N)}=\P (Z^N=i)$ denote the probability that exactly $k$ initial susceptibles out of $N$ get infected during the outbreak. Reasoning in terms of subsets among the initial susceptibles as described earlier, and using the exchangeability of individuals, it can be shown (\cite{B86I}) that for any $i\le k\le N$,
\begin{equation}
\frac{p_i^{(N)}}{\binom{N}{i}} = \frac{p_i^{(k)}}{\binom{k}{i}} \E\left(e^{-(N-k)A^{k}}|Z^{k}=i\right).\label{subset-prob}
\end{equation}
The equation is explained as follows. On the left-hand side is the probability that a specific group of size $i$ (out of $N$) get infected and no one else. On the right-hand side this event is divided into two sub events.This is done by considering another group of size $k\ge i$, containing the earlier specified group of size $i$ as a subset. The first factor is then the probability that exactly the subgroup of size $i$ get infected within the bigger group of size $k$. The second factor, the expectation, is the probability that all individuals outside the bigger subgroup avoid getting infected. The notation $A^{k}$ and $Z^{k}$ hence denote the total pressure and final size starting with $k$ susceptibles.

We use the following steps to show Wald's identity recalling that $\psi_I(b)=\E (e^{b I})$ is the moment generating function of the infectious period (so $\psi_I(-b)$ is the Laplace transform)
\begin{align*}
(\psi _I(-\theta \bar \lambda))^{k+1}&= \E\left[ \exp\left(-\theta \bar \lambda\sum_{i=0}^{k}I_i \right)\right]
\\
&= \E\left[ \exp\left(-\theta \left(A^{k}+\bar \lambda\sum_{i=Z^{k}+1}^{k}I_i\right) \right)\right]
\\
&= \E\left[ e^{-\theta A^{k}}(\psi_I (-\theta \bar \lambda)^{k-Z^{k}} \right].
\end{align*}
The last identity follows since the $k-Z^{k}$ infectious periods
$I_{Z^{k}+1},\dots I_{k}$, are mutually independent and
also independent jointly of the total pressure $A^{k}$ (which only depends on the first $Z^{k}$ infectious periods and the contact processes of these individuals).
If we now divide both sides by $(\psi_I (-\theta \bar \lambda))^{k+1}$ we obtain Wald's identity for $Z^{k}$ and $A^{k}$:
\begin{equation}
\E\left(\frac{e^{-\theta A^{k}}}{(\psi_I (\theta \bar \lambda))^{1+Z^{k}}}\right)
= 1,\qquad \theta \ge 0.\label{Wald}
\end{equation}

If we apply Wald's identity with $\theta =N-k$ and condition on the value of $Z^{k}$ we get
\begin{equation}
\sum_{i=0}^k\frac{\E\left(e^{-(N-k)A^{k}}|Z^{k}=i\right)}{(\psi_I
  (-(N-k) \bar \lambda))^{i+1}}p_i^{(k)}=1.
\label{wald-cond}
\end{equation}
If we now use Equation (\ref{subset-prob}) in the equation above we
get
\[
\sum_{i=0}^k \frac{\binom{k}{i}p_i^{(N)}}{\binom{N}{i}(\psi_I
  (-(N-k) \bar \lambda))^{i+1}} =1.
\]
Simplifying the equation, returning to $\lambda=\bar \lambda N$ and
putting $p_k^{(N)}$ on one side, we obtain the recursive formula for the final size distribution $p_k^{(N)}, k=0,\dots ,N$.

\begin{theorem}
The exact final size distribution is given by the recursive formula
\begin{equation}
p_k^{(N)} = \binom{N}{k}[\psi_I (-(N-k) \lambda/N)]^{k+1}
-\sum_{i=0}^{k-1} \binom{N-i}{k-i}[\psi_I
  (-(N-k) \lambda/N)]^{k-i}p_i^{(N)}.
  \label{fin-size}
\end{equation}
\end{theorem}

For example,
solving Equation (\ref{fin-size}) for $k=0$ (when the sum is vacuous) and then for $k=1$ gives, after some algebra,
\begin{align*}
p_0^{(N)}&=\psi_I\left(\lambda\right),
\\
p_1^{(N)}& = N\psi_I\left(\frac{(N-1)\lambda}{N}\right)
\\ &\qquad \times \left[\left( \psi_I\left(\frac{(N-1)\lambda}{N}\right) \right) -
\psi_I\left( \lambda\right)  \right] .
\end{align*}
In order to compute
$p_k^{(N)}$ using (\ref{fin-size}) it is required to sequentially compute $p_0^{(N)}$ up to
$p_{k-1}^{(N)}$. Further, the formula is not very
enlightening and it may be numerically very unstable when $k$ (and
hence $N\ge k$) is large. For this reason we
devote the major part of these notes to approximations assuming $N$ is large.

In Section 1.9 of Part II of this volume the exact results above are generalized to a model allowing for heterogeneous spreading, meaning that the transmission rate depends on the two individuals involved.

\begin{exercise} \label{x3.2}
Compute the final size distribution $\{p_k^{(N)}\}$ numerically using some suitable software for $N=10,\ 50$ and 100, for $\lambda=2$ and $I\equiv 1$ (the Reed--Frost model) and $I\sim \Gamma (3, 1/3)$ (having mean 1 and variance 1/3).
\end{exercise}

\section{The Sellke construction}\label{sec-Sellke}
We now present the Sellke construction (Sellke \cite{S1983I}), which is an ingenious way to define the epidemic outbreak in continuous time using two sets of i.i.d.\ random variables. This elegant construction is made use of in many new epidemic models, as proven by having more than 50 citations in the past decade.

We number the individuals from $0$ to $N$:
 $  0\ 1\ 2\ 3\ \ldots\ N$.
 Index $0$ denotes the initially infected individual, and the individuals numbered from $1$ to $N$ are all
  susceptible at time $0$.

 Let
 \begin{description}
 \item $Q_1,Q_2,\ldots,Q_N$ be i.i.d.\ random variables, with the law $\mathrm{Exp}(1)$;
 \item $(L_0,I_0), (L_1,I_1),\ldots,(L_N,I_N)$ be i.i.d.\ random variables, with the law $\P_{(L,I)}$.

 In the Markov model, $L_i$ and $I_i$ are independent, hence
 $\P_{(L,I)}=\P_L\otimes\P_I$,\footnote{This notation stands for the product of the two probability measures $\P_L$ and $\P_I$. The fact that the law of the pair is the product of the two marginals is equivalent to the fact that the two random variables $L$ and $I$ are independent.}  where $\P_L$ is the law of the
  latency period and $\P_I$ that of the infectious period. But this need not be the case in more general non-Markov models.
 \end{description}
 Individual $0$ has the latency period $L_0$ and the infectious period $I_0$.
 We denote below
 \begin{align*}
 L(t)&\text{  the number of individuals in state $E$ at time $t$};\\
 I(t)&\text{  the number of individuals in state $I$ at time $t$}.
 \end{align*}
 Note that for each $i$, the two random variables $L_i$ and $I_i$ could be dependent, which typically is not the case in a Markov model.

 We define the cumulative force of infection experienced by an individual, between times $0$ and $t$ as
 \[  \Lambda_C(t)=\frac{\lambda}{N}\int_0^t I(s)ds.\]
 For $i=1,\ldots,N$, individual $i$ is infected at the time when $\Lambda_C(t)$ achieves the
 value $Q_i$ (which might be considered as the ``level of resistance to infection of individual $i$'').
 The $j$-th infected susceptible has the latency period $L_j$ and the infectious period $I_j$.
 The epidemic stops when there is no individual in either the latent or infectious state, after which $\Lambda_C(t)$
 does not grow any more, $\Lambda_C(t)=\Lambda_C(\infty)$.
 The individuals such that $Q_i>\Lambda_C(\infty)$ escape infection.

 We put the $Q_i$s in increasing order: $Q_{(1)}<Q_{(2)}<\cdots<Q_{(N)}$. It is the order in which individuals are infected in Sellke's model. Note that Sellke's model  respects the durations of latency and infection.
 In order to show that Sellke's construction gives a process which has the same law as the process from Definition \ref{def:SEIR}, it remains to verify that the rates at which infections happen are the correct ones.

In the initial model, we assume that each infectious meets other individuals at rate $c$. Since  each individual
 has the same probability of being the one who is met, the probability that a given individual is that one is $1/N$. Hence the rate at which a given individual is met by a given infectious one is $c/N$. Each encounter between a susceptible and an infectious individual achieves an infection with probability $p$. Hence the rate at which a given individual is infected by a given infectious individual is $\lambda/N$, where we have set $\lambda=cp$. The rate at which an  infectious individual infects  susceptibles is then $\lambda S(t)/N$.
 Finally the epidemic  propagates at rate $\lambda S(t)I(t)/N$.

 Let us go back to Sellke's construction. At time $t$, $S(t)$ susceptibles have not yet been infected.
 Each of those corresponds to a $Q_i>\Lambda_C(t)$. At time $t$, the slope of the curve which represents
 the function $t\mapsto\Lambda_C(t)$ is $\lambda I(t)/N$.
 If $Q_i>\Lambda_C(t)=x$, then
 \pagebreak
  \begin{align*}
 \P(Q_i>x+y|Q_i>x)&=e^{-y},\\
\text{hence }\ \P(Q_i>\Lambda_C(t+s)|Q_i>\Lambda_C(t))&=\exp\left(-\frac{\lambda}{N}\int_t^{t+s}I(r)dr\right)\\
 &=\exp\left(-\frac{\lambda}{N}I(t)s\right),
 \end{align*}
 if $I$ is constant on the interval $[t,t+s]$.
 Consequently, conditionally upon $Q_i>\Lambda_C(t)$,
 \[  Q_i-\Lambda_C(t)\sim \mathrm{Exp}\left(\frac{\lambda}{N}I(t)\right).\]
 The same is true for all $S(t)$ of those $Q_i$ which are $>\Lambda_C(t)$. The next individual to get infected corresponds to the minimum of those $Q_i$, hence the waiting time after  $t$ for the next infection follows the law $\mathrm{Exp}\left(\frac{\lambda}{N}I(t)S(t)\right)$,
 if no removal of an infectious individual happens in the mean time, which would modify $I(t)$.

 Thus in Sellke's construction, at time $t$ the next infection comes at rate
 \[  \frac{\lambda}{N}I(t)S(t),\]
 as in the model described above.

\section{LLN and CLT for the final size of the epidemic}\label{TB-EP_sec_LLNCLT_final_size}
Define, for $0\le w\le N+1$, with the notation $[w]$ = integer part of $w$,
and the convention that a sum over an empty index set is zero,
\[  \J(w)=\frac{\lambda}{N}\sum_{i=0}^{[w]-1}I_{i}.\]
Note that $i=0$ is the index of the initially infected individual, $I_{i}$ denotes here the length of the infectious period of individual  whose resistance level is $Q_{(i)}$ (who is not that of the $i$-th individual of the original list, but of the individual having the $i$-th smallest resistance).

$\J(w)$ is the infection pressure produced by the first $[w]$ infected individuals (including number $0$).
For any integer $k$, $\J$ is of course constant on the interval $[k,k+1)$. Define for $v>0$ the number of individuals who do not
resist to the infectious pressure $v$:
\[  \q(v)=\sum_{i=1}^N {\bf1}_{\{Q_i\le v\}}.\]
The total number of infected individuals in the epidemic is
\begin{align}\label{tailletot}
Z&=\min\left\{k\ge0;\ Q_{(k+1)}>\frac{ \lambda}{N}\sum_{i=0}^kI_i\right\}\\
&=\min\left\{k\ge0;\ Q_{(k+1)}>\J(k+1)\right\} \nonumber\\
&=\min\left\{w\ge0;\ \q(\J(w+1))=w\right\}. \nonumber
\end{align}
Suppose indeed that $Z=i$. Then according to \eqref{tailletot},
\begin{align*}
\J(j)&>Q_{(j)},\ \text{ hence } \q(\J(j))\ge j,\quad \text{for all } j\le i,\\
\text{and }\J(i+1)&<Q_{(i+1)}\ \text{ hence } \q(\J(i+1))<i+1.
\end{align*}
In other words  $Z=i$ if and only if $i$ is the smallest integer such that
\[
\q(\J(i+1))<i+1,\ \text{ hence }\q(\J(i+1))=i.\]

\subsection{Law of Large Numbers} \label{subsec_LLN}
Let us index $\J$ and $\q$ by $N$, the population size, so that they become $\J_N$ and $\q_N$.
We now define
\begin{align*}
\overline{\J}_N(w)&=\J_N(N w)\\
\overline{\q}_N(v)&=\frac{\q_N(v)}{N}.
\end{align*}
It follows from the strong law of large numbers that as $N\to\infty$,
\begin{align*}
\overline{\J}_N(w)&\to \lambda\E(I) w=R_0 w\ \text{ almost surely, and}\\
\overline{\q}_N(v)&\to 1-e^{-v}\ \text{ a.s}.
\end{align*}
Hence, with the notation $f\circ g(u):=f(g(u))$, as $N\to\infty$,
\[  \overline{\q}_N\circ \overline{\J}_N(w)\to 1-e^{-R_0 w}\]
a.s., uniformly on $[0,1]$ (the uniformity in $w$ follows from the second Dini theorem, as in the proof of Proposition \ref{lgn_input}).
We have (replacing now $Z$ by $Z^N$)
\begin{align*}
\frac{Z^N}{N}&=\min\left\{\frac{w}{N}\ge0;\ \q_N(\J_N(w+1))=w\right\}\\
&=\min\left\{ s\ge0;\ \frac{1}{N}\q_N\left(\J_N\left(N\left(s+\frac{1}{N}\right)\right)\right)=s \right\}\\
&=\min\left\{ s\ge0;\  \overline{\q}_N\left(\overline{\J}_N\left(s+\frac{1}{N}\right)\right)=s\right\}.
\end{align*}
Recall from \eqref{R_0} that $R_0=\lambda\iota$, where $\iota=\E(I)$.
Note that when $R_0>1$, the equation
\begin{equation}\label{eq:LLNfinal-size}
z=1-e^{-R_0 z}
\end{equation}
 (which is equation \eqref{final-size} from Section \ref{Sec_Det-mod}) has a unique solution $z^\ast\in(0,1)$ (besides the zero solution). Indeed, $f(z)=1-e^{-R_0 z}$ is concave, $f(1)<1$, and $f'(0)=R_0$.
For the proof of the next theorem, we follow an argument from Andersson and Britton \cite{AB00I} (see also Ball and Clancy \cite{BC93I}).

\begin{theorem}\label{th:LN_finsiz}
If $R_0\le1$, then  $Z^N/N\to0$ a.s., as $N\to\infty$.

If $R_0>1$,
as $N\to\infty$,  $Z^N/N$ converges in law to the random variable $\zeta$ which is such that
$\P(\zeta=0)=z_\infty =1-\P(\zeta=z^\ast)$, where

$z_\infty$, the probability of a minor outbreak (i.e.\  that the epidemic does not get off), is the solution in $(0,1)$ of \eqref{eq:P-min-epid} below, and $z^\ast$
is the positive solution of \eqref{eq:LLNfinal-size}.
\end{theorem}
Let us explain how one can characterize $z_\infty$. It follows from Theorem \ref{Th-br-epid-coupling}
 that the probability $z_\infty$ that the epidemic does not get off equals the probability that the associated branching process goes extinct, which is the probability that the associated discrete time branching process (where we consider the infected by generation) goes extinct.
 According to Proposition \ref{BGW} from Appendix A, the probability that this happens is the solution in the interval $(0,1)$ of the equation
 $g(s)=s$, where $g$ is the generating function of the random number $\xi$ of individuals that one infected infects. As explained in Section \ref{sec-early-stage}, the law of $\xi$ is MixPoi$(\lambda I)$, so if we denote by $\psi_I(\mu)=\E[\exp(-\mu I)]$ the Laplace transform
  of $I$, which is well defined for $\mu>0$, then $g(s)=\psi_I(\lambda(1-s))$. Hence $z_\infty$ is the unique solution in $(0,1)$ of the equation
  \begin{equation}\label{eq:P-min-epid}
  \psi_I(\lambda(1-s)) = s\,.
  \end{equation}

  \bpf
 If $R_0\le1$, then from Corollary \ref{Th-early-stage-sub},
$Z^N$ remains bounded, hence $Z^N/N\to0$.

If $R_0>1$, then $Z^N$ remains bounded with probability $z_\infty$, which is the probability of extinction in the branching process which approximates the early stage of the epidemic. We now need to see what happens on the complementary event. For that sake, we first  choose an arbitrary sequence of integers $t_N$, which satisfies both $t_N/N\to0$ and $t_N/\sqrt{N}\to\infty$, as $N\to\infty$. We note that
on the event $\{Z^N\le t_N\}$, each infective infects susceptibles at a rate which is bounded below by $\lambda_N=\lambda\frac{N+1-t_N}{N}$.  Let $Z(\lambda_N,I)$ denote the total progeny of a single ancestor in a branching process, where each individual has children
according to a rate $\lambda_N$ Poisson process, during his life whose length is $I$. It is plain that for ant $t\in\Z_+$, and $N$ large enough such that $t\le t_N$,
\[ \P(B(\lambda,I)\le t)\le \P(Z^N\le t)\le \P(Z^N\le t_N)\le \P(B(\lambda_N,I)<\infty).\]
Define as in the statement $z_\infty=\P(B(\lambda,I)<\infty)$ the probability of extinction of the branching process approximating the early stage of the epidemic, and $z_{N,\infty}=\P(B(\lambda_N,I)<\infty)$. It is not hard to show that $z_{N,\infty}\to z_\infty$ as $N\to\infty$, as a consequence of the fact that $\lambda_N\to\lambda$ (since $t_N/N\to0$). Hence for any $\eps>0$, we can choose $t$ large enough such that $\P(B(\lambda,I)\le t)\le z_\infty-\eps$, and $N$ large enough such that $z_{N,\infty}\le z_\infty+\eps$. We have shown that
\begin{equation}\label{tN}
\P(Z^N\le t_N) \to z_\infty,\ \text{ as }N\to\infty.
\end{equation}
This shows that a.s. on the event that the epidemic goes off, $Z^N$ tends to $\infty$ faster than $t_N$. We will next prove that
\begin{equation}\label{Prob_deviat}
\lim_{c\to\infty}\lim_{N\to\infty}\P\left(\left\{t_N< Z^N<Nz^\ast-c\sqrt{N}\right\}\bigcup\left\{Z^N>Nz^\ast+c\sqrt{N}\right\}\right)=0\, .
\end{equation}
Recalling the last formula preceding the statement of the present theorem,
\begin{align}
&\left\{\frac{Z^N}{N}\in\left(t_N,z^\ast-\frac{c}{\sqrt{N}}\right)\bigcup\left(z^\ast-\frac{c}{\sqrt{N}},1\right]\right\}\nonumber\\
&\quad\subset \left\{\exists s\in \left(t_N,z^\ast-\frac{c}{\sqrt{N}}\right)\bigcup\left(z^\ast-\frac{c}{\sqrt{N}},1\right];\
\overline{\q}_N\left(\overline{\J}_N\left(s+\frac{1}{N}\right)\right)=s\right\}\nonumber\\
&\quad\subset \left\{\sup_{0\le s\le 1}\left| \overline{\q}_N\left(\overline{\J}_N\left(s+\frac{1}{N}\right)\right)-1+e^{-R_0s}\right|>\frac{\phi(c)}{\sqrt{N}}\right\},\label{eventN}
\end{align}
where $\phi(c)\to\infty$, as $c\to\infty$, for $N$ large enough. We have exploited the facts that $t_N/\sqrt{N}\to\infty$ as $N\to\infty$, and
$f'(0)>1$, $f'(z^\ast)<1$.
However, we shall see in the next subsection (see \eqref{ComFunctCV}) that
\[ \left\{\sqrt{N}\left(\overline{\q}_N\left(\overline{\J}_N\left(s+\frac{1}{N}\right)\right)-1+e^{-R_0s}\right),\, s\in[0,1]\right\}\]
converges weakly, for the sup--norm topology, to a centred Gaussian process with finite covariance, hence the limit as $N\to\infty$ of the probability of the event
\eqref{eventN} tends to $0$, as $c\to\infty$, which establishes \eqref{Prob_deviat}. It is easily seen that the second part of the Theorem follows from the combination of \eqref{tN} and \eqref{Prob_deviat}.
\epf

We see that $z^*$ is the size, measured as the proportion of the total population,
of a ``significant'' epidemic, if it takes off, which happens with probability $1-z_\infty$.

We notice that $z^*$ depends on the particular model only through the quantity $R_0$. In particular it depends on the law of the infectious period $I$ only through its mean. In the case where both $E$ and $I$ are exponential random variables, we know from Section \ref{TB-EP_sec_LLN}
that the model has a law of large numbers limit, which is a system of ODEs. The same value for $z^*$ has been deduced from an analysis of this deterministic model in Section \ref{Sec_Det-mod}. The last theorem holds for a larger class of models.
\pagebreak

\subsection{Central Limit Theorem}

From the classical CLT, as $N\to\infty$,
\begin{align*}
A_N(\omega):=\sqrt{N}(\overline{\J}_N(w)-R_0 w)&=\frac{\lambda\sqrt{w}}{\sqrt{Nw}}\sum_{i=0}^{[Nw]}[I_i-\E(I_i)] + O(1/\sqrt{N})\\
&\Rightarrow A(w),
\end{align*}
where $A(w)\sim \NN(0,p^2c^2\text{Var}(I)w)$.
One can in fact show that, as processes
\[ \{\sqrt{N}(\overline{\J}_N(w)-R_0 w),\ 0\le w\le 1\}\Rightarrow\{A(w),\ 0\le w\le 1\}\]
for the topology of uniform convergence,
where $\{A(w),\ 0\le w\le 1\}$ is a Brownian motion
(i.e.\  a centered Gaussian process with independent increments and continuous trajectories)
such that Var$(A(w))=r^2R_0^2w$, where $r^2=(\E I)^{-2}\text{Var}(I)$.
It is easy to show that for all $k\ge1$, all
$0<w_1<\dots<w_k\le 1$, if we define $A_N(w):=\sqrt{N}(\overline{\J}_N(w)-R_0 w)$,
\[ (A_N(w_1),\ldots,A_N(w_k))\Rightarrow(A(w_1),\ldots,A(w_k)).\]
This means the convergence of the finite dimensional distributions. Combining this with the techniques exposed in Section \ref{tight} of the Appendix yields the above functional weak convergence.

Consider now $\overline\q_N$. Again from the usual CLT,
\begin{align*}
B_N(v)&=\sqrt{N}(\overline{\q}_N(v)-[1-e^{-v}])\\
&=\frac{1}{\sqrt{N}}\sum_{i=1}^N\left[{\bf1}_{\{Q_i\le v\}}-(1-e^{-v})\right]\\
&\Rightarrow B(v),
\end{align*}
where $B(v)\sim \NN(0,e^{-v}(1-e^{-v}))$. We have again a functional convergence, according to the
Kolmogorov--Smirnov theorem, towards a time changed Brownian bridge. In simpler words,
$\{B(v),\ v\ge0\}$ is a centred Gaussian process with continuous trajectories
whose covariance function is specified by the identity
$\E[B(u)B(v)]=e^{-u\vee v}-e^{-(u+v)}$, where $u\vee v:=\sup(u,v)$.

Let us now combine the two functional central limit theorems which we have just derived.
We have

\begin{align*}
\sqrt{N}&\left(\overline{\q}_N(\overline{\J}_N(w))-1+e^{-R_0w}\right)\\
&\qquad=\!\sqrt{N}\!\left(\overline{\q}_N(\overline{\J}_N(w))-1+\exp(-\overline{\J}_N(w))\right)
\!+\!\sqrt{N}\!\left(e^{-R_0 w}-e^{-\overline{\J}_N(w))}\right)\\
&\qquad\sim B_N(\overline{\J}_N(w))-R_0e^{-R_0w}A_N(w).
\end{align*}
Consequently
\begin{equation}\label{ComFunctCV}
\sqrt{N}\left(\overline{\q}_N(\overline{\J}_N(w))-1+e^{-R_0w}\right) B(R_0w)-R_0e^{-R_0w}A(w),
\end{equation}
which is the
functional central limit theorem which was used in the proof of Theorem \ref{th:LN_finsiz}.

Recall that the above Law of Large Numbers has been obtained by taking the limit
in the equation
\[ \overline{\q}_N\left(\overline{\J}_N\left(z+N^{-1}\right)\right)=z.\]
Making use of the above two CLTs, we get
\begin{align*}
z&=1-e^{-\overline{\J}_N\left(z+N^{-1}\right)}+N^{-1/2}B_N(\overline{\J}_N(z+N^{-1}))\\
&=1-\exp\left(-R_0(z+N^{-1})-N^{-1/2}A_N(z+N^{-1})\right)\\
&\quad +N^{-1/2}B_N\left(R_0(z+N^{-1})+N^{-1/2}A_N(z+N^{-1})\right).
\end{align*}
Let $z=z^\ast+z_NN^{-1/2}+\circ(N^{-1/2})$, where $z^\ast$ satisfies
$e^{-R_0z^\ast}=1-z^\ast$.
We obtain
\begin{align*}
&z^\ast+z_NN^{-1/2}+\circ(N^{-1/2})\\&=
1-\exp\left(-R_0z^\ast-R_0z_NN^{-1/2}-A_N(z^\ast)N^{-1/2}+\circ(N^{-1/2})\right)\\
&\quad+N^{-1/2}B_N(R_0 z^\ast)+\circ(N^{-1/2})\\
&=1-e^{-R_0z^\ast}+N^{-1/2}e^{-R_0z^\ast}\left(R_0z_N+A_N(z^\ast)\right)
+N^{-1/2}B_N(R_0z^\ast)+\circ(N^{-1/2}).
\end{align*}
We simplify this relation by making use of the equation which specifies $z^\ast$. Multiplying the
remaining terms by $N^{1/2}$, we deduce
\[ [1-(1-z^{\ast})R_0]z_N=B_N(R_0 z^\ast)+(1-z^\ast) A_N(z^\ast)+\circ(1).\]
Hence $z_N\Rightarrow \Xi$, where (note that $e^{-R_0z^\ast}(1-e^{-R_0z^\ast})=z^\ast(1-z^\ast)$)
\[\Xi\sim \NN\left(0,\frac{z^\ast(1-z^\ast)}{(1-(1-z^\ast)R_0)^2}\left(1+r^2(1-z^\ast)R_0^2\right)\right),\]
where we have exploited the independence of the two processes $A(\cdot)$ and $B(\cdot)$, which follows from that of the two collections of random variables $(I_i,\ i\ge0)$ and $(Q_i,\ i\ge1)$.

Finally we can conclude with the following theorem. We refer to Scalia-Tomba \cite{ST85I} and \cite{ST90I} for a more complete justification.
\begin{theorem}\label{CLT_final_size}
 As $N\to\infty$, conditionally upon the event that the epidemic takes off, the law of $N^{-1/2}(ZN-Nz^\ast)$ converges towards the Gaussian distribution
\[\NN\left(0,\frac{z^\ast(1-z^\ast)}{(1-(1-z^\ast)R_0)^2}\left(1+r^2(1-z^\ast)R_0^2\right)\right).\]
\end{theorem}
\begin{exercise} \label{xCLT1}
Compute numerically the limiting mean and standard deviation of the final size $Z^N$ in case of a major outbreak and $N=1000$, $\lambda=1.5$ and $\iota=1$, for the following two situations. The first scenario is when $I\equiv 1$ (fixed infectious period), and the second when $I\sim \mathrm{Exp}(1)$ (Markovian SIR).
\end{exercise}

\section{The duration of the stochastic SEIR epidemic}\label{sec:duration-epid}

Recall that $L^N(t)$ and $I^N(t)$ denote the numbers of latent and infectious individuals at time $t$ respectively, and introduce $Z^N(t)=N-L^N(t)-I^N(t)-R^N(t)$ to denote the number of individuals who have been infected by time $t$ (i.e.\ who are no longer susceptible). We now study how long it takes for the epidemic to first grow big, and then later to end, i.e.\ for the end of the epidemic we will study properties of $\tau^N=\inf \{  t; L^N(t)+I^N(t)=0\}$ as $N\to\infty$. It will only be a sketch since it is quite technical to prove the results rigorously. For detailed results we refer to Barbour \cite{B75I}. From an applied point of view, this question has clear practical relevance, since for instance hospitals are on highest pressure when the epidemic peaks, and knowing how long until the outbreak is over indicates how long preventive measures should be enforced.

If the epidemic does not take off we know from branching process theory that the time to extinction is finite, so $\tau^N=O_p(1)$ on this part of the sample space  ($O_p(1)$ denotes bounded in probability). We hence focus on the situation where the epidemic takes off resulting in a major outbreak, hence implicitly assuming that $R_0>1$.

We divide the duration of the whole epidemic $\tau^N$ into three parts: the beginning, the main part and the end of the epidemic. Pick $\epsilon>0$ small. Formally we define these parts by defining two intermediate times (inspired by Sir Winston Churchill): the end of the beginning $\tau_{Beg}^N=\inf \{ t\le \tau^N; Z^N(t)\ge \epsilon N\}$, and the beginning of the end $\tau_{End}^N=\inf \{ t\le \tau^N; Z^N(t)\ge (1-\epsilon) z^\ast N\}$, where $z^\ast$ is the positive solution to the final size equation from Section \ref{Sec_Det-mod}. Each of these times are equal to $\tau^N$ in the case when the event never occurs.

With these definitions the beginning of the epidemic is the time interval $[ 0, \tau_{Beg}^N )$, the main part $[ \tau_{Beg}^N, \tau_{End}^N )$ and the end part $[ \tau_{End}^N, \tau^N ]$.

During the beginning we can sandwich the epidemic between two branching processes. The upper bound is the branching process $Z(t)$ described in Section \ref{sec-early-stage}. Similarly, we can construct a lower bound using a very similar branching process $Z^-(t)$, the only difference being that the birth rate is $\lambda (1-\epsilon)$ as opposed to $\lambda$ for the upper branching process. This is true because before $\tau_{Beg}^N$ the rate of new infections in the epidemic equals $\lambda (1-Z^N(t)/N)$ which lies between $\lambda (1-\epsilon)$ and $\lambda$. Since $Z^{-}(t)\le Z^N(t)\le Z(t)$ for $t\le \tau_{Beg}^N$ it follows that $\tau_{+}^N\le \tau_{Beg}^N\le \tau_{-}^N$, where $\tau_{+}^N = \inf \{ t; Z(t)\ge \epsilon N\}$ and $\tau_{-}^N = \inf \{ t; Z^-(t)\ge \epsilon N\}$.

From Section \ref{Sec_Cont_Br_Pr} we know the rate at which a branching process grows. More specifically, we know that when a branching process $Z'(t)$ takes off, it grows exponentially: $Z'(t)\sim e^{r't}$, where $r'$ is the unique solution to $1=\int_0^\infty e^{-r's}\lambda (s)ds =1$, where $\lambda (s)$ is the average (expected) rate at which an individual gives birth at age $s$ (cf.\ Equation (\ref{Malthus-eq})). For our two branching processes $Z(t)$ and $Z^-(t)$ we have $\lambda (s)= \lambda \P (L<s<L+I)$ and $\lambda^-(s)= \lambda (1-\epsilon) \P (L<s<L+I)$ respectively. From this it follows that the exponential growth rates $r$ and $r^-$ can be made arbitrary close to each other by choosing $\epsilon$ small enough ($r^-=r(1+o(\epsilon ))$). The particular form of $r$ and $r^-$ depends on the distribution of $L$ and $I$ (see Exercise \ref{xDur.1} below). Recall that $\tau_{+}^N = \inf \{ t; Z(t)\ge \epsilon N\}$, so the fact that $Z(t)\sim e^{rt}$ implies that $\tau_{+}^N = \frac{\log (\epsilon N)}{r} + O_p(1)$. Similarly, $\tau_{-}^N = \frac{\log (\epsilon N)}{r^-} + O_p(1)$. As a consequence, the two stopping times are arbitrary close to each other  on the logarithmic scale. From this we have $\tau_{Beg}^N= \frac{\log (N)}{r} (1+o(\epsilon))+ O_p(1)$.

We now turn to the duration of the main part of the epidemic: $\tau_{End}^N-\tau_{Beg}^N$ which is positive only if the epidemic takes off, which we hence condition upon. During this part of the epidemic, the Markovian SEIR epidemic can be approximated by the deterministic SEIR model. This means that for the Markovian SEIR model, the duration of the main part of the epidemic $\tau_{End}^N-\tau_{Beg}^N$ can be well approximated by the corresponding duration of the deterministic system $\tau_{End}^{Det}-\tau_{Beg}^{Det}$. The deterministic system is started at $\tau^{Det}_{Beg}=0$ with initial conditions  $(s(0), e(0), i(0),r(0))=(1-\epsilon, a\epsilon, b\epsilon, (1-a-b)\epsilon)$ for some positive numbers $a$ and $b$ with $0<a+b\le 1$ (there is no closed form expression for how the infected individuals are divided into exposed, infectives and recovereds). The system is then run until $\tau^{Det}_{End}=\inf \{ t; e(t)+i(t)+r(t)\ge (1-\epsilon)z\}$. We know that $z(t)=e(t)+i(t)+r(t)\to z$ and $z(t)$ is monotonically increasing (since $s(t)=1-z(t)$ is decreasing). This implies that $\tau^{Det}_{End}-\tau^{Det}_{Beg}=\tau^{Det}_{Beg}$ is just a constant for any fixed positive $\epsilon$. It will depend slightly on $a$ and $b$, but when $\epsilon$ is small the dependence is weak and there is a uniform bound. From this we conclude that the main part of the epidemic is bounded:
$$
\tau_{End}^N-\tau_{Beg}^N= \tau^{Det}_{Beg} +o_p(\epsilon)=O_p(1).
$$
 If the latent and infectious periods are not exponentially distributed, then the stochastic SEIR epidemic is not Markovian, and the deterministic approximating system is a difference-delay-system which we will not study more closely. The qualitative properties of this system coincide with those of the Markovian SEIR system; in particular, the duration of the main part is bounded in probability.

Just like the main part of the epidemic the duration of end of the epidemic, $\tau^N-\tau^N_{End}$ is only positive if the epidemic takes off, which we hence condition upon. At the beginning of the end part, the number of infected (either exposed, infectious or recovered) equals $Z^N(\tau^N_{End})=(1-\epsilon)z^\ast N$ and $S^N(\tau^N_{End})=(1-z^\ast )N +\epsilon z^\ast N$. Since $\epsilon$ is assumed to be small, infectious individuals give birth at rate $\lambda (1-z^\ast  +\epsilon z^\ast )\approx (1-z^\ast )$ during the rest of the epidemic (we know the final fraction infected converges to $z^\ast $ in probability). Further, at the start of the beginning the fractions exposed and infectious will both close to that of the deterministic system which are both small, having size $c_E\epsilon$ and $c_I \epsilon$ say (cf.\ Figure \ref{fig-det-syst} where it is seen that $e(t)$ and $i(t)$ are both small for large $t$). So, from the beginning of the end part, the epidemic behaves like a branching process with childhood duration $L$, adult duration $I$ and birth rate $\lambda (1-z^\ast )$ during the adult life stage, and this part is started with $c_E\epsilon N$ children (exposed) and $c_I\epsilon N$ adults (infectious). The mean off-spring distribution for this branching process equals $\lambda \E (I) (1-z^\ast )=R_0(1-z)$ where $z^\ast $ is the positive solution to $1-z=e^{R_0z}$. It can be shown (cf.\ Exercise \ref{xSubcrit} below) that $R_0(1-z^\ast )<1$ implying that the branching process is subcritical (otherwise the epidemic would not be on decline).

The duration $\tau^N-\tau^N_{End}$ of the end part can hence be approximated by the time until extinction of a subcritical branching process, starting with $c_E\epsilon N$ children\linebreak (exposed) and $c_I\epsilon N$ adults (infectious). This branching process will have \emph{negative} drift $r^*<0$ being the solution to the corresponding equation $\int_0^\infty e^{-rs}\lambda (s)ds=1$ where now $\lambda (s)=\lambda (1-z^\ast )\P (L<s<L+I)$. So, $E(t)+I(t)\sim (E(0)+I(0))e^{r^*t}=(c_E+c_I)\epsilon N e^{r^*t}$. The time until this branching process goes extinct (i.e.\ $E(t)+I(t)<1$) is hence of order $-\log ((c_E+c_I)\epsilon N)/r^*=-\log N/r^*+O_p(1)$.

To sum up, the duration of the epidemic $\tau^N=O_p(1)$ if the epidemic does not take off, whereas it has the following structure in case it does take off:
\begin{align}
\tau^N=\tau^N_{Beg}+\left( \tau^N_{End}- \tau^N_{Beg}\right) + \left(\tau^N- \tau^N_{End} \right) = \frac{\log N}{r} + O_p(1) + \frac{-\log N}{r^*}.
\end{align}
Note that the last term is also positive since $r^*<0$.

\begin{exercise} \label{xSubcrit}
Show that $R_0(1-z^\ast )<1$ and compute it numerically for $R_0=1.5$.
\end{exercise}

\begin{exercise} \label{xDur.1}
Consider the stochastic SEIR epidemic with infection rate $\lambda=1.5$ per time unit. Compute the two leading terms of the duration of a major outbreak for the following three case: $L\equiv 0$ and $I\sim \mathrm{Exp} (1)$ (Markovian SIR), $L\equiv 0$ and $I\equiv 1$ (continuous time Reed--Frost), and $L\sim \mathrm{Exp}(1)$ and $I\sim \mathrm{Exp} (1)$ (Markovian SEIR).
\end{exercise}

\chapter{Open Markov Models}\label{TB-EP_chap_OpenMarkovMod}

In this chapter, contrary to the situation considered in earlier chapters, we study models where
there is a constant supply of susceptibles (either by births, immigration or loss of immunity of the removed individuals) giving rise to endemic type situations.
We study how the random fluctuations in the model can drive the system out of the basin
of attraction of the stable endemic equilibrium of the deterministic model, such that the disease goes extinct.

As we shall see in Section \ref{Sec_Open-pop},  in the case of a moderate population size, one may expect that the Gaussian fluctuations described by the central limit theorem are strong enough to stop the endemy in a SIR model with demography.
For larger population sizes, following Freidlin and Wentzell \cite{FWI}, we describe in Section \ref{TB-EP_sec_LD} how long it will take for the random perturbations to stop the endemy. We apply this approach successively to the SIRS, the SIS and the SIR model with demography. In the case of the SIS model, we  compute explicitly the constant which appears in the Freidlin--Wentzell theory, see Proposition \ref{explicitV} below. This is unfortunately the only case where we have such a simple and explicit formula in terms of the coefficients of the model.

\section[Open populations: time to extinction and critical population size]{Open populations: time to extinction and critical population size}\label{Sec_Open-pop}

Up until now we have (mainly) considered the stochastic SEIR epidemic model in a fixed community of size $N$, where $N$ has been assumed large (except in Section \ref{Sec_Exact} when $N$ was assumed small). This is of course an approximation of reality, but when considering outbreaks of a few months (e.g.\ influenza outbreaks) it seems like a fair approximation; recall that the time to extinction of our model was $O_p(\log N)$. For other diseases including childhood diseases, the disease is present in the community constantly -- such diseases are said to be \emph{endemic}. When trying to understand the behaviour of such diseases it is necessary to also allow people to die and new people entering the population (by birth or immigration). In the current section we do this and derive approximations for two important quantities: the time to extinction of the endemic disease $T_E^N$, and the critical community size $N_c$. These two quantities have received much attention in the literature over the years. In particular, the critical community size $N_c$ and how it depends on properties of the disease and the community have been studied both in the mathematical and applied communities (e.g.\ Lindholm and Britton \cite{LB07I} and Keeling and Grenfell \cite{KG97I}).

Let us first describe the population model, which is the simplest model for a population which fluctuates randomly in time with a mean size of $N$ individuals, and where individuals have life time distribution with mean $1/\mu$ (cf.\ Example \ref{ex:SEIRS}). The population $N(t)$ is defined to be a Markovian birth-death process with constant birth rate $\mu N$ and linear death rate, all individuals dying at rate $\mu$. This process $N(t)$ will fluctuate around $N$, a parameter we denote by the mean population size. If $N$ is large, it is known that $N(t)$ will be approximately normally distributed with mean $N$ and standard deviation proportional to $\sqrt{N}$, so for practical purposes we will later approximate $N(t)$ by $N$.

\begin{exercise}\label{exer4.1.1}
Assuming that $N(0)=N$, write $N(t)$ as the solution of an SDE of the same form
as the SDE appearing at the beginning of Section \ref{TB-EP_sec_LLN}. Define $Q^N_t=N(t)/N$ and show that, as a consequence of Theorem \ref{th:LLN}, $Q^N_t\to1$ a.s., locally uniformly in $t$.
Then deduce from Theorem \ref{th:CLT} that $\sqrt{N}(Q^N_t-1)$ converges weakly, as $N\to\infty$, towards an Ornstein--Uhlenbeck process of the form
\[U_t=\sqrt{2\mu}\int_0^te^{-\mu(t-s)}dB_s,\]
where $B_t$ is a standard Brownian motion. Prove that $\E(U_t)=0$ and Var$(U_t)\to1$ as
$t\to\infty$. Deduce that for large $N$ and $t$, $N(t)$ is approximately normally distributed with mean $N$ and standard deviation proportional to $\sqrt{N}$.
\end{exercise}

For this population model, we assume that the Markovian SIR epidemic spreads (this can easily be extended to the Markovian SEIR model). By this we mean that individuals who get infected immediately become infectious and remain so for an $\mathrm{Exp}(\gamma)$ time, unless they happen to die before by chance. In the fixed population size model, the contact rate was $\lambda$ which implied that it was $\lambda/N$ to each specific individual. Now, in the open population model, we assume that the infection rate to a specific individual is unchanged, $\lambda/N$. More appropriate would perhaps have been to instead have  $\lambda/N(t)$ but since $N(t)\approx N$ for all $t$ we use the simpler choice $\lambda/N$.  So, new individuals enter the community at constant rate $\mu N$ and all individuals die, irrespective of being susceptible, infectious or recovered, at rate $\mu$, susceptible individuals get infected at rate $\lambda I^N(t)/N$, and infectious individuals recover at rate $\gamma$. The rate at which susceptibles get infected and infected recover hence equals $\lambda I^N(t)S^N(t)/N$, and $\gamma I^N(t)$ respectively. If we study the limiting deterministic system for the \emph{fractions} in each state we get the following system of differential equations:
\pagebreak
\begin{align}
s'(t) &=\mu -\lambda s(t)i(t) - \mu s(t),\nonumber
\\
i'(t)&= \lambda s(t)i(t) - \gamma i(t)-\mu i(t), \label{det-SIR-demog}
\\
r'(t) &= \gamma i(t)-\mu r(t),\nonumber
\end{align}
which is identical to those of Example \ref{Ex_SEIRS_3} with $\rho=0$ and $\nu=+\infty$.
From this we can compute the endemic state where all derivatives are 0. First we note that the basic reproduction number $R_0$ (the expected number of infectious contact while infectious and alive) and the expected relative time of a life an individual is infected, $\varepsilon$, are given by
\begin{equation}
R_0=\frac{\lambda}{\gamma + \mu} \hskip2cm \varepsilon =\frac{1/(\gamma + \mu)}{1/\mu}=\frac{\mu}{\gamma + \mu}.
\end{equation}
The rate of recovery $\gamma$ is much larger than the death rate $\mu$ (52 compared to 1/75 for a one week infectious period and 75 year life length) so for all practical purposes the two expressions can be approximated by $R_0\approx \lambda/\gamma$ and $\varepsilon\approx \mu/\gamma$.

If we solve the system of differential equations (\ref{det-SIR-demog}) by setting all derivatives equal to 0, and replace $\mu$, $\lambda$ and $\gamma$ by the dimensionless quantities $R_0$ and $\varepsilon$ (three parameters can be replaced by two because the unit of time for the rates is arbitrary and one rate can be set to unity), we obtain the endemic level which is given by
\begin{equation}
(\hat s, \hat i, \hat r)= \left( \frac{1}{R_0},\ \varepsilon \left(1-\frac{1}{R_0}\right) ,\ 1- \frac{1}{R_0} - \varepsilon \left(1-\frac{1}{R_0}\right) \right) \label{end-level}
\end{equation}

\begin{exercise} \label{xopen-end-level}
Show that this is the endemic level, i.e.\ that the solution solves Equation (\ref{det-SIR-demog}) with all derivatives being 0.
\end{exercise}

This state is only meaningful if $R_0>1$ (otherwise some fraction is negative), so the endemic level only exists if $R_0>1$. Another solution to the equation system is of course the disease free equilibrium $(s,i,r)=(1,0,0)$. It is well known that when $R_0>1$ (which we from now on assume), then the endemic state is globally stable whereas the disease free state is locally unstable, meaning the system converges to the endemic level irrespective of starting value as long as $i(0)>0$.

Using the theory of Section \ref{TB-EP_sec_LLN} it can be shown that the current Markov model (for an open population) converges to the above deterministic model as $N\to\infty$, if the starting point is such that the fraction initially infectious is strictly positive ($I^N(0)/N\to i(0)>0$).

This suggests that the stochastic model (for the fractions in different states) can be approximated by the corresponding deterministic function \[(S^N(t)/N, I^N(t)/N, R^N(t)/N)\approx (s(t), i(t), r(t))\] which solves Equation (\ref{det-SIR-demog}) and having the same initial condition as the stochastic system. And, since we know that $(s(t), i(t), r(t))\to (\hat s, \hat i, \hat r)$ as $t\to\infty$ this suggests that $(S^N(t)/N, I^N(t)/N, R^N(t)/N)\approx  (\hat s, \hat i, \hat r)$ when $N$ and $t$ are large. This is indeed true in some sense, but it is only true depending on the relation between $N$ and $t$. For any finite $N$, the stochastic epidemic, which fluctuates randomly around the endemic equilibrium, will eventually go extinct, meaning that for some random $T^N_{Ext}$ (the extinction time) it will happen that $I^N(T^N_{Ext})=0$. When this happens the rate of new infections is 0 so the stochastic epidemic will remain disease free ever after (and eventually all removed will have died so all individuals are susceptible. Using large deviation theory (cf.\ Section \ref{TB-EP_sec_LD} below) it can be shown that the time to extinction grows exponentially with $N$, $T^N_{Ext}\approx e^{cN}$ for some $c>0$ as $N\to\infty$.

On the other hand, for any arbitrary but fixed time horizon $[0,t_{\max}]$ the stochastic epidemic will converge to the deterministic process as $N\to\infty$. It also follows from Theorem \ref{th:CLT} that the scaled process
\[\sqrt{N} ( S^N(t)/N-s(t),\ I^N(t)/N-i(t),\ R^N(t)/N-r(t))\]
converges to an Ornstein--Uhlenbeck process $(\tilde S(t),\tilde I(t), \tilde R(t))$. This Ornstein--\linebreak Uhlenbeck process is a Gaussian process with stationary distribution being Normally distributed. In particular, the variance of $\tilde I(t)$ in stationarity is well approximated by $1/R_0-1/R_0^2$, see N\aa sell \cite{N99I}.
\begin{exercise}\label{Nasell}
Show this as a consequence of Theorem \ref{th:CLT}, Lemma \ref{le:cov_inv}, and Exercise \ref{exerOU}
\end{exercise}
This suggests that $I^N(t)$ will be approximately Gaussian with mean $N\hat i$ and standard deviation $\sqrt{N/R_0}$ when $N$ is large and $t$ is moderately large (smaller than $T^N_{Ext}$ but still large since we assume the Ornstein--Uhlenbeck is close to stationary).

From above we know that $T^N_{Ext}$ will grow exponentially with $N$ as $N\to\infty$. On the other hand, if $N$ is small or moderate, the disease will go extinct very quickly, e.g.\ within a year. We now use the Gaussian approximation above to define
a sort of threshold, the \emph{critical population size} $N_c$,  between these two scenarios (quick extinction and very long time before extinction). Of course, there is no unique exact such value, so it will involve some arbitrary choice(s).

Above we noted that $I^N(t)$ was approximately Gaussian with mean $N\hat i$ and standard deviation $\sqrt{N/R_0}$. If we want to be above the critical population size, then we want to avoid quick extinction for which it is necessary that this approximately Gaussian process avoids extinction for a fairly long time. Extinction occurs when $I^N(t)=0$, and if we want to avoid this we want the value 0 to be far enough away from the mean, e.g.\ at least 3 standard deviations away. The choice 3 is of course arbitrary but if we instead choose 2 the \emph{process} will hit 0 fairly quickly with large enough probability, and if we choose 4 it seems extremely unlikely that it will hit extinction within e.g.\ a life time, so 3 seems like a reasonable compromise when it is unlikely but not completely impossible. This choice then suggests that the threshold is for the case $N\hat i -3\sqrt{N/R_0}=0$. This is equivalent to $\sqrt{N}=3/\hat i\sqrt{R_0}$, i.e.\ $N=9/\hat i^ 2R_0$. Inserting that $\hat i=\varepsilon (1-1/R_0)$ (remember that $\varepsilon=\mu/(\gamma+\mu)$ is the relative length of the infectious period compared to life-length), then we arrive at our definition of the \emph{critical population size} $N_c$:
\begin{equation}
N_c= \frac{9}{ \varepsilon^2( 1-\frac{1}{R_0} )^2 R_0} .\label{def-critpopsize}
\end{equation}
The conclusion is that, for a given infectious disease, i.e.\ given $R_0$ and $\varepsilon$, the disease will die out quickly in a community of size $N\ll N_c$ whereas it will persist for a very long time if $N\gg N_c$, during which the disease is endemic. As an illustration, consider measles prior to vaccination. If we assume that $R_0\approx 15$ and the infectious period is 1 week (1/52 years) and life duration 75 years, implying that $\varepsilon \approx \frac{1/75}{1/(1/52) + 1/75}\approx 1/3750$ we arrive at $N_c\approx 9(3750)^2/15 \approx 8\cdot 10^ 6$. So, if the population is a couple of million (or less) the disease will go extinct quickly, whereas the disease will become endemic (for a very long time) in a population being larger than e.g.\ 20 million people. This confirms the empirical observation that measles was continuously endemic in UK whereas it died out quickly in Iceland (and was later reintroduced by infectious people visiting the country).

\begin{exercise} \label{xopen-end-level1}
Which parameter affects $N_c$ the most? Compute $N_c$ using the measles example but making $R_0$ 50\% bigger/smaller and the same for the duration of the infectious period (assuming we live equally long).
\end{exercise}

\begin{exercise} \label{xopen-end-level2}
Suppose that a vaccine giving 100\% life long immunity is available, and that a fraction $v$ of all infants are continuously vaccinated. How does this affect the critical community size, i.e.\ give an expression for $N_c$ also containing $v$. (\emph{Hint}: Vaccinating people affects both the relevant population size $N_v$, the non-vaccinated population, and the reproduction number $R_v$, but other than that nothing has changed.)
\end{exercise}

\section[Large deviations and extinction of an endemic disease]{Large deviations and extinction of an endemic disease}\label{TB-EP_sec_LD}

\subsection{Introduction}\label{introLD}
In Section \ref{TB-EP_sec_LLN}, we have proved that, under appropriate conditions,
the solution of the SDE
\begin{equation}\label{SDE}
 Z^N_t=x_N+\sum_{j=1}^k\frac{h_j}{N}P_j\left(N\int_0^t\beta_{j}(s,Z^N_s)ds\right)
 \end{equation}
converges a.s., locally uniformly in $t$, towards the unique solution of the ODE
\begin{equation}\label{ODE}
\frac{dz_t}{dt}=b(t,z_t),\quad z_0=x,
\end{equation}
see Theorem \ref{th:LLN}, where $b(t,x)=\sum_{j=1}^k h_j\beta_j(t,x)$. Consequently the above SDE \eqref{SDE} can be considered for
large $N$ as a small random perturbation of the ODE \eqref{ODE}. Small random perturbations
of ODEs by Brownian motion have been studied by many authors, starting with Freidlin and Wentzell \cite{FWI}. Our aim is to study the above type of random perturbations of an ODE like
\eqref{ODE}. The starting point is the estimation of a large deviation from the law of large numbers, which has been studied for our type of Poisson driven SDEs by Shwartz and Weiss \cite{SWI}.
The difficulty is the fact that some of the rates in the SDE \eqref{SDE} vanish when the solution hits part of the boundary. This makes the estimate a bit delicate, since the logarithms of the rates enter the rate function in our large deviations estimate. This situation has been addressed first by Shwartz and Weiss \cite{SW05I}, but their assumptions are not quite satisfied in our framework. Recently Kratz and Pardoux \cite{KPI} and Pardoux and Samegni-Kepgnou \cite{PSI}
 have developed an approach to Large Deviations which is well adapted to the
epidemics models which are considered in these Notes. In fact the main difficulty concerns the lower bound. In the following, we present a new approach to the lower bound, based upon a quasi--continuity result, Proposition \ref{quasicont} below,  which mimics a similar result for Brownian motion driven SDEs due to Azencott \cite{RA}.  The same approach, for other types of Poisson driven SDEs,
will soon appear in Kouegou-Kamen and Pardoux \cite{KKP1}, \cite{KKP2}.

The main application we have in mind is to estimate the time needed for the small random perturbations to drive the system from a stable endemic equilibrium to the disease free equilibrium (i.e.\ extinction). This applies to the classical SIS and SIRS models, as well as to an SIR model with demography, as well as to models with vaccination and to models with several levels of susceptibility, thus predicting the time it will take for the random perturbation to
end an endemic disease.

We rewrite our model as
\[ Z^{N,x_N}_t=x_N+\sum_{j=1}^k h_j\int_0^t\int_0^{\beta_j(s,Z^{N,x_N}_s-)}Q^N_j(ds,du),\]
where
\[ Q^N_j(ds,du)=\frac{1}{N}Q_j(ds,Ndu),\]and the $Q_j$'s are i.i.d.\ Poisson random measures
on $[0,T]\times\R_+$, with mean $\lambda^2$, the $2$-dimensional Lebesgue measure.

\begin{description}
\item{(A.1)} We shall assume in all of this section that the $\beta_j$'s are locally Lipschitz with respect to $x$, uniformly for $t\in[0,T]$.
\end{description}

\subsection{The rate function}\label{sec:rate}
We want to establish a large deviations principle for trajectories in the space $D([0,T];\R^d)$ of $\R^d$-valued right-continuous functions which have a left limit at any time $t\in(0,T]$. We shall also consider the sets $C([0,T];\R^d)$
of continuous functions from $[0,T]$ into $\R^d$, and the subset of absolutely continuous functions, which we will denote $\mathcal{AC}_{T,d}$.
For any $\phi\in\mathcal{AC}_{T,d}$, let $\mathcal{A}_{k}(\phi)$ denote the (possibly empty) set of functions $c\in L^1(0,T;\R^k_+)$ such that $c_j(t)=0$ a.e. on the set $\{t,\, \beta_j(\phi_t)=0\}$ and
\begin{equation*}\label{allowed}
  \frac{d\phi_{t}}{dt}=\sum_{j=1}^{k} c_j(t) h_{j}, \quad\text{t a.e}.
\end{equation*}
We define the rate function
\begin{equation*}
I_{T}(\phi) :=
 \begin{cases}
\inf_{c\in\mathcal{A}_{k}(\phi)}I_{T}(\phi|c),& \text{ if } \phi\in\mathcal{AC}_{T,A}; \\
\infty ,& \text{ otherwise.}
\end{cases}
\end{equation*}
where as usual the infimum over an empty set is $+\infty$, and
\begin{equation*}
  I_{T}(\phi|c)=\int_{0}^{T}\sum_{j=1}^{k}g(c_j(t),\beta_{j}(\phi_{t}))dt
\end{equation*}
with $g(\nu,\omega)=\nu\log(\nu/\omega)-\nu+\omega$.
 We assume in the definition of $g(\nu,\omega)$ that for all $\nu>0$, $\log(\nu/0)=\infty$ and $0\log(0/0)=0\log(0)=0$.

We consider $I_T$ as a functional defined on the space $D([0,T];\R^d)$ equipped with Skorokhod's topology.
We first give two other possible definitions of the functional $I_T$. Let $\ell:\R^{3d}\mapsto\R$ be defined as
\[\ell(x,y,\theta)=\langle y,\theta\rangle-\sum_{j=1}^k\beta_j(x)\left(e^{\langle h_j,\theta\rangle}-1\right).\]
We define the map $L:\R^{2d}\mapsto(-\infty,+\infty]$  as
\[ L(x,y)=\sup_{\theta\in\R^d} \ell(x,y,\theta)\, .\] We let
\[\hat{I}_T(\phi)=\int_0^T L(\phi_t,\dot{\phi}_t)dt.\]
It is not hard to see that the following is an equivalent definition of $\hat{I}_T(\phi)$:
\[  \hat{I}_T(\phi)=\sup_{\theta\in C^1([0,T];\R^d)}\int_0^T\ell(\phi_t,\dot{\phi}_t,\theta_t)dt\, .\]

We first establish
\begin{proposition}\label{pro:equivRatef}
For any $\phi\in D([0,T]:\R^d)$, $I_{T}(\phi)=\hat{I}_T(\phi)$.
\end{proposition}
\bpf
We note that if $y=\sum_{j=1}^kc_jh_j$ with some $c\in\R^k_+$,
\[ \ell(x,y,\theta)=\sum_{j=1}^k\left[ c_j\langle h_j,\theta\rangle-\beta_j(x)\left( e^{\langle h_j,\theta\rangle}-1\right)      \right]\, .\] But for any $1\le j\le k$,
\begin{align*}
c_j\langle h_j,\theta\rangle-\beta_j(x)\left( e^{\langle h_j,\theta\rangle}-1\right)
&\le\displaystyle \sup_{r\in\R}\left[c_j r-\beta_j(x)\left(e^r-1\right)\right] \\
&= c_j\log\left(\frac{c_j}{\beta_j(x)}\right)-c_j+\beta_j(x)\\
&= g(c_j,\beta_j(x)).
\end{align*}
The inequality $\hat{I}_T(\phi)\le I_T(\phi)$ for any $\phi\in D([0,T];\R^d)$ follows readily.

%
In order to prove the converse inequality, we fix $x,y\in\R^d$ such that $L(x,y)<\infty$ (otherwise there is nothing to prove). Let
$\theta_n$ be a sequence in $\R^d$ such that $L(x,y)=\lim_{n\to\infty}\ell(x,y,\theta_n)$. It is clear that
for any $1\le j\le k$ such that $\beta_j(x)>0$, the\linebreak sequence $\langle\theta_n,h_j\rangle$ is bounded from above. Hence we can and do assume that, after the extraction of a subsequence, for any $1\le j\le k$ such that $\beta_j(x)>0$, the sequence $e^{\langle\theta_n,h_j\rangle}\to  s_j$, for some $s_j\ge0$.  Consequently,
as $n\to\infty$,
\begin{equation}\label{conv}
 \langle\theta_n,y\rangle\to L(x,y)+\sum_{j=1}^k\beta_j(x)\left(s_j-1\right).
\end{equation}
Differentiating $\ell(x,y,\theta_n)$ with respect to its last variable, we get
\begin{align*}
\nabla_\theta\ell(x,y,\theta_n)&=y-\sum_{j=1}^k\beta_j(x)e^{\langle h_j,\theta_n\rangle}h_j\\
&\to y-\sum_{j=1}^k\beta_j(x)s_jh_j,
\end{align*}
as $n\to\infty$.
But since $\theta_n$ is a maximizing sequence and the gradients converge, then since $L(x,y)<\infty$, their limit must be zero. Consequently
\[ y=\sum_{j=1}^k\beta_j(x)s_jh_j\, .\]
Hence, with $c_j=\beta_j(x)s_j$, we have
\begin{align*}
\langle\theta_n,y\rangle&=\sum_{j=1}^kc_j\langle\theta_n,h_j\rangle\\
&\to\sum_{j=1}^kc_j\log(s_j),
\end{align*}
with the convention that $c_j\log(s_j)=0$ if both $c_j=0$ and $s_j=0$.
This, combined with \eqref{conv}, yields that
\[ L(x,y)=\sum_{j=1}^kg(c_j,\beta_j(x))\]
which entails that $\hat{I}_T(\phi)\ge I_T(\phi)$. The proposition is established. \epf

We have the
\begin{proposition}\label{I=0}
For any $T>0$, $\phi\in D([0,T];\R^d)$, $I_T(\phi)\ge0$, and $I_T(\phi)=0$ iff $\phi$ solves the ODE \eqref{ODE}.
\end{proposition}
\bpf It suffices to show that $L(x,y)\ge L(x,\sum_j\beta_j(x)h_j)=0$, with strict inequality if
$y\not=\sum_j\beta_j(x)h_j$. We first note that
\begin{align*}
L\left(x,\sum_j\beta_j(x)h_j\right)&=\sup_\theta \left\{\sum_j \beta_j(x)(\langle h_j,\theta\rangle
-\exp\langle h_j,\theta\rangle+1)\right\}=0,
\end{align*}
since $z-e^z +1\le0$, with equality at $z=0$. Let now $y$ be such that $L(x,y)=0$. Then
\[ \langle y,\theta\rangle-\sum_j\beta_j(x)(\exp\langle h_j,\theta\rangle-1)\le0 \ \text{ for all }\theta\in\R^d.\]
Choosing $\theta=\eps e_i$ (where $e_i$ is the $i$-th basis vector of $\R^d$) yields
\[ \eps y_i\le\sum_j\beta_j(x)(\exp(\eps h_j^i)-1).\]
Dividing by $\eps$, then letting $\eps\to0$ yields $y_i\le\sum_j\beta_j(x)h_j^i$, while the opposite inequality follows if we start with $\theta=-\eps e_i$. The result follows.  \epf

In the next statement, we use the notion of a lower semi-continuous real-valued function, which is defined in Definition \ref{de:sc} below. In the proof we use the notion of an equicontinuous collection of functions, which is defined in Definition \ref{de:equicont}.
\begin{theorem}\label{th:goodratef}
$\phi\to I_{T}(\phi)$ is lower semi-continuous on $D([0,T];\R^d)$, and for any $R$, $K>0$, the set
$\{\phi\in D([0,T];\R^d),\, \sup_{0\le t\le T}|\phi_t|\le R,\, I_T(\phi)\le K\}$
is compact.
\end{theorem}
\bpf The lower semicontinuity property is an immediate consequence of the fact that, from its second definition, $\hat{I}_T$ is a supremum over continuous functions. To finish the proof, it suffices from the Arzel\`a--Ascoli theorem (see e.g.\ Theorem 7.2 in 
Billingsley \cite{PBI}) to show that the set of functions satisfying
$\sup_{0\le t\le T}|\phi_t|\le R$ and $I_T(\phi)\le K$ is equicontinuous. It is clear that if $\bar{h}=\sup_{1\le j\le k}|h_j|$
and $\bar{\beta}_R=\sup_{1\le j\le k}\sup_{0\le t\le T,\, |x|\le R}\beta_j(t,x)$,
\begin{align*}
 L(x,y)&\ge\ell\left(x,y,\frac{y\log(|y|)}{\bar{h}|y|}\right)\\
 &\ge\frac{|y|\log(|y|)}{\bar{h}}-k\bar{\beta}_R|y|.
 \end{align*}
 Now let $0\le s<t\le T$, with $t-s\le\delta$.
 \begin{align*}
 |\phi_t-\phi_s|&\le\int_s^t|\dot{\phi}_r|dr\\
 &\le\delta^{-1/2}\int_s^t{\bf1}_{|\dot{\phi}_r|\le\delta^{-1/2}}dr+
\int_s^t{\bf1}_{|\dot{\phi}_r|>\delta^{-1/2}}\frac{L(\phi_r,\dot{\phi}_r)}{L(\phi_r,\dot{\phi}_r)/|\dot{\phi}_r|}dr\\
&\le \delta^{1/2}+\frac{K}{f(\delta^{-1/2})},
\end{align*}
where $f(a)=\inf_{|x|\le R,\, |y|\ge a}\frac{L(x,y)}{|y|}$. The result follows from the fact that from the above lower bound of  $L(x,y)$, $f(a)\to\infty$, as $a\to\infty$. \epf

\subsection{The lower bound}
Let $\eta=(\eta_1,\ldots,\eta_k)$ be a vector of locally finite measures on $[0,T]\times\R_+$. We shall say that
$\eta\in\mathcal{M}^k$. To $x\in\R^d$ and $\eta\in\mathcal{M}^k$, we associate $\Phi^x_t(\eta)$, solution (if it exists)
of the ODE
\[ \Phi^x_t(\eta)=x+\sum_{j=1}^kh_j\int_0^t\int_0^{\beta_j(s,\Phi^x_{s-})}\eta_j(ds,du).\]
If $\eta_j(ds,du)=f_j(s,u)dsdu$, $1\le j\le k$, the above ODE has at least one solution (possibly up to an explosion time, as the solution of an ODE with continuous coefficients). If moreover
\[ \sup_{u\ge0} f_j(\cdot,u)\in L^1[0,T],\ 1\le j\le k,\]
then the above ODE has  a unique solution (as the solution of an ODE with locally Lipschitz coefficients).

Let $\phi\in C([0,T];\R^d)$ be an absolutely continuous function. We define
\begin{equation}\label{Kphi}
 K_\phi:=\inf_{c\in\mathcal{A}_k(\phi)}\sum_{j=1}^k\int_0^T\frac{c_j(t)}{\beta_j(t,\phi_t)}dt.
 \end{equation}
To a pair $(\phi,c)$ with $c\in\mathcal{A}_k(\phi)$, we associate for $1\le j\le k$ the measure $\eta_j(ds,du)$ with the density
\[ f_j(s,u)=\frac{c_j(s)}{\beta_j(s,\phi_s)}{\mathbf1}_{[0,\beta_j(s,\phi_s)]}(u)+{\bf1}_{(\beta_j(s,\phi_s),+\infty)}(u).\]
Then, with $x=\phi_0$, $\phi_t= \Phi^x_t(\eta)$.

Moreover, given $\phi\in C([0,T];\R^d)$ and $L>0$, we consider the set
\[ A_{\phi,L}=\{(t,x),\ 0\le t\le T,\ |x-\phi_t|\le L+1\},\] and define
\[ \overline{\beta}(\phi,L)=\sup_{1\le j\le k}\sup_{(t,x)\in A_{\phi,L}} \beta_j(t,x).\]

We can now prove the following.
\begin{proposition}\label{quasicont}
Let $T>0$ be arbitrary. Given $(\phi,\eta)$ as above, such that in particular $K_\phi<\infty$,
if $x_N=Z^N_0$, for any $R,\, L>0$, there exists a $\delta, r>0$ (depending upon $K_\phi$) and $N_0$ such that
whenever $|x-x_N|\le r$,  $N\ge N_0$,
\[ \P\left(\|Z^N-\phi\|_T>L,\, d_{T,\overline{\beta}}(Q^N,\eta)\le\delta\right)\le e^{-NR},\]
where
\[ d_{T,\overline{\beta}}(\nu,\eta)=\sum_{j=1}^k\sup_{0\le t\le T,\, 0\le u\le \overline{\beta}}
|\nu_j([0,t]\times[0,u])-\eta_j([0,t]\times[0,u])|,\]
and $ \overline{\beta}:=\overline{\beta}(\phi,L)$.
\end{proposition}
\bpf It is clear that
\begin{align*}
|Z^N_t-\phi_t|&\le|x_N-x|+\sum_{j=1}^k|h_j|\left|\int_0^t\int_0^{\beta_j(s,Z^N_{s-})}[Q^N_j(ds,du)-\eta_j(ds,du)]\right|\\
&\qquad\qquad +
\sum_{j=1}^k|h_j|\left|\int_0^t\int_{\beta_j(s,Z^N_{s-})\wedge\beta_j(s,\phi_s)}^{\beta_j(s,Z^N_{s-})\vee\beta_j(s,\phi_s)}f_j(s,u)\vee1 \, du\,ds\right|\\
&\le r+\sum_{j=1}^k|h_j|\left|\int_0^t\int_0^{\beta_j(s,Z^N_{s-})}[Q^N_j(ds,du)-\eta_j(ds,du)]\right|\\
&\qquad\qquad+\sum_{j=1}^k|h_j| C\int_0^t\left(\frac{c_j(s)}{\beta_j(s,\phi_s)}\vee1\right)|Z^N_s-\phi_s|ds,
\end{align*}
where $C$ is an upper bound of the Lipschitz constants of the $\beta_j$'s in $[0,T]\times[0,\overline{\beta}]$.
Subdividing  $[0,T]$ into $\left[\frac{T}{\rho}\right]+1$ intervals of the form $\displaystyle [(i-1)\rho, i\rho\wedge T]$ and denoting
\[\overline{\beta^i_j} := \sup_{(i-1)\rho\leq s\leq i\rho}\beta_j(s,Z^{N,x^N}_{s-}) \,  , \qquad
   \underline{\beta^i_j}:= \inf_{(i-1)\rho\leq s\leq i\rho}\beta_j(s,Z^{N,x^N}_{s-}) \,  , \]
   we define the random sets
   \[ A^{\rho,i}_j := [(i-1)\rho,i\rho]\times[0\:,\:\underline{\beta^i_j}], \quad      B^{\rho,i}_j :=  [(i-1)\rho,i\rho]\times[\underline{\beta^i_j}\:,\:\overline{\beta^i_j}  ] \, .\]
   For all $i$ and $j$,
   \[\sum_{j=1}^k|Q^N_j(A^{\rho,i}_j)-\eta_j(A^{\rho,i}_j)|\le 2d_{T,\bar{\beta}}(Q^N,\eta),\;
   \sum_{j=1}^k|Q^N_j(B^{\rho,i}_j)-\eta_j(B^{\rho,i}_j)|\le 4d_{T,\bar{\beta}}(Q^N,\eta).\]
Consequently for all $0\le t\le T$, if $\bar{h}:=\sup_{1\le j\le k} |h_j|$, then on the event\linebreak $\{d_{T,\overline{\beta}}(Q^N,\eta)\le\delta\}$,

\begin{align*}
\sum_{j=1}^k&|h_j|\Bigg| \int_{0}^{t}\int_{0}^{\beta_j\left(s,Z^{N,s-}\right)} \left[Q^N_j\left(ds,du\right) -\eta_j(ds,du) \right] \Bigg| \\
&\le\bar{h}\sum_{j=1}^k\left(\sum_{i=1}^{\left[\frac{t}{\rho}\right]+1}\left| Q^N_j(A^{\rho,i}_j)-\eta_j(A^{\rho,i}_j) \right| + \sum_{i=1}^{\left[\frac{t}{\rho}\right]+1}\left\{Q^N_j\left(B^{\rho,i}_j\right)+\eta_j\left(B^{\rho,i}_j\right)\right\}\right)
\\
&\le\bar{h}\sum_{j=1}^k\left(\sum_{i=1}^{\left[\frac{t}{\rho}\right]+1}\left| Q^N_j(A^{\rho,i}_j)-\eta_j(A^{\rho,i}_j) \right| + \sum_{i=1}^{\left[\frac{t}{\rho}\right]+1}\left|Q^N_j\left(B^{\rho,i}_j\right)-\eta_j\left(B^{\rho,i}_j\right)\right|\right.\\
&\qquad\qquad\left.+2\sum_{i=1}^{\left[\frac{t}{\rho}\right]+1}\eta_j\left(B^{\rho,i}_j\right)\right)\\
&\le 6\left(\frac{t}{\rho}+1\right)\bar{h}\delta +2\bar{h}\sum_{j=1}^k\sum_{i=1}^{\left[\frac{t}{\rho}\right]+1}\eta_j\left(B^{\rho,i}_j\right).
\end{align*}

It  follows from the two above inequalities and Gronwall's Lemma that
\begin{equation}\label{useGron}
\sup_{0\le t\le T}|Z^N_t-\phi_t|\le\Big(r+6\!\left(\!\frac{T}{\rho}\!+\!1\!\right)\bar{h}\delta +2\bar{h}\sum_{j=1}^k\sum_{i=1}^{\left[\frac{t}{\rho}\right]+1}\eta_j\left(B^{\rho,i}_j\right)\!\Big)\exp\left[C(K_\phi+kT)\bar{h}\right].
\end{equation}
Since the $\left(B^{\rho,i}_j\right)_i$ are disjoints we have for all $j$
\begin{align*}
 \sum_{i=1}^{\left[\frac{T}{\rho}\right]+1}\eta_j(B^{\rho,i}_j) = \eta_j\left(  \bigcup_{i=1}^{\left[\frac{T}{\rho}\right]+1}B^{\rho,i}_j  \right) &\leq \sum_{i=1}^{\left[\frac{T}{\rho}\right]+1}(\bar{\beta^i_j} - \underline{\beta^i_j})\int_{(i-1)\rho}^{i\rho} \frac{c_j(s)}{\beta_j(\phi_s)}\vee 1 \: \:ds\\
&\leq \max_{1\leq i \leq \left[\frac{T}{\rho}\right]+1}\left(\overline{\beta^i_j} - \underline{\beta^{i}_j}\right) \int_0^T \frac{c_j(s)}{\beta_j(\phi_s)}\vee 1 \: ds\\
&\leq (K_{\phi} + T)\max_{1\leq i \leq \left[\frac{T}{\rho}\right]+1}\left(\overline{\beta^i_j} - \underline{\beta^{i}_j}\right).
\end{align*}
We note that  for every $i, j$
\[ \overline{\beta^i_j} -  \underline{\beta^i_j}  \leq C\frac{X_i}{N} \]
where $X_i$ is a Poisson random variable of mean $\rho N\bar{\beta}$.  For  any $a>0$, we have
with $\bar{a}=\frac{a}{k(K_\phi+T)}$, using Cram\'er's Theorem \ref{th:cramer} for the fourth inequality,
\begin{align}
 \P&\left[ \sum_{j=1}^{k}\sum_{i=1}^{\left[\frac{T}{\rho}\right]+1}\eta_j(B^{\rho,i}_j) >  a\right]\leq k\max_{j} \: \P\left[ \sum_{i=1}^{\left[\frac{T}{\rho}\right]+1}\eta_j(B^{\rho,i}_j) >  \frac{a}{k} \right]\nonumber\\
 &\qquad\qquad\leq k \: \mathbb{P}\left[ \max_{1\leq i\leq \left[\frac{T}{\rho}\right]+1} \frac{X_i}{N} >\bar{a}\right]\nonumber\\
 &\qquad\qquad\leq k \: \mathbb{P}\left[ \bigcup_{1\leq i\leq \left[\frac{T}{\rho}\right]+1} \left\{\frac{X_i}{N} >\bar{a}\right\}\right] \label{estimprob}\\
 &\qquad\qquad\leq  k\left(\frac{T}{\rho}+1\right)\exp\left( -N\left[ \bar{a}\log\frac{\bar{a}}{\rho\bar{\beta}} + \bar{a} - \rho\bar{\beta}\right] \right)\nonumber\\
 &\qquad\qquad=\exp\left( -N\left[ \bar{a}\log\frac{\bar{a}}{\rho\bar{\beta}} + \bar{a} -\frac{1}{N}\log\left(k\left[\frac{T}{\rho}+1\right]\right) - \rho\bar{\beta}\right] \right).\nonumber
\end{align}
 We choose $\rho=\sqrt{\delta}$. Let $\delta_0$ be such that
\begin{align*}
6\left(T\sqrt{\delta_0}+\delta_0\right)\bar{h}&\le\frac{L}{3}\exp\left[-C(K_\phi+kT)\bar{h}\right],\ \text{and}\\
r&=\frac{L}{3}\exp\left[-C(K_\phi+kT)\bar{h}\right],\\
 a&=\frac{L}{6\bar{h}}\exp\left[-C(K_\phi+kT)\bar{h}\right],
 \end{align*}
 so that from \eqref{useGron},
 \begin{equation}\label{inclusion}
 \left\{\sum_{j=1}^{k}\sum_{i=1}^{\left[\frac{T}{\rho}\right]+1}\eta_j(B^{\rho,i}_j) \le  a\right\}\subset
 \Big\{\|Z^N-\phi\|_T\le L\Big\}.
 \end{equation} $R>0$ being arbitrary, we now choose
 \begin{align*}
  \delta&=\min\left\{ \delta_0,\left(\frac{\bar{a}}{\bar{\beta}}\right)^2e^{-2R/\bar{a}},\frac{\bar{a}}{2\bar\beta}\right\},
  \ \text{ and}\\
 N_0&=\left\lceil\frac{2}{\bar{a}}\log\left(k\left[\frac{T}{\rho}+1\right]\right)\right\rceil.
 \end{align*}
The result  follows from those choices,  \eqref{estimprob} and \eqref{inclusion}.
  \epf

Before we establish the lower bound, we need to formulate an assumption.
\begin{description}
\item{(A.2)} We assume that for any $\phi\in C([0,T];\R^d)$ such that $I_T(\phi)<\infty$ and any $\eps>0$, there exists a $\phi^\eps$ such that $\phi^\eps_0=\phi_0$,
$K_{\phi^\eps}<\infty$, $\|\phi-\phi^\eps\|_T\le \eps$ and $I_T(\phi^\eps)\le I_T(\phi)+\eps$.
\end{description}
\begin{exercise}
Consider the SIRS model with fixed population size, and let $A:=\{(x,y),\, 0\le x, 0\le y, x+y\le1\}$. Show that if $\phi\in C([0,T];A)$ hits the boundary, then for any $\eps>0$,
one can find $\phi^\eps$ such that $\phi^\eps_0=\phi_0$,
$K_{\phi^\eps}<\infty$, $\|\phi-\phi^\eps\|\le \eps$ and $I_T(\phi^\eps)\le I_T(\phi)+\eps$, where $\phi^\eps$ can either remain in the interior of $A$, or else can hit the boundary.
\end{exercise}
We now have, with the notation $I_{T,x}(O)=\inf_{\phi\in O,\, \phi_0=x}I_T(\phi)$,
\begin{theorem}\label{th:LB}
If the assumptions (A.1) and (A.2) are satisfied, then for any open subset $O\subset D([0,T];\R^d)$,
if $x_N\to x$ as $N\to\infty$,
\[\liminf_{N\to\infty}\frac{1}{N}\log\P\left(Z^{N,x_N}\in O\right)\ge -I_{T,x}(O).\]
\end{theorem}
\bpf It clearly suffices to treat the case where $I_{T,x}(O)<\infty$. Then for any $\eps>0$ there exists a $\phi\in O$ such that $\phi_0=x$ and
\[ I_T(\phi)\le I_{T,x}(O)+\frac{\eps}{4}.\]
It follows from assumption (A.2) that there exists a $\hat{\phi}\in O$ such that $\hat{\phi}_0=\phi_0$, $K_{\hat{\phi}}<\infty$,
$\|\hat{\phi}-\phi\|_T\le\eps$ and
\[I_T(\hat{\phi})\le I_T(\phi)+\frac{\eps}{4}.\]
Now there exists a $c\in\mathcal{A}_k(\phi)$ such that $\sum_{j=1}^k\int_0^T\frac{c_j(t)}{\beta_j(t,\phi_t)}dt<\infty$, and
\[I_T(\hat{\phi}|c)\le I_T(\hat{\phi})+\frac{\eps}{4}.\]
If $\eps$ has been chosen small enough, there exists an $L>0$ be such that\linebreak $\{\psi;\, \|\psi-\hat{\phi}\|_T<L\}\subset O$.
From Proposition \ref{quasicont}, if $\eta^c$ denotes the vector of measures associated to $c$, $|x-x_N|$ is small enough and $N$ large enough,
for any $R>0$, there exists a $\delta>0$ such that with $\hat{\beta}=\overline{\beta}(\hat{\phi},L)$,
\begin{align}
\P\left(Z^{N,x_N}\in O\right)&\ge
\P\left(\|Z^{N,x_N}-\phi\|_T<L\right)\nonumber\\
&\ge\P\left(d_{T,\hat{\beta}}(Q^N,\eta^c)<\delta\right)\nonumber\\
&\qquad\qquad-\P\left(\|Z^{N,x_N}-\phi\|_T>L,d_{T,\hat{\beta}}(Q^N,\eta^c)<\delta\right)\nonumber\\
&\ge\P\left(d_{T,\hat{\beta}}(Q^N,\eta^c)<\delta\right)-e^{-NR}.\label{borninf}
\end{align}
Let us admit for a moment the next lemma.
\begin{lemma}\label{le:part}
 There exists a sequence of partitions $\{A^i_n,\, 1\le i\le a_n\}$ of $[0,T]\times[0,\hat{\beta}]$
such that $\sup_i\lambda^2(A^i_n)\to 0$ as $n\to\infty$, and a sequence $\delta_n\downarrow0$ and $n_0$
such that for all $n\ge n_0$,
\[\bigcap_{j=1}^k\bigcap_{i=1}^{a_n}\left\{Q_j^N(A^i_n)\in(\eta_j^c(A^i_n)-\delta_n,\eta_j^c(A^i_n)+\delta_n)\right\}\subset
\{ d_{T,\hat{\beta}}(Q^N,\eta^c)<\delta\}.  \]
\end{lemma}
As a consequence of this lemma, making use of Cram\'er's Theorem \ref{th:cramer} for the second inequality,
\begin{align*}
\liminf_{N\to\infty}\frac{1}{N}&\log\P\left(d_{T,\hat{\beta}}(Q^N,\eta^c)<\delta\right)\\&\ge
\sum_{j=1}^k\sum_{i=1}^{a_n}\liminf_{N\to\infty}\frac{1}{N}\log\P\left(Q_j^N(A^i_n)\in(\eta_j^c(A^i_n)-\delta_n,\eta_j^c(A^i_n)+\delta_n)\right)
\\
&\ge -\sum_{j=1}^k\sum_{i=1}^{a_n}\left(\eta_j^c(A^i_n)\log\frac{\eta_j^c(A^i_n)}{\lambda^2(A^i_n)}-\eta_j^c(A^i_n)+\lambda^2(A^i_n)\right)\\
&\ge-\sum_{j=1}^{k}\int_0^T\int_0^{\bar{\beta}}\left[ f_j^c(s,u)\log[f_j^c(s,u)]-f_j^c(s,u)+1\right]dsdu-\frac{\eps}{4}\\
&=-\sum_{j=1}^{k}\int_0^T\left[c_j(s)\log\frac{c_j(s)}{\beta_j(s,\phi_s)}-c_j(s)+\beta_j(s,\phi_s)\right] ds-\frac{\eps}{4}\\
&=-I_T(\hat{\phi}|c)-\frac{\eps}{4}\\
&\ge-I_{T,x}(O)-\eps,
\end{align*}
where \[f^c_j(s,u)=\frac{c_j(s)}{\beta_j(s,\phi_s)}{\mathbf1}_{[0,\beta_j(s,\phi_s)]}(u)+{\mathbf1}_{(\beta_j(s,\phi_s),+\infty)}(u)\]
and the second inequality holds true for $n$ chosen large enough as a function of $\eps$.
We let $\eps\to0$, and to combine the resulting inequality with \eqref{borninf}, hence
\begin{align*}
-I_{T,x}(O)&\le\liminf_{N\to\infty}\frac{1}{N}\log\left(\P\left(Z^{N,x_N}\in O\right)+e^{-NR}\right)\\
&\le\left(\liminf_{N\to\infty}\frac{1}{N}\log\P\left(Z^{N,x_N}\in O\right)\right)\vee(-R).
\end{align*}
The result finally follows by letting $R\to\infty$.
\epf

We now need to pass to the
\smallskip

\begin{proof}[Proof of Lemma \ref{le:part}]
For convenience, we replace the partition $\{A^i_n,\ 1\le i\le a_n\}$ by a partition
$\{A^{i,j}_n,\ 1\le i,j\le n\}$, which we construct as follows. We first choose $0=\beta_n^0<\beta_n^1<\cdots<\beta_n^n=\hat{\beta}$ such that
\[ \sup_{1\le j\le n}\eta^c([0,T]\times(\beta_n^{j-1},\beta_n^j])\le\frac{2}{n}\eta^c([0,T]\times[0,\hat{\beta}]).\]
We next choose a sequence $0=t^0_n<t^1_n<\cdots t_n^n=T$ such that, if
$A_n^{i,j}=(t_n^{i-1},t_n^i]\times(\beta_n^{j-1},\beta_n^j]$,
\[\sup_{1\le i\le n}\eta^c(A_n^{i,j})\le \frac{2}{n}\eta^c([0,T]\times(\beta_n^{j-1},\beta_n^j])
\le \frac{4}{n^2}\eta^c([0,T]\times[0,\hat{\beta}]):=\frac{C}{n^2}.\]
For an arbitrary $0\le t\le T$ and $0\le\alpha\le\hat{\beta}$, we define the set
\[ \partial_{t,\alpha}=\{t\}\times[0,\hat{\beta}]\cup[0,T]\times\{\alpha\},\]
which is the ``boundary'' of $[0,t]\times[0,\alpha]$. We note that $|\{i,j,\ A_n^{i,j}\cap\partial_{t,\alpha}\not=\emptyset\}|\le 2n$.
We need to bound
\begin{align*}
\Big|Q^N([0,t]&\times[0,\alpha])-\eta^c([0,t]\times[0,\alpha])\Big|\\&
\le\!\!\!\!\sum_{i,j,\ A_n^{i,j}\subset[0,t]\times[0,\alpha]}\!\!\!\!\!\!\!\left|Q^N(A_n^{i,j})-\eta^c(A_n^{i,j})\right|
+\!\!\!\!\!\!\!\sum_{i,j,\ A_n^{i,j}\cap\partial_{t,\alpha}\not=\emptyset}\!\!\!\!\!\!\!\left(Q^N(A_n^{i,j})+\eta^c(A_n^{i,j})\right)\\
&\le n^2\delta_n+2n\left(\frac{2C}{n^2}+\delta_n\right)\\
&\le \delta,
\end{align*}
for all $n\ge n_0$,
provided we choose first $n_0\ge\frac{8C}{\delta}$, and then a sequence $\delta_n$ such that
$\delta_n\le[2(n^2+2n)]^{-1}\delta$ for each $n\ge n_0$. \epf

We now establish a slightly stronger result. Here and below we shall use the following notation concerning the initial condition of $Z^N$. We fix $x\in\R^d$ and start $Z^N$ from the point $Z^N_0=x_N$, where the $i$-th coordinate $x_N^i$ of $x_N$ is given by $x^i_N=\frac{[x^i N]}{N}$. Here we assume that the process $Z^{N}$ lives in a closed subset $A\subset\R^d$. We shall need the following
\begin{definition}\label{def:adapted}
We shall say that the compact set of initial conditions $\mathcal{K}$ is {\bf adapted} to the open set of trajectories
$O\subset D([0,T];A)$ if
\begin{enumerate}
\item $\mathcal{K}\subset\{\phi_0,\, \phi\in O\}$.
\item For any $\eps>0$, the following holds. For any $x\in\mathcal{K}$, there exists a $\phi^x\in O$ such that $\phi^x_0=x$,
$I_T(\phi^x)\le I_{T,x}(O)+\eps$ and moreover $\sup_{x\in\mathcal{K}}K_{\phi^x}<\infty$.
\end{enumerate}
\end{definition}
It follows readily from the proof of Theorem \ref{th:LB} that the following reinforced version holds.
\begin{theorem}\label{th:LBunif}
For any open subset $O\subset D([0,T];A)$ and any compact subset
$\mathcal{K}$ of initial conditions which is adapted to $O$,
\[\liminf_{N\to\infty}\frac{1}{N}\log\inf_{x\in\mathcal{K}}\P(Z^{N,x_N}\in O)\ge -\sup_{x\in \mathcal{K}}I_{T,x}(O)\,.\]
\end{theorem}

\subsection{The upper bound}
In this subsection, we shall again use the notation $x_N$ for the vector whose $i$-th coordinate is given by $x^i_N=\frac{[x^i N]}{N}$.
We want to prove that for any closed $F$, $F\subset D([0,T];\R^d)$,
\begin{equation}\label{eq:ub}
\limsup_{N\to\infty}\log \P(Z^{N,x_N}\in F)\le -I_{T,x}(F).
\end{equation}
Let us recall the concept of exponential tightness.
\begin{definition}\label{de:expt}
The sequence $Z^N$ is said to be exponentially tight if for any $\alpha>0$, there exists a compact $K^\alpha$ such that
\[\limsup_N\frac{1}{N}\log\P(Z^N\in K_\alpha^c)\le -\alpha.\]
\end{definition}
We have the following lemma.
\begin{lemma}\label{le:expt}
If \eqref{eq:ub} holds for any compact subset $F=K\subset\subset D([0,T];A)$, and $Z^N$ is exponentially tight, then \eqref{eq:ub}
holds for any closed subset $F\subset D([0,T];A)$.
\end{lemma}
\bpf Let $F$ be closed and $\alpha:=I_{T,x}(F)$. We assume w.l.o.g. that $\alpha>0$ (unless the conclusion below would be obvious).
Let $K_\alpha$ be the compact set associated to $\alpha$ by Definition \ref{de:expt}. It is clear that $F\cap K_\alpha$ is compact and
$I_{T,x}(F\cap K_\alpha)\ge\alpha$. Hence from our assumption
\[ \limsup_{N\to\infty} \frac{1}{N}\log\P(Z^N\in F\cap K_\alpha)\le-\alpha.\]
Also from the choice of $K_\alpha$,
\[ \limsup_{N\to\infty} \frac{1}{N}\log\P(Z^N\in  K_\alpha^c)\le-\alpha.\]
But $\P(Z^N \in F)\le\P(Z^N\in F\cap K_\alpha)+\P(Z^N\in  K_\alpha^c)$, hence
\[\log\P(Z^N \in F)\le\log2+\sup(\log\P(Z^N\in F\cap K_\alpha),\log\P(Z^N\in  K_\alpha^c)),\]
and we clearly deduce that
\[\limsup_{N\to\infty} \frac{1}{N}\log\P(Z^N\in  F)\le-\alpha,\]
as desired. \epf

 Let us first establish
\begin{theorem}\label{th:UB-compact}
Let $T>0$ and $x\in\R^d$ be fixed. Let $x_N\to x$ as $N\to\infty$. For any compact set $K\subset D([0,T];\R^d)$,
\[ \limsup_{N\to\infty}\frac{1}{N}\log\P\left(Z^{N,x_N}\in K\right)\le -I_{T,x}(K)\, .\]
\end{theorem}
\bpf
Recall the formula
\begin{align*}
I_T(\phi)&=\sup_{\theta\in C^1([0,T];\R^d)}\int_0^T\ell(\phi_t,\dot{\phi}_t,\theta_t)dt\\
&=\sup_{\theta\in C^1([0,T];\R^d)}\L(\phi,\theta),
\end{align*}
where
\[ \L(\phi,\theta)=\langle\phi_T,\theta_T\rangle-\langle\phi_0,\theta_0\rangle-\int_0^T\langle\phi_t,\dot{\theta}_t\rangle dt
-\sum_{j=1}^k\int_0^T\beta_j(\phi_t)\left[e^{\langle h_j,\theta_t\rangle}-1\right]dt\, .\]
For any $\theta\in C^1([0,T];\R^d)$, $0\le s<t\le T$, we define
\begin{align*}
M^{N,\theta}_{s,t}&=\!\langle Z^{N,x_N}_t,\theta_t\rangle\!-\!\langle Z^{N,x_N}_s,\theta_s\rangle\!
-\!\int_s^t\!\!\langle Z^{N,x_N}_r,\dot{\theta}_r\rangle dr\!-\!\sum_{j=1}^k\!\!\int_s^t\!\!\langle h_j,\theta_r\rangle\beta_j(r,\Z^{N,x_N}_r)dr,\\
\Xi^{N,\theta}_{s,t}&=\exp\left(NM^{N,\theta}_{s,t}-N\sum_{j=1}^k\int_s^t\tau(\langle h_j,\theta_r\rangle)
\beta_j(r,\Z^{N,x_N}_r)dr\right),
\end{align*}
where $\tau(a)=e^a-1-a$, are such that $M^{N,\theta}_{0,t}$ and $\Xi^{N,\theta}_{0,t}$ are local martingales, the second being also a supermartingale such that $\E[\Xi^{N,\theta}_{0,t}]\le1$.

We assume that $I_{T,x}(K)>0$, since otherwise the result is trivial. We also assume that $I_{T,x}(K)<\infty$. The case $I_{T,x}(K)=\infty$ can be treated in a way which is very similar to what follows, and we will not repeat the argument. Since $\phi\mapsto I_T(\phi)$ is lower semicontinuous and $K_x=\{\phi\in K,\, \phi_0=x\}$ is compact, there exists a $\hat{\phi}\in K$ such that $\hat{\phi}_0=x$ and
$ I_T(\hat{\phi})= I_{T,x}(K)$.
Let now $\phi\in K_x$ be arbitrary. First assume that $I_T(\phi)<\infty$. Then there exists a
$\theta_\phi\in C^1([0,T];\R^d)$ such that
\[ I_T(\phi)\le \L(\phi,\theta_\phi)+\frac{\eps}{2}.\]
Since $\psi\mapsto \L(\psi,\theta_\phi)$ is continuous on $D([0,T];\R^d)$ equipped with the Skorokhod topology,
there exists  a neighbourhood $\V_{\phi,\theta_\phi}(\eps)$ of $\phi$ in $D([0,T];\R^d)$ such that for any
$\psi\in\V_{\phi,\theta_\phi}(\eps)$,
\[|\L(\phi,\theta_\phi)-\L(\psi,\theta_\phi)|\le \frac{\eps}{2}\,.\]
Now
\begin{align}
\P\left(Z^{N,x_N}\in\V_{\phi,\theta_\phi}(\eps)\right)&=\E\left({\mathbf1}_{Z^{N,x_N}\in\V_{\phi,\theta_\phi}(\eps)}\right)
\nonumber\\
&= e^{-N\L(\phi,\theta_\phi)}\E\left(e^{N\L(\phi,\theta_\phi)}{\mathbf1}_{Z^{N,x_N}\in\V_{\phi,\theta_\phi}(\eps)}\right)
\nonumber\\
&\le e^{-N[\L(\phi,\theta_\phi)-\frac{\eps}{2}]}\E\left(e^{N\L(Z^{N,x_N},\theta_\phi)}\right)\nonumber\\
&\le e^{-N[\L(\phi,\theta_\phi)-\frac{\eps}{2}]}\nonumber\\
&\le e^{-N I_{T}(\phi)+N\eps},\label{ineqUB}
\end{align}
where the before last inequality follows the fact that $N\L(Z^{N,x_N},\theta_\phi)=\log(\Xi^{N,\theta_\phi}_T)$
and $\E[\Xi^{N,\theta_\phi}_T]\le1$.

The second case is the one where $I_T(\phi)=+\infty$.
Then there exists $M>I_{T,x}(K)+1$ and
$\theta_\phi\in C^1([0,T];\R^d)$ such that $\L(\phi,\theta_\phi)>M+\eps$. From the same argument as above, we deduce that
\[ \P\left(Z^{N,x_N}\in\V_{\phi,\theta_\phi}(\eps)\right)\le e^{-NM}.\]

Let $K_x=\{\phi\in K,\, \phi_0=x\}$. Since $K_x\!\subset\!\bigcup_{\phi\in K, \phi_0=x} \V_{\phi,\theta_\phi}(\eps)$ and $K_x$ is compact, there exists $m=m(\eps)\ge1$ and
$\phi_1,\ldots,\phi_m\in K_x$ where we assume that $\phi_1=\hat{\phi}$, such that
\[ K_x\subset\bigcup_{i=1}^m\V_{\phi_i,\theta_{\phi_i}}(\eps)\,.\]
Now there exists a finite set of functions $\{\phi_{m+1},\ldots,\phi_{m+n}\}\subset K\backslash K_x$, such that
\[ K\subset \bigcup_{i=1}^{m+n}\V_{\phi_i,\theta_{\phi_i}}(\eps)\,.\]
We choose $\eps$ small enough for $i\ge m+1$ such that $x\not\in\V_{\phi,\theta_{\phi_i}}(\eps)$. Then for $N$ large enough,
$\P(Z^{N,x_N}\in \V_{\phi_i,\theta_{\phi_i}}(\eps))=0$ if $i\ge m+1$. Hence
\begin{align*}
\limsup_{N\to\infty}\frac{1}{N}\log\P(Z^{N,x_N}\in K)&\le
\limsup_{N\to\infty}\frac{1}{N}\log\left(\sum_{i=1}^{m+n}\P\left(Z^{N,x_N}\in \V_{\phi_i,\theta_{\phi_i}}(\eps)\right)\right)\\
&\le \max_{1\le i\le m}\limsup_{N\to\infty}\frac{1}{N}\log\P\left(Z^{N,x_N}\in \V_{\phi_i,\theta_{\phi_i}}(\eps)\right)\\
&\le-\inf_{1\le i\le m}I_T(\phi_i)+\eps\\
&\le-I_{T,x}(K)+\eps,
\end{align*}
where we have used \eqref{ineqUB} in the third inequality.
It remains to let $\eps\to0$.\epf

It remains to establish exponential tightness. Now we need to impose a growth condition on the $\beta_j$'s.
One natural assumption would be to assume that for some $C>0$, all $1\le j\le m$ and $x\in\R^d$, $\beta_j(t,x)\le C(1+|x|)$. However, this condition is not satisfied in most of our examples, because one of the $\beta_j$'s is quadratic.
We shall instead formulate an assumption which is satisfied in our epidemic models. We shall write $\mathds{1}$ for the vector in $\R^d$ whose coordinates are all equal to $1$, and we exploit the fact that for those $j$'s such that
$\beta_j$ is quadratic, $\langle h_j,\mathds{1}\rangle=0$.

\begin{description}
\item{(A.3)} We assume that for all starting points $x_N\in \Z_+^d/N$, $Z^{N,x_N}$ takes its values in $\R^d_+$ a.s., and moreover that there exists a $C_\beta>0$ such that
for any $0\le j\le k$ such that $\langle h_j,{\mathds{1}}\rangle\not=0$, $\beta_j(t,x)\le C_\beta(1+|x|),\ 0\le t\le T,\ x\in\R^d$.
\end{description}

We now prove
\begin{proposition}\label{exptight}
Assume that Conditions (A.1) and (A.3) are satisfied. Let $T>0$ and $x\in\R^d$ be given, as well as a sequence
$x_N\to x$ as $N\to\infty$, such that for all $N\ge1$, $x_N\in\Z_+^d/N$. Then or all $\xi>0$,
\[\lim_{\delta\downarrow0}\limsup_{N\to\infty}\frac{1}{N}\log
\P\left(\sup_{0\le s,t\le T,\ |t-s|\le\delta}\left|Z^{N,x_N}_t-Z^{N,x_N}_s\right|>\xi\right)=-\infty.\]
\end{proposition}
\bpf
$\xi>0$ and $T>0$ will be fixed throughout this proof. Consider the stopping time
\[ \sigma^{N,x_N}_R=\inf\{t\in[0,T],\ |Z^{N,x_N}_t|>R\}.\]
It is clear that
\begin{align*} \P\left(\sup_{\ |t-s|\le\delta}\left|Z^{N,x_N}_t-Z^{N,x_N}_s\right|>\xi\right)
&\le \P\left(\sup_{\ |t-s|\le\delta}\left|Z^{N,x_N}_{t\wedge\sigma^{N,x_N}_R}-Z^{N,x_N}_{s\wedge\sigma^{N,x_N}_R}\right|>\xi\right)\\&\quad+\P\left( \sigma^{N,x_N}_R<T\right).
\end{align*}
We first consider the first term of the above right-hand side. For that purpose, we divide $[0,T]$ into subintervals of length $\delta$, and let $\underline{i}(s)\le s<\overline{i}(s)$ denote the points of the grid nearest to $s$.
\begin{align*}
\mathbb{P}&\left[ \sup_{|s-t|\leq \delta} \left|Z^{N,x_N}_{t\wedge\sigma^{N,x_N}_R} - Z^{N,x_N}_{s\wedge\sigma^{N,x_N}_R} \right| >\xi  \right]\\
&=  \mathbb{P}\left[ \exists\, 0\le s<t\le T,\: t-s\leq \delta,\:
\left|Z^{N,x_N}_{t\wedge\sigma^{N,x_N}_R} - Z^{N,x_N}_{s\wedge\sigma^{N,x_N}_R}  \right| >\xi   \right]\\
&\le   \mathbb{P}\left[ \exists\, 0\le s<t\le T,\! t-s\leq \delta,\!
\left|Z^{N,x_N}_{t\wedge\sigma^{N,x_N}_R} - Z^{N,x_N}_{\underline{\mathit{i}}(s)}  \right|  \!+\!
\left|Z^{N,x_N}_{\underline{\mathit{i}}(s)} - Z^{N,x_N}_{s\wedge\sigma^{N,x_N}_R}  \right|  \!>\!\xi   \right] \\
&\le 2\left(\frac{T}{\delta} + 1 \right)\sup_{s\in [0,T]}\mathbb{P}\left[
\sup_{ t\in [s,s+2\delta[} \left|Z^{N,x_N}_{t\wedge\sigma^{N,x_N}_R} - Z^{N,x_N}_{s\wedge\sigma^{N,x_N}_R} \right| >\xi /2 \right].
\end{align*}
Let $\left\{\theta_i \:, \: 1\leq i\leq d \right\}$ (resp.\ $\left\{\theta_i \:, \: d+1\leq i\leq 2d \right\}$) denote  the standard  basis of $\mathbb{R}^d_+$  (resp.\ of $\mathbb{R}^d_-$).
Thus for every $\lambda>0$,  assuming w.l.o.g. that $|z|$ stands here for $\sup_{1\le i\le d}|z_i|$,
\smallskip

\resizebox{0.95\linewidth}{!}{
  \begin{minipage}{\linewidth}
\begin{align*}
\mathbb{P}&\left[
\sup_{ t\in [s,s+2\delta[ } \left|Z^{N,x_N}_{t\wedge\sigma^{N,x_N}_R} - Z^{N,x_N}_{s\wedge\sigma^{N,x_N}_R} \right| >\xi /2
\right] \\
&\leq  \sum_{i=1}^{2d}
\mathbb{P}\left[ \sup_{t\in [s,s+2\delta[}  \: \langle Z^{N,x_N} _{t} - Z^{N,x_N}_{s},\lambda\theta_i\rangle \: \: >\lambda\xi/2   \right]
\end{align*}
\end{minipage}
}

\resizebox{0.95\linewidth}{!}{
  \begin{minipage}{\linewidth}
\begin{align*}
&\leq   \sum_{i=1}^{2d}
\mathbb{P}\left[\sup_{t\in [s,s+2\delta[}  \: \mathbf{M}^{N,\lambda\theta_i}_{(s,t)\wedge\sigma^{N,x_N}_R} + \sum_{j=1}^{k}\int_{s\wedge\sigma^{N,x_N}_R}^{t\wedge\sigma^{N,x_N}_R} \langle h_j,\lambda\theta_i\rangle \beta_j\left(r,Z^{N,x_N}_r\right)dr \: > \lambda\xi/2  \right]\\
&\leq   \sum_{i=1}^{2d}
\mathbb{P}\left[\sup_{t\in [s,s+2\delta[}  \!\! \exp\!\left(\!\!N\mathbf{M}^{N,\lambda\theta_i}_{(s,t)\wedge\sigma^{N,x_N}_R} \!+\! N\sum_{j=1}^{k}\int_{s\wedge\sigma^{N,x_N}_R}^{t\wedge\sigma^{N,x_N}_R} \langle h_j,\lambda\theta_i\rangle \beta_j\left(r,Z^{N,x_N}_r\right)dr \!\right) \!>\! e^{N\lambda\xi/2}  \right]\\
&\leq   \sum_{i=1}^{2d}
\mathbb{P}\left[ \sup_{t\in [s,s+2\delta[}\! \Xi^{N,\lambda\theta_i}_{(s,t)\wedge\sigma^{N,x_N}_R} \!\!>\!\!\exp\!\left(\!N\lambda\xi/2 \!-\!N\!\sum_{i=1}^k\!
\left(e^{\langle h_j,\lambda\theta_i\rangle}-1\right)\!\int_{s\wedge\sigma^{N,x_N}_R}^{t\wedge\sigma^{N,x_N}_R}\!\!\!\beta_j(s,Z^{N,x_N}_r)dr\right)\right]  \\
&\leq   \sum_{i=1}^{2d}
\mathbb{P}\left[ \sup_{t\in [s,s+2\delta[}\: \Xi^{N,\lambda\theta_i}_{(s,t)\wedge\sigma^{N,x_N}_R} >\exp\left(N\lambda\xi/2 - 2\delta Nk\bar{\beta}_Re^{\lambda\bar{h}}\right)\right]  \\
&\leq  2d\exp\left(   -N  \lambda\xi/2  + 2\delta Nk\bar{\beta}_Re^{\lambda\bar{h}}\right),
\end{align*}
\end{minipage}
}
\smallskip

\noindent where $\bar{\beta}_R=\sup_{1\le j\le k}\sup_{0\le t\le T,\, |x|\le R}\beta_j(t,x)$.
Optimizing over $\lambda>0$ yields
\begin{align*}
\limsup_{N\to\infty}\frac{1}{N}\log\mathbb{P}\left[\sup_{|s-t|\leq \delta} \left|Z^{N,x_N}_{t\wedge\sigma^{N,x_N}_R} - Z^{N,x_N}_{s\wedge\sigma^{N,x_N}_R} \right| >\xi \right] \le
- \frac{\xi}{2\bar{h}}\left(\log\left(\frac{\xi\delta^{-1}}{4\bar{h}k\bar{\beta}_R}\right)-1\right).
\end{align*}
Consequently for any fixed $R>0$,
\[ \lim_{\delta\to0}\limsup_{N\to\infty}\frac{1}{N}\log\mathbb{P}\left[\sup_{|s-t|\leq \delta} \left|Z^{N,x_N}_{t\wedge\sigma^{N,x_N}_R} - Z^{N,x_N}_{s\wedge\sigma^{N,x_N}_R} \right| >\xi \right]=-\infty.\]
It remains to show that
\begin{equation}\label{-inf}
\lim_{R\to\infty}\limsup_{N\to\infty}\frac{1}{N}\log\P\left( \sigma^{N,x_N}_R<T\right)=-\infty\, .
\end{equation}
Combing the fact that  for all $t\leq T$
\[\displaystyle \sup_{s\leq t}|Z^{N,x_N}_s| \leq \sup_{s\leq t}|\langle Z^{N,x_N}_s, \mathds{1}\rangle |  \]
and that by Gronwall's Lemma \ref{le:Gron}, with $\bar{h}=\sup_{1\le j\le k}|h_j|$ and $C_\beta$ the constant from assumption (A.3),
\[\displaystyle \sup_{s\leq t}|\langle Z^{N,x_N}_s, \mathds{1}\rangle| \leq \left( |\langle x,\mathds{1}\rangle +k\bar{h}Ct + \sup_{s\leq t}|M^{N,\mathds{1}}_s| \right)e^{k\bar{h}C_{\beta}t},    \]
we deduce that
\begin{equation}\label{GronwallZNxMN1}
\sup_{s\leq t}|Z^{N,x_N}_s| \leq \left( |\langle x,\mathds{1}\rangle| + k\bar{h}Ct + \sup_{s\leq t}|M^{N,\mathds{1}}_s| \right)e^{k\bar{h}C_{\beta}t}.
\end{equation}
By It\^o's formula we have, with $\mathcal{M}^{N,\mathds{1}}_t$ a local martingale, and defining
$A^N_s:=1\vee(\sup_{0\le r\le s}|M^{N,\mathds{1}}_r|$,
\smallskip

\resizebox{0.95\linewidth}{!}{
  \begin{minipage}{\linewidth}
\begin{align*}
&\left(M^{N,\mathds{1}}_{t} \right)^{2N} \\
&=  \! N\!\!\!\!\!\sum_{j;\langle h_j,\mathds{1}\rangle\not=0}\!\!\int_0^t \!\!\!\beta_j\left(s,Z^{N,x_N}_{s}\right)\!\! \Bigg[ \!\!
\left(M^{N,\mathds{1}}_{s} + \frac{\langle h_j,\mathds{1}\rangle}{N} \right)^{2N}\!\!\!\!\!\!\! - \!\!\left(M^{N,\mathds{1}}_{s-} \right)^{2N}\!\! \!\!\!\!\!
- \!2N\!\!\left(M^{N,\mathds{1}}_{s} \right)^{2N-1}\!\!\frac{\langle h_j,\mathds{1}\rangle}{N}\Bigg] \!\ ds \!+\! \mathcal{M}^{N,\mathds{1}}_t
\end{align*}
\end{minipage}
}

\begin{align}
&\le NC_\beta\sum_j\frac{N(2N-1)}{N^2}\langle h_j,\mathds{1}\rangle^2\int_0^t(1+|Z^{N,x_N}_s|)
\left(|M^{N,\mathds{1}}_{s}|+\frac{\langle h_j,\mathds{1}\rangle}{N}\right)^{2N-2}ds+ \mathcal{M}^{N,\mathds{1}}_t
\nonumber\\
&\le NC_\beta C_T\left(1+\frac{\bar{h}}{N}\right)^{2N}\int_0^t\frac{1+|Z^{N,x_N}_s|}{A^N_s}(A^N_s)^{2N-1}ds+ \mathcal{M}^{N,\mathds{1}}_t\nonumber \\
&\le NC'_T\int_0^t(A^N_s)^{2N}ds+ \mathcal{M}^{N,\mathds{1}}_t,\label{pwN}
\end{align}
 where we have used \eqref{GronwallZNxMN1} and the inequality $a+b\le a(1+b)$ for $a\ge1$, $b\ge0$.
From Doob's inequality,
\begin{equation}\label{Doob}
 \E\left[\sup_{s\le t\wedge\sigma^{N,x_N}_R}(M^{N,\mathds{1}}_s)^{2N}\right]\le\left(\frac{2N}{2N-1}\right)^{2N}
\E\left[(M^{N,\mathds{1}}_{t\wedge\sigma^{N,x_N}_R})^{2N}\right].
\end{equation}
Since $\mathcal{M}^{N,\mathds{1}}_{t\wedge\sigma^{N,x_N}_R}$ is a martingale, we can take the expectation in the inequality \eqref{pwN} at time $t\wedge\sigma^{N,x_N}_R$, and deduce from the resulting inequality,
 \eqref{Doob} and $\sup_{N\ge1}\left(\frac{2N}{2N-1}\right)^{2N}<\infty$
\[\E\left[\sup_{s\le t\wedge\sigma^{N,x_N}_R}(M^{N,\mathds{1}}_s)^{2N}\right]\le NC''_T\int_0^t
\E\left[\left(A^N_{s\wedge\sigma^{N,x_N}_R}\right)^{2N}\right]ds.\]
Since for $a\ge0$, $(1\vee a)^{2N}\le 1+a^{2N}$, it follows that for all $0\le t\le T$,
\[ \E\left[\left(A^N_{t\wedge\sigma^{N,x_N}_R}\right)^{2N}\right]\le 1+NC''_T\int_0^t
\E\left[\left(A^N_{s\wedge\sigma^{N,x_N}_R}\right)^{2N}\right]ds.\]
Hence it follows from Gronwall's lemma that
\begin{equation}\label{mom2N}
\E\left[\sup_{t\le T\wedge\sigma^{N,x_N}_R}(M^{N,\mathds{1}}_t)^{2N}\right]\le \exp\left(C_T NT\right).
\end{equation}
For any $0<\kappa<R$, denoting
\[ \displaystyle C(R,\kappa) :=(R-\kappa)e^{-k\bar{h}C_T T} - |\langle x,\mathds{1}\rangle| -k\bar{h}C_TT\:, \]
we have
\begin{align*}
\limsup_{N\to +\infty}\frac{1}{N} \log\mathbb{P}&\left[ \sigma^{N,x_N}_R\leq T\right] \leq \displaystyle
\limsup_{N\to +\infty}\frac{1}{N} \log\mathbb{P}\left[ \sup_{t\leq T\wedge\sigma^{N,x_N}_R}|Z^{N,x}_t| >R -\kappa\right] \\
&\leq   \limsup_{N\to +\infty}\frac{1}{N}\log\mathbb{P}\left[ \sup_{t\leq T\wedge\sigma^{N,x_N}_R} |M^{N,\mathds{1}}_t| > C(R,\kappa)\right]\\
&\leq  \limsup_{N\to +\infty}\frac{1}{N}\log\mathbb{P}\left[ \sup_{t\leq T\wedge\sigma^{N,x_N}_R} (M^{N,\mathds{1}}_t)^{2N} >\left[ C(R,\kappa)\right]^{2N}\right]\\
&\leq  -2\log\left[C(R,\kappa)\right] + \limsup\limits_{N\to +\infty}\frac{1}{N} \log\mathbb{E}\left[  \sup_{t\leq T\wedge \sigma^{N,x_N}_R} \left(M^{N,\mathds{1}}_t\right)^{2N} \right] \\
&\le -2\log\left[C(R,\kappa)\right] + CT,
\end{align*}
where we have used \eqref{mom2N} for the last inequality.
We deduce \eqref{-inf} by letting $R$ tend to $+\infty$. \epf

We shall also need the following lemma, where we use the notation
\[ w'_{Z^N}(\delta)=\inf_{\{t_i\}}\max_{1\le i\le n}w_x([t_{i-1},t_i)),\]
with $w_x([t_{i-1},t_i))=\sup_{t_{i-1}\le s<t<t_i)}|x_t-x_s|$ and the infimum is taken over all sequences $0=t_0<t_1<\ldots<t_n=T$ satisfying
$\inf_{1\le i\le n}(t_i-t_{i-1})\ge\delta$.
\begin{lemma}\label{le:tightN}
If Conditions (A.1) and (A.3) are satisfied, then
for any $N\ge1$, $\rho>0$,
\[ \lim_{\delta\to0}\P(w'_{Z^N}(\delta)>\rho)=0.\]
\end{lemma}
\bpf Since the space $D([0,T];\R^d)$ is separable and complete, the law of $Z^N$ on this space is tight, see Theorem 1.3 in Billingsley \cite{PBI},
which implies the lemma, from Theorem 13.2 of the same reference. \epf

We can now deduce the following theorem from Proposition \ref{exptight} and Lemma \ref{le:tightN}.
\begin{theorem}\label{th:exptight} If Conditions (A.1) and (A.3) is satisfied, then
the sequence\linebreak $\{Z^{N,z_N},\, N\ge1\}$ is exponentially tight in $D([0,T];\R^d)$.
\end{theorem}
\bpf
Given $R>0$ and a sequence $\{\delta_\ell>0,\, \ell\ge1\}$ the following is a compact subset of $D([0,T];\R^d)$ (see Theorem 12.3 in Billingsley \cite{PBI}):
\[ K_{R,\{\delta_\ell\}}=\left\{x,\, \|x\|_T\le R\right\}\bigcap\bigcap_{\ell\ge1}\left\{x,\, w'_x(\delta_\ell)\le \ell^{-1}\right\}.\]
For any $\alpha>0$, we need to find $R_\alpha$ and $\{\delta^\alpha_\ell,\, \ell\ge1\}$ such that
\begin{equation}\label{eq:exptight}
\limsup_{N\to\infty}\frac{1}{N}\log\P\left(\left\{\|Z^{N,z_N}\|_T>R_\alpha\right\}\bigcup\bigcup_{\ell\ge1}\left\{w'_{Z^N}(\delta^\alpha_\ell)> \ell^{-1}\right\}\right)\le-\alpha.
\end{equation}
It is not hard to find $R_\alpha$ such that $\P(\|Z^N\|_T>R_\alpha)\le e^{-N\alpha}$, for 	all $N\ge1$. Since $w'_x(\delta)\le w_x(2\delta)$, it follows from Proposition \ref{exptight} that for each $\ell\ge1$, there exists a $\delta_\ell>0$ such that
\[\limsup_{N\to\infty}\frac{1}{N}\log\P\left(w'_{Z^N}(\delta_\ell)> \ell^{-1}\right)\le-(\alpha+\ell).\]
Consequently, there exists an $N_\ell$ such that for $N\ge N_\ell$,
\[\P\left(w'_{Z^N}(\delta_\ell)> \ell^{-1}\right)\le e^{-N(\alpha+\ell)}.\]
Combining this with Lemma \ref{le:tightN}, we deduce that there exists $0<\delta^\alpha_\ell\le\delta_\ell$ such that for all $N\ge1$,
\[\P\left(w'_{Z^N}(\delta^\alpha_\ell)> \ell^{-1}\right)\le e^{-N(\alpha+\ell)}.\]
It follows that for all $N\ge1$,
\begin{align*}
\P\left(\left\{\|Z^{N,z_N}\|_T>R_\alpha\right\}\bigcup\bigcup\left\{w'_{Z^N}(\delta^\alpha_\ell)> \ell^{-1}\right\}\right)&\le e^{-N\alpha}\sum_{\ell\ge0}e^{-N\ell}\\
&\le (1-e^{-N})^{-1}e^{-\alpha N},
\end{align*}
from which \eqref{eq:exptight} follows. \epf

It is not hard to see that a combination of the exact same arguments as used in the proofs of Theorem \ref{th:UB-compact}, Proposition \ref{exptight}
and Theorem \ref{th:exptight} yields the following result.
\begin{theorem}\label{th:UBunif}
Assume that assumptions (A.1) and (A.3) are satisfied. Then for any closed subset $F\subset D([0,T];\R^d)$ and any compact $\mathcal{K}\subset \R^d$, we have
\[\limsup_{N\to\infty}\frac{1}{N}\log\sup_{x\in \mathcal{K}}\P(Z^{N,x_N}\in F)\le -\inf_{x\in \mathcal{K}}I_{T,x}(F)\,.\]
\end{theorem}

\subsection{Time of extinction in the SIRS model}\label{sec:textSIRS}
We shall denote by $T^N_{\text{Ext}}$ the time of extinction of the disease, and we want to learn what large deviations can tell us about it. In order to simplify the presentation, we start with to the two most simple examples of the SIRS model and the SIS model. These are models with fixed population size $N$. We treat the SIRS model in this section, and the SIS model in the next one. In this section, we shall follow the arguments from Kratz and Pardoux \cite{KPI}, which itself follows closely the arguments in Dembo and Zeitouni \cite{DZI}.

The deterministic SIRS Model can be reduced to a $2$-dimensional ODE for the pair $(s(t),i(t))$ which reads
\begin{equation}\label{eq:SIRS}
\left\{
\begin{aligned}
i'(t)&= \lambda s(t)i(t)-\gamma i(t),\\
s'(t)&=-\lambda s(t)i(t)+\rho(1-s(t)-i(t)).
\end{aligned}
\right.
\end{equation}
This process lives in the compact set $A=A_{SIRS}=\{(x,y),\, 0\le x,y,\, x+y\le1\}$. Provided again $R_0=\frac{\lambda}{\gamma} >1$, there is a unique stable endemic equilibrium
$(i^*,s^*)=\left(\frac{\rho}{\lambda}\frac{\lambda-\gamma}{\rho+\gamma},\frac{\gamma}{\lambda}\right)\in A$,
while the disease free equilibrium $(1,0)$ is unstable. Here $h_1=\begin{pmatrix}-1\\1\end{pmatrix}$,
$\beta_1(x,y)=\lambda xy$, $h_2=\begin{pmatrix}0\\-1\end{pmatrix}$, $\beta_2(x,y)=\gamma y$,
$h_3=\begin{pmatrix}1\\0\end{pmatrix}$, $\beta_3(x,y)=\rho(1-x-y)$.

  The stochastic process $(I^N(t),S^N(t))$ may hit  $\{0\}\times[0,1]$, and then stays there for ever (this is how the disease goes extinct). On the other hand, if it hits $\partial A\backslash
  \{0\}\times[0,1]$, the process comes back to $\mathring{A}$. Similarly, starting form
   $\{0\}\times[0,1]$, the ODE stays there for ever (and converges to $(0,1)$),
  while starting from $\partial A\backslash\{0\}\times[0,1]$), it enters
  $\mathring{A}$ instantaneously. We thus define
  \[ T^N_{\text{Ext}}=\inf\{t\ge0,\, I^N(t)=0\}.\]
  Unfortunately, the theory of Large Deviations will not give us directly results on $T^N_{\text{Ext}}$, but rather on
  \[ T^N_\delta=\inf\{t\ge0,\, I^N(t)\le\delta\}, \text{ for any }\delta>0.\]   An ad hoc argument, which we shall present at the end, allows us to deduce the desired result concerning $T^N_{\text{Ext}}$.  We are interested in the exit time from $A_\delta:=\{(x,y)\in A,\, x\ge\delta\}$ through the boundary $\partial A_\delta:=\{(x,y)\in A,\, x=\delta\}$.

  We shall write $D_{T,A}:=D([0,T]; A)$. In order to formulate our results, we shall need the following notations (below $z$ stands for $(x,y)$)
  \begin{align*}
V(z,z',T)&= \inf_{\phi\in D_{T,A}, \phi_{0}=z, \phi_{T}=z'} I_{T}(\phi) \\
V(z,z')&= \inf_{T>0} V(z,z',T)  \\
 \overline{V}_\delta&= \inf_{z \in\partial A_\delta} V(z^{*},z),\\
 \overline{V}&= \inf_{z \in\{0\}\times[0,1]} V(z^{*},z).
\end{align*}
We want to prove the
\begin{theorem}\label{th:exitT}
Let $T^{N,z}_{\text{Ext}}$ denote the extinction time in the SIRS model starting from $z_N=\frac{[zN]}{N}$.
Given $\eta>0$, for all $z\in A$,
\begin{equation*}
\lim_{N\to\infty}\mathbb{P}\big(\exp\{N(\overline{V}-\eta)\}<T^{N,z}_{\text{Ext}}<\exp\{N(\overline{V}+\eta)\}\big)=1.
\end{equation*}
Moreover, for all $\eta>0$, $z\in A$ and $N$ large enough,
\begin{equation*}
\exp\{N(\overline{V}-\eta)\}\leq\mathbb{E}(T^{N,z}_{\text{Ext}})\leq\exp\{N(\overline{V}+\eta)\}.
\end{equation*}
\end{theorem}
We shall first establish
\begin{proposition}\label{pro:exitT}
Given $\eta>0$, for all $z\in \mathring{A}_\delta$,
\begin{equation*}
\lim_{N\to\infty}\mathbb{P}\big(\exp\{N(\overline{V}_\delta-\eta)\}<T^{N,z}_\delta<\exp\{N(\overline{V}_\delta+\eta)\}\big)=1.
\end{equation*}
Moreover, for all $\eta>0$, $z\in \mathring{A}_\delta$ and $N$ large enough,
\begin{equation*}
\exp\{N(\overline{V}_\delta-\eta)\}\leq\mathbb{E}(T^{N,z}_\delta)\leq\exp\{N(\overline{V}_\delta+\eta)\}.
\end{equation*}
\end{proposition}
Let us now formulate a set of assumptions which are satisfied in our case, under which we will prove
Proposition \ref{pro:exitT}. For that sake, we shall rewrite the ODE \eqref{eq:SIRS} as
\begin{align}\label{ODEsirs}
 \frac{dz_t}{dt}=b(z_t),\ z_0=z.
 \end{align}
\begin{assumption}\label{AssExit}$\mbox{ }$
\begin{enumerate}
\item[(E1)]
$z^*$ is the only stable equilibrium point of~\eqref{ODEsirs} in $A_\delta$ and the solution $z_t^x$ of~\eqref{ODEsirs} satisfies, for all $z_0=z \in A_\delta$,
\[
z_t^z \in \mathring{A_\delta} \text{ for all } t>0 \text{ and } \lim_{t\rightarrow \infty} z_t^z=z^*.
\]
\item[(E2)]
$\bar{V} <\infty$.
\item[(E3)]
For all $\rho>0$ there exist constants $T(\rho)$, $\epsilon(\rho)>0$ with $T(\rho), \epsilon(\rho) \downarrow 0$ as $\rho \downarrow 0$ such that for all $z \in \partial A_\delta \cup \{z^*\}$ and all $x,y \in \overline{B(z,\rho)} \cap A$ there exists
\[
\phi=\phi(\rho,x,y):[0,T(\rho)] \mapsto A  \text{ with } \phi_0=x, \phi_{T(\rho)}=y \text{ and } I_{T(\rho)}(\phi)<\eps(\rho).
\]
\item[(E4)]
For all $z\in {\partial A_\delta}$ there exists an $\eta_0>0$ such that for all $\eta<\eta_0$ there exists a $\tilde z=\tilde z(\eta)\in A\backslash A_\delta$ with $|z-\tilde z|>\eta$.
\end{enumerate}
\end{assumption}
Note that the conditions $(E1)$ and $(E4)$ would not be satisfied if we replace $A_\delta$ by $A$.

The proof of Proposition \ref{pro:exitT} relies upon the following sequence of lemmas, whose proofs will be given below, after the proof of the proposition.
\begin{lemma}\label{le:1}
For any $\eps>0$, there exists a $\rho_0>0$ such that for all $\rho<\rho_0$,
\[ \sup_{z\in\partial A_\delta\cup\{z^*\}}\sup_{|z'-z|\vee|z''-z|\le\rho}\inf_{0\le T\le 1}V(z',z'',T)<\eps.\]
\end{lemma}

\begin{lemma}\label{le:2}
For any $\eta>0$, there exists a $\rho_0>0$ such that for all $\rho<\rho_0$, there exists a $T_0<\infty$ such that
\[ \liminf_{N\to\infty}\frac{1}{N}\log\inf_{|z-z^*|\le\rho}\P(T^{N,z}_\delta\le T_0)\ge -(\bar{V}+\eta).\]
\end{lemma}

Let us define for some $\rho>0$ small enough,  $B_\rho:=\overline{B(z^*,\rho)}$ and
\[ \sigma^N_\rho=\inf\{t\ge0,\, Z^N_t\in B_\rho\cup\{z,\, z_1\le\delta\}\}.\]
\begin{lemma}\label{le:3}
If $\rho>0$ is such that $B_\rho\subset\mathring{A_\delta}$, then
\[ \lim_{t\to\infty}\limsup_{N\to\infty}\frac{1}{N}\log\sup_{x\in A_\delta}\P(\sigma^{N,z}_\rho>t)=-\infty.\]
\end{lemma}

\begin{lemma}\label{le:4}
Let $C$ be a closed subset of $A\backslash \mathring{A_\delta}$. Then
\[ \lim_{\rho\to0}\limsup_{N\to\infty}\frac{1}{N}\log\sup_{2\rho\le |z-z^*|\le 3\rho}\P(Z^{N,z}_{\sigma^N_\rho}\in C)
\le - \inf_{z'\in C}V(z^*,z').\]
\end{lemma}

\begin{lemma}\label{le:5}
If $\rho>0$ is such that $B_\rho\subset\mathring{A_\delta}$ and $z\in\mathring{A_\delta}$,
\[ \lim_{N\to\infty}\P(Z^{N,z}_{\sigma^N_\rho}\in B_\rho)=1.\]
\end{lemma}

\begin{lemma}\label{le:6}
For all $\rho,c>0$, there exists a constant $T=T(c,\rho)<\infty$ such that
\[ \limsup_{N\to\infty}\frac{1}{N}\log\sup_{z\in A_\delta}\P(\sup_{0\le t\le T}|Z^{N,z}_t-z|\ge\rho)\le -c.\]
\end{lemma}

We first give the
\smallskip

\begin{proof}[Proof of Proposition \ref{pro:exitT}]
{\sc Step 1: upper bound of $T^N_\delta$}
We choose $\eta=\eps/2$, and $\rho$, $T_0$ as in Lemma \ref{le:2}.
By Lemma \ref{le:3},  for any arbitrarily fixed $a>0$, there exists a $T_1$ such that
\[ \limsup_{N\to\infty}\frac{1}{N}\log\sup_{z\in A_\delta}\P(\sigma^{N,z}_\rho>T_1)<-2a<0.\]
Let $T=T_0+T_1$. There exists an $N_0\ge1$ such that for all $N\ge N_0$,
\begin{align}\label{qinf}
q:=\inf_{z\in A_\delta}\P(T^{N,z}_\delta\le T)&\ge\inf_{z\in A_\delta}\P(\sigma^{N,z}_\rho\le T_1)
\inf_{z\in B_\rho}\P(T^{N,z}_\delta\le T_0)\nonumber\\
&\ge e^{-N(\bar{V}_\delta+\eta)},
\end{align}
since the second factor is bounded from below by say $e^{-N(\bar{V}_\delta+\eta/2)}$ from Lemma \ref{le:2},
and from the previous estimate, we deduce that for $N$ large enough,
\begin{align*}
 \inf_{z\in A_\delta}\P(\sigma^{N,z}_\rho\le T_1)&=1-\sup_{z\in A_\delta}\P(\sigma^{N,z}_\rho>T_1)\\
&\ge1-e^{-Na}\\
&\ge e^{-N\eta/2}\,.
\end{align*}
Next, by the strong Markov property,
\begin{align*}
\P(T^{N,z}_\delta>(k+1)T)&=[1-\P(T^{N,z}_\delta\le(k+1)T|T^{N,z}_\delta>kT)]\P(T^{N,z}_\delta>kT)\\
&\le (1-q)\P(T^{N,z}_\delta>kT).
\end{align*}
Iterating, we get
\[ \sup_{z\in A_\delta}\P(T^{N,z}_\delta>kT)\le(1-q)^k.\]
Therefore
\[\sup_{z\in A_\delta}\E(T^{N,z}_\delta)\le T[1+\sum_{k=1}^\infty\sup_{z\in A_\delta}\P(T^{N,z}_\delta>kT)]
\le T\sum_{k=0}^\infty(1-q)^k=\frac{T}{q},\]
so from \eqref{qinf},
\begin{align}\label{upExp}
\sup_{x\in A_\delta}\E[T^{N,z}_\delta]\le Te^{N(\bar{V}_\delta+\eta)},
\end{align}
and the upper bound for $\E[T^{N,z}_\delta]$ follows. From Chebycheff,
\[ \P(T^{N,z}_\delta\ge e^{N(\bar{V}_\delta+\eps)})\le e^{-N(\bar{V}_\delta+\eps)}\E[T^{N,z}_\delta]\le Te^{-N\eps/2},\]
which tends to $0$ as $N\to\infty$, hence the upper bound for $T^{N}_\delta$.

\noindent{\sc Step 2: lower bound of $T^N_\delta$}
Let $\rho>0$ be small enough such that $B_{2\rho}:=B(z^*,2\rho)\subset \mathring{A_\delta}$. We define a sequence of stopping times as follows. $\theta_0=0$ and for $m\ge0$,
\begin{align*}
\tau_m&=\inf\{t\ge\theta_m,\, Z^{N}_t\in B_\rho\cup\{z,\, z_1\le\delta\}\},\\
\theta_{m+1}&=\inf\{t>\tau_m,\, Z^{N}_t\in (B_{2\rho})^c\},
\end{align*}
with the convention that $\theta_{m+1}=\infty$ in case $Z^{N}_{\tau_m}\in\{z,\, z_1\le\delta\}$.

In case $\bar{V}_\delta=0$, the lower bound is an easy consequence of Lemmas \ref{le:5} and \ref{le:6}. So we assume from now on that $\bar{V}_\delta>0$ and fix $\eps>0$ arbitrarily small. Since $\{z,\, z_1\le\delta\}$ is a closed set,
from Lemma \ref{le:4}, for $\rho>0$ small enough,
\[ \limsup_{N\to\infty}\frac{1}{N}\log\sup_{2\rho\le|z-z^*|\le3\rho}\P(Z^{N,z}_{\sigma^N_\rho}\in \{z,\, z_1\le\delta\})
\le - \bar{V}_\delta+\frac{\eps}{3}.\]
Now with $c=\bar{V}_\delta$, we let $T_0=T(c,\rho)$ be as in Lemma \ref{le:6}. Then there exists an $N_0$ such that for $N\ge N_0$, and all $m\ge1$,
\[\sup_{z\in A_\delta}\P(T^{N,z}_\delta=\tau_m)\le
\sup_{2\rho\le|z-z^*|\le3\rho}\P(Z^{N,z}_{\sigma^N_\rho}\in\{z,\, z_1\le\delta\})\le e^{-N(\bar{V}_\delta-\eps/2)},\]
while
\[ \sup_{z\in A_\delta}\P_z(\theta_m-\tau_{m-1}\le T_0)\le\sup_{z\in A_\delta}\P(\sup_{0\le t\le T_0}|Z^{N,z}_t-z|\ge\rho)
\le e^{-N(\bar{V}_\delta-\eps/2)}.\]
The event $\{T^{N}_\delta\le kT_0\}$ implies that either one of the first $k+1$ events $\{T^{N}_\delta=\tau_m\}$
occurs, or else at least one of the first $k$ excursions $[\tau_m,\tau_{m+1}]$ away from $B_\rho$ is of length at most $T_0$.
Consequently, from the two preceding estimates,
\begin{align*}
\P(T^{N}_\delta\le kT_0)&\le \sum_{m=0}^k\P(T^{N}_\delta=\tau_m)+\P(\min_{1\le m\le k}(\theta_m-\tau_{m-1})\le T_0)\\
&\le \P(T^{N}_\delta=\tau_0)+2ke^{-N(\bar{V}_\delta-\eps/2)}.
\end{align*}
Choosing now $k=[T_0^{-1}e^{N(\bar{V}_\delta-\eps)}]+1$ yields
\[ \P(T^{N}_\delta\le e^{N(\bar{V}_\delta-\eps)})\le\P(Z^{N}_{\sigma^N_\rho}\not\in B_\rho)+3T_0^{-1}e^{-N\eps/2}.\]                       By Lemma \ref{le:5}, the right-hand side tends to $0$ as $N\to\infty$. We have completed the proof of the first statement in Proposition \ref{pro:exitT}. This result  combined with Chebycheff's inequality and \eqref{upExp}
yields the second result. \epf

We now turn to the proofs of the lemmas.
\smallskip

\begin{proof}[Proof of Lemma \ref{le:1}]
This lemma is a direct consequence of the assumption $(E3)$. \epf

\begin{proof}[Proof of Lemma \ref{le:2}]
We make use of Lemma \ref{le:1} with $\eps=\eta/4$ and choose $\rho<\rho_0$. Let $z\in B_\rho$. There exists a continuous path  $\psi^z$ such that $\psi^z_0=z$, $\psi^z_{t_z} =z^*$ for some $t_z\le1$ and $I_{t_z}(\psi^z)\le\eta/4$.
From assumption $(E2)$, there exists a continuous path $\phi\in C([0,T_1];A)$ such that $\phi_0=z^*$, $\phi_{T_1}=z'
\in \partial A_\delta$, and $I_{T_1}(\phi)\le \bar{V}+\eta/4$. From Lemma \ref{le:1}, there exists a continuous path
$\tilde{\psi}$ such that $\tilde{\psi}_0=z'$ and $\tilde{\psi}_{s_{z'}}=z''\in A\backslash A_\delta$, with $s_{z'}\le 1$,
$I_{s_{z'}}(\tilde{\psi}) \le\eta/4$ and $d(z'',A_\delta)=\Delta>0$, where $\Delta<\delta$. Finally let $\{\xi_t,\, 0\le t\le 2-t_z-s_{z'}\}$ be a solution of \eqref{ODEsirs} starting from $\xi_0=z''$. From Proposition \ref{I=0}, $I_{}(\xi)=0$. Concatenating the paths
$\psi^z$, $\phi$, $\tilde{\psi}$ and $\xi$, we obtain a path $\phi^z\in C([0,T_0];A)$ (with $T_0=T_1+2$) starting from $z$, with $I_{T_0}(\phi^z)\le\bar{V}+3\eta/4$. Let now
\[ \Psi=\bigcup_{z\in B_\rho}\{\psi\in D([0,T_0];A),\, \|\psi-\phi^z\|_{T_0}<\Delta/2\}.\]
$\Psi$ is an open subset of $D([0,T_0];A)$, such that
 $B_\rho$ is adapted to $\Psi$ in the sense of Definition \ref{def:adapted}. Hence we can make use of Theorem \ref{th:LBunif}, hence
\begin{align*}
\liminf_{N\to\infty}\frac{1}{N}\log\inf_{z\in B_\rho}\P(Z^{N,z}\in\Psi)&\ge-\sup_{z\in B_\rho}\inf_{\phi\in\Psi, \phi_0=z}I_{T_0}(\phi)\\
&\ge-\sup_{z\in B_\rho}I_{T_0}(\phi^z)\\
&>-(\bar{V}+\eta).
\end{align*}
The results follows from this and $\{Z^{N}\in\Psi\}\subset \{T^N_\delta\le T_0\}$. \epf

\begin{proof}[Proof of Lemma \ref{le:3}]
Since $\sigma^{N,z}_\rho=0$ if $z\in B_\rho$, it suffices to restrict ourselves to $z\in A_\delta\backslash B_\rho$.
For each $t>0$, we define the closed set
\[ \Psi_t:=\{\phi\in D([0,t];A),\, \phi_s\in\overline{A_\delta\backslash B_\rho}\text{ for all }0\le s\le t\},\]
so that $\{\sigma^{N,z}_\rho>t\}\subset\{Z^{N,z}\in\Psi_t\}$. Hence by Theorem \ref{th:UBunif},
\[\limsup_{N\to\infty}\frac{1}{N}\log\sup_{z\in\overline{A_\delta\backslash B_\rho}}\P(\sigma^{N,z}_\rho>t)
\le -\inf_{\phi\in\Psi_t}I_t(\phi).\]
It then suffices to show that
\begin{equation}\label{hyp}
\inf_{\phi\in\Psi_t}I_t(\phi)\to\infty\text{ as }t\to\infty.
\end{equation}
Starting from any $z\in\overline{A_\delta\backslash B_\rho}$, the solution $z^z_t$ of \eqref{ODEsirs} hits
$\mathring{B}_{\rho/2}$ in finite time $T_z$ which is upper semicontinuous (recall Definition \ref{de:sc}) in $z$, so by compactness $T:=\sup_{z\in\overline{A_\delta\backslash B_\rho}}T_z<\infty$, and $z_\cdot\not\in\Psi_T$ as soon as $z_\cdot$ solves \eqref{ODEsirs}.

Now if \eqref{hyp} does not hold, there would exist $M>0$ and for each $n\ge1$ $\phi_n\in\Psi_{nT}$ such that
$I_{nT}(\phi_n)\le M$ or all $n\ge1$. Now if $\phi_{n,k}(t)=\phi_n(kT+t)$, $0\le t\le T$, $0\le k\le n-1$, we have that
\[n\min_{0\le k\le n-1}I_T(\phi_{n,k})\le\sum_{k=0}^{n-1}I_T(\phi_{n,k})=I_{nT}(\phi_n)\le M.\]
Hence we would produce a sequence $\psi_n\in\Psi_T$ such that  $I_T(\psi_n)\to0$ as $n\to\infty$.
From Theorem \ref{th:goodratef}, the sequence $\psi_n$ belongs to a compact set, and $I_T$ is lower semicontinuous (recall Definition \ref{de:sc}), so that
along a subsequence, $\psi_n\to\psi^*$, where $\psi^*\in\Psi_T$ and $I_T(\psi^*)\le\liminf_nI_T(\psi_n)=0$, and those two last statements are contradictory from Proposition \ref{I=0} and the fact that $\Psi_T$ contains no solution of \eqref{ODEsirs}. \epf

\begin{proof}[Proof of Lemma \ref{le:4}]
We need only consider the case $\inf_{z\in C}V(z^*,z)>0$, since in the other case the result is trivial.
So we can choose $\eps>0$ such that
\[ V^\eps_C:=\left(\inf_{z\in C}V(z^*,z)-\eps\right)\wedge\eps^{-1}>0.\]
By Lemma \ref{le:1}, there exists a $\rho_0>0$ such that for all $0<\rho<\rho_0$,
\[\sup_{z\in\overline{B_{3\rho}\backslash B_{2\rho}}}V(z^*,z)<\eps,\]
hence
\[\inf_{z'\in\overline{B_{3\rho}\backslash B_{2\rho}},z\in C}V(z',z)\ge \inf_{z\in C}V(z^*,z)-
\sup_{z'\in\overline{B_{3\rho}\backslash B_{2\rho}}}V(z^*,z')>V^\eps_C.\]
For $T>0$, consider the closed set $\Phi^T\subset D([0,T]; A)$ defined as
\[\Phi^T=\{\phi\in D([0,T]; A),\, \phi_t\in C\text{ for some }0\le t\le T\}.\]
For $z'\in\overline{B_{3\rho}\backslash B_{2\rho}}$,
\begin{equation}\label{upbound}
\P(Z^{N,z'}_{\sigma^N_\rho}\in C)\le \P(\sigma^{N,z'}_\rho>T)+\P(Z^{N,z'}\in\Phi^T).
\end{equation}
We next bound from above the two terms of the last right-hand side. Concerning the second term,
\[\inf_{\phi\in\Phi^T,\phi_0\in\overline{B_{3\rho}\backslash B_{2\rho}}}I_T(\phi)\ge
\inf_{z'\in\overline{B_{3\rho}\backslash B_{2\rho}},z\in C}V(z',z)\ge V^\eps_C.\]
Hence from Theorem \ref{th:UBunif},
\[\limsup_{N\to\infty}\frac{1}{N}\log\sup_{z'\in\overline{B_{3\rho}\backslash B_{2\rho}}}\P(Z^{N,z'}\in\Phi^T)
\le -V^\eps_C.\]
For the first term, we deduce from Lemma \ref{le:3} that for some $T_0>0$, all $T\ge T_0$,
\[\limsup_{N\to\infty}\frac{1}{N}\log\sup_{z'\in\overline{B_{3\rho}\backslash B_{2\rho}}}
\P(\sigma^{N,z'}_\rho>T)<-V^\eps_C.\]
\eqref{upbound} together with the last two estimates produces an inequality which, after letting $\eps\to0$,
yields the result. \epf

\begin{proof}[Proof of Lemma \ref{le:5}]
$z^z_t$ denoting the solution of \eqref{ODEsirs} starting from $z\in\mathring{A_\delta}$, let $T_\rho=\inf\{t>0,\, z^z_t\in B_{\rho/2}\}$. From $(E1)$ it follows that
$T_\rho<\infty$ and $\Delta:=\inf_{0\le t\le T_\rho}d(z^z_t,\partial A_\delta)>0$. Consequently
\[ \P\left(Z^{N,z}_{\sigma^N_\rho}\not\in B_\rho\right)\le \P\left(\sup_{0\le t\le T_\rho}|Z^{N,z}_t-z^z_t|\ge\frac{\Delta\wedge\rho}{2}\right),\]
which tends to $0$ as $N\to\infty$ from Theorem \ref{th:LLN}. Q.E.D. \epf

\begin{proof}[Proof of Lemma \ref{le:6}]
Let $\rho, c>0$ be fixed. For $T>0$, $N\ge1$ and $z\in A_\delta$,
\begin{align*}
\P\left(\sup_{0\le t\le T}|Z^{N,z}_t-z|\ge\rho\right)&=
\P\left(\sup_{0\le t\le T}\left|\sum_jh_jP_j\left(N\int_0^t\beta_j(Z^{,x}_s)ds\right)\right|\ge\rho\right)\\
&\le\P\left(\sum_jP_j(N\bar{\beta}T)\ge N\rho\bar{h}^{-1}\right)\\
&\le k\P(P(N\bar{\beta}T)\ge N\rho\bar{h}^{-1}k^{-1}).
\end{align*}
Now from Cram\'er's Theorem \ref{th:cramer}
\[\limsup_{N\to\infty}\frac{1}{N}\log\sup_{z\in A_\delta}\P(\sup_{0\le t\le T}|Z^{N,z}_t-z|\ge\rho)\le
-\frac{\rho}{\bar{h}k}\log\left(\frac{\rho}{\bar{h}k\bar{\beta} T}\right)+\frac{\rho}{\bar{h}k}-\bar{\beta}T,\]
and the absolute value of the right-hand side can be made arbitrarily large by choosing $T$ arbitrarily small.
\epf

It remains finally to turn to the
\smallskip

\begin{proof}[Proof of Theorem \ref{th:exitT}]
Since $V_\delta\uparrow V$ as $\delta\downarrow0$, it is clear that the lower bounds for $T^N_{\text{Ext}}$ and its expectation follow from Proposition \ref{pro:exitT}. It remains to establish the upper bound. Analyzing carefully the proof of the upper bound, we notice that the key step, which relies upon Lemmas \ref{le:2} and \ref{le:3} whose proof do not extend to our new situation, is the derivation of the inequality \eqref{qinf}. The upper bound both for
the time of exit and its expectation are a direct consequence of \eqref{qinf}, without any further reference to those assumptions which are not valid any more.
We fix $\eta>0$. Let $t>0$ be arbitrary. From Lemma \ref{le:ext-bp} below, if
$c_t:=\log\left(\frac{\lambda-\gamma e^{(\gamma-\lambda)t}}{\gamma-\gamma e^{(\gamma-\lambda)t}}\right)$,
\begin{equation}\label{extinct}
 \inf_{z\in A\backslash A_\delta}\P(T^{N,z}_{Ext}\le t)\ge e^{-\lceil N\delta\rceil c_t}\ge e^{-N(\delta+N_0^{-1})c_t},
 \end{equation}
provided $N\ge N_0$. Choose $N_0$ large enough and $\delta>0$ small enough such that $(\delta+N_0^{-1})c_t\le \eta/2$.
From \eqref{qinf}, there exists a $T_\delta>0$ such that, possibly increasing $N_0$ if necessary, if $N\ge N_0$,
\begin{equation}\label{qinfb}
 \inf_{z\in A_\delta}\P(T^{N,z}_\delta\le T_\delta)\ge e^{-N(\bar{V}_\delta+\eta/2)}\ge e^{-N(\bar{V}+\eta/2)}\,.
 \end{equation}
 We deduce from \eqref{extinct}, \eqref{qinfb} and the strong Markov property that, with $T=T_\delta +t$,
 \[ \inf_{z\in A}\P(T^{N,z}_{Ext}\le T)\ge e^{-N(\bar{V}+\eta)}\,,\]
 which is the wished extension of \eqref{qinf}.
 \epf

\begin{lemma}\label{le:ext-bp}
For any $t>0$, if $c_t:=\log\left(\frac{\lambda-\gamma e^{(\gamma-\lambda)t}}{\gamma-\gamma e^{(\gamma-\lambda)t}}\right)$,
\[ \inf_{z\in A\backslash A_\delta}\P(T^{N,z}_{\text{Ext}}<t)
\ge\exp\left\{-\lceil N\delta\rceil c_t\right\}.\]
\end{lemma}
\bpf
Since $z\in A\backslash A_\delta$ implies that $z_1\le \delta$, the first component of the process $NZ^{N,x}(t)$ is dominated by the  process
\[\lceil N\delta\rceil+P_1\left(N\lambda\int_0^tZ^{N,z}_1(s)ds\right)-P_1\left(N\gamma\int_0^tZ^{N,z}_1(s)ds\right),
\] which is a continuous time binary branching process with birth rate $\lambda$ and death rate $\gamma$. This process goes extinct before time $t$ with probability
\[ \left(\frac{\gamma-\gamma e^{(\gamma-\lambda)t}}{\lambda-\gamma e^{(\gamma-\lambda)t}}\right)^{\lceil N \delta\rceil},\]
as can be seen by combining formula (1) from section III.4 of Athreya and Ney \cite{ANI} with the formula in section 5 for $F(0,t)$ in the birth and death case.
The result follows readily.
\epf

We shall need below the following additional results.
\begin{proposition}\label{exitLocdelta}
Under the assumptions of Proposition \ref{pro:exitT}, if $C\subset\partial A_\delta$ is a closed set such that $V_C:=\inf_{z\in C} V(z^*,z)>\bar{V}_\delta$,
then for any $z\in \mathring{A}_\delta$, all $\eps>0$ small enough
\[ \lim_{N\to\infty}\P(d(Z^{N,z_N}_{T^N_\delta}, C)\le\eps)=0.\]
\end{proposition}
\bpf
Fix $\eta<(V_C-\bar{V}_\delta)/3$. From Lemma \ref{le:4}, for  $\eps>O$ small enough,
there exists a $\rho>0$ small enough and $N_0$ large enough such that for all $N\ge N_0$,
\[ \sup_{2\rho\le|z-z^*|\le3\rho}\P(d(Z^{N,z_N}_{\sigma^N_\rho},C)\le\eps)\le e^{-N(V_C-\eta)}.\]
Let $c=V_C-\eta$ and $T_0=T(c,\rho)$ given by Lemma \ref{le:6}. Then, increasing $N_0$ if necessary, we deduce from that Lemma that
for any $N\ge N_0$, $\ell\ge1$,
\[ \P(\tau_\ell\le\ell T_0)\le \ell\sup_{z\in A_\delta}\P\left(\sup_{0\le t\le T_0}|Z^{N,z_N}_t-z|\ge\rho\right)\le\ell e^{-N(V_C-\eta)}.\]
For all $z\in B_\rho$, $\ell\ge1$,
\begin{align*}
\P(d&(Z^{N,z_N}_{T^N_\delta}, C)\le\eps)\\&\le\P(T^{N,z}_\delta>\tau_\ell)
+\sum_{m=1}^\ell\P(T^{N,z}_\delta>\tau_{m-1})\P(d(Z^{N,z_N}_{\tau_m},C)\le\eps|T^{N,z}_\delta>\tau_{m-1})\\
&\le\P(T^{N,z}_\delta>\ell T_0)+\P(\tau_\ell\le \ell T_0)\\&\qquad
+\sum_{m=1}^\ell \P(T^{N,z}_\delta>\tau_{m-1})\E[\P(d(Z^{N,Z^N_{\theta_m}}_{\sigma^N_\rho},C)\le\eps|T^{N,z}_\delta>\tau_{m-1}]\\
&\le \P(T^{N,z}_\delta>\ell T_0)+\P(\tau_\ell\le \ell T_0)+\ell\sup_{2\rho\le|z-z^*|\le3\rho}\P(d(Z^{N,z_N}_{\sigma^N_\rho},C)\le\eps)\\
&\le \P(T^{N,z}_\delta>\ell T_0)+ 2\ell e^{-N(V_C-\eta)}.
\end{align*}
Increasing further $N_0$ if necessary, we have that \eqref{upExp} holds for some $T>0$ and all $N\ge N_0$. We choose
$\ell=\left[e^{N(\bar{V}_\delta+2\eta)}\right]$, hence from our choice of $\eta$,
\[\limsup_{N\to\infty}\sup_{z\in B_\rho}\P(d(Z^{N,z_N}_{T^N_\delta}, C)\le\eps)\le\limsup_{N\to\infty}\left(\frac{T}{\ell T_0}e^{N(\bar{V}_\delta+\eta)}
+2\ell e^{-N(V_C-\eta)}\right)=0.\]
It remains to combine Lemma \ref{le:5} and the inequality
\[ \P(d(Z^{N,z_N}_{T^N_\delta}, C)\le\eps)\le\P(Z^{N,z_N}_{\sigma^N_\rho}\not\in B_\rho)+\sup_{y\in B_\rho}\P(d(Z^{N,y_N}_{T^N_\delta}, C)\le\eps).\]
\epf

The proof of the next important result is a bit lengthy, and we refer to Pardoux and Samegni-Kepgnou \cite{PS3I} for it.
\begin{corollary}\label{exitLoc}
If $C\subset\{z,\, z_1=0\}$ is such that $V_C:=\inf_{z\in C}V(z^*,z)>\bar{V}$, then for any $z\in\mathring{A}$,
\[ \lim_{N\to\infty}\P(Z^{N,z_N}_{T^N_{\text{Ext}}}\in C)=0.\]
\end{corollary}
%
%
%
%

\subsection{Time of extinction in the SIS model}
While the above results are rather precise, it is frustrating that it does not seem possible to express the important constant $\overline{V}$
explicitly in terms of the few constants of the model. One can only do a numerical evaluation of $\overline{V}$. We now simplify the problem, and consider the SIS model, where when an infectious individual cures, he immediately becomes susceptible again: there is no immunity. The advantage of this simplified model is that it can be written in dimension one and, as we shall see now, we can deduce from the Pontryagin maximum principle, see
Section \ref{se:Pont} in the Appendix, a very simple explicit formula for $\overline{V}$.

The deterministic SIS model can be reduced to the following one--dimen\-si\-onal equation for the proportion of infected individuals
\[ \dot{x}_t=\lambda x_t(1-x_t)-\gamma x_t.\]
Here the process lives in the interval $A_{SIS}=[0,1]$. Provided $R_0=\frac{\lambda}{\gamma} >1$, there is a unique stable endemic equilibrium $x^*=1-\frac{\gamma}{\lambda}\in (0,1)$, while the disease free equilibrium $x^0=0$ is unstable. Here $h_1=1$, $\beta_1(x)=\lambda x(1-x)$, $h_2=-1$, $\beta_2(x)=\gamma x$.

We assume that $\lambda>\gamma$, i.e.\  $R_0>1$. As the reader can easily verify,
Theorem \ref{th:exitT} applies to this situation, and now $\overline{V}$ is the minimal value of the following control problem. With the notations of
Section \ref{se:Pont} below, we are in the situation $d=1$, $k=2$, $\beta_1(x)=\lambda x(1-x)$, $\beta_2(x)=\gamma x$,
$B=\begin{pmatrix} 1 & -1\end{pmatrix}$. The identity \eqref{eq:Ha} reads here
\[ \lambda x_t(1-x_t)(1-e^{p_t})+\gamma x_t(1-e^{-p_t})=0.\]
Hence either $p_t=0$, or else $p_t=\log\frac{\gamma}{\lambda(1-x_t)}$. It is easy to convince oneself that $p_t=0$ does not produce a control which does the wished job. Hence  $p_t=\log\frac{\gamma}{\lambda(1-x_t)}$, $\hat{u}_1(t)=e^{p_t}\beta_1(x_t)=\gamma x_t$, $\hat{u}_2(t)=e^{-p_t}\beta_2(x_t)
=\lambda x_t(1-x_t)$. The optimal trajectory reads
\begin{equation}\label{eq:optx}
 \dot{x}_t=\gamma x_t-\lambda x_t(1-x_t).
 \end{equation}
From the right-hand side of the identity \eqref{eq:IC},
\begin{align*}
\overline{V}&=\int_0^{\hat{T}}\left[\gamma x_t-\lambda x_t(1-x_t)\right]\log\frac{\gamma}{\lambda (1-x_t)}dt\\
&= \int_0^{\hat{T}}\log\frac{\gamma}{\lambda (1-x_t)}\dot{x}_tdt\\
&= \int_{0}^{\frac{\lambda-\gamma}{\lambda}} \log\frac{\gamma}{\lambda (1-x)} dx\\
&= \log\frac{\lambda}{\gamma}-1+\frac{\gamma}{\lambda}.
\end{align*}
Finally
\begin{proposition}\label{explicitV}
We have the identities
\[\overline{V}=\log R_0-1+R_0^{-1},\quad e^{N\overline{V}}=R_0^Ne^{-N(R_0-1)/R_0}.
\]
\end{proposition}
Combining this result with Theorem \ref{th:exitT} adapted to the SIS model yields the following.
\begin{corollary}
Suppose that $R_0>1$ and define
\[T^{N,z}_{Ext}=\inf\{t>0,\ Z^{N,z}_t=0\}.\]
Then for any $0<z\le1$, and $c>1$,
\begin{equation*}
\lim_{N\to\infty}\mathbb{P}\left(\left(R_0/c\right)^Ne^{-N(R_0-1)/R_0}<T^{N,z}_{\text{Ext}}<\left(cR_0\right)^Ne^{-N(R_0-1)/R_0}\right)=1,
\end{equation*}
\begin{equation*}
\text{and }\left(R_0/c\right)^Ne^{-N(R_0-1)/R_0} \leq\mathbb{E}(T^{N,z}_{\text{Ext}})\le\left(cR_0\right)^Ne^{-N(R_0-1)/R_0}\end{equation*}
for $N$ large enough.
\end{corollary}

\begin{remark}
In fact, the pair $(\hat{u}_1(t),\hat{u}_2(t))$ is not an optimal control for the above control problem. Such an optimal control does not exist! The optimal trajectory, which is the original ODE time reversed, would take an infinite time to leave $x^\ast$, and an infinite time to reach $0$. However, our $(\hat{u}_1(t),\hat{u}_2(t))$ is the limit of a minimizing sequence obtained by choosing a suboptimal control to drive the system from $x^\ast$ to $x^\ast-\delta$, then the optimal control to drive the system from $x^\ast-\delta$ to $\delta$, and finally a suboptimal control to drive the system from $\delta$ to $0$.  $\log\frac{\lambda}{\gamma}-1+\frac{\gamma}{\lambda}$ is indeed the minimal cost. Note that $\hat{T}=+\infty$.
\end{remark}

\subsection{Time of extinction in the SIR model with demography}

We now turn to the SIR model with demography, which is the model which has been formally presented in Example \ref{ex:SEIRS}, but where we let $\nu=+\infty$ (we suppress the stage E between S and I), and $\gamma=0$ (there is no loss of immunity). The limiting ODE reads
\begin{align*}
\dot{x}_t&=\lambda x_ty_t-\gamma x_t-\mu x_t,\\
\dot{y}_t&= - \lambda x_ty_t +\mu -\mu y_t.
\end{align*}
We assume that $\lambda>\gamma+\mu$, in which case there is a unique stable endemic equilibrium,
namely $z^\ast=(x^\ast,y^\ast)=(\frac{\mu}{\gamma+\mu}-\frac{\mu}{\lambda},\frac{\gamma+\mu}{\lambda})$.
The extinction in such a model has been studied using the Central Limit Theorem for moderate population size in Section \ref{Sec_Open-pop}. We now finally apply Large Deviations to this model. In this model, $Z^N_t=(I^N_t,S^N_t)$ lives in all of $\R^2_+$. We note that in the proof of Proposition \ref{pro:exitT}, the compactness of the set of possible values for $Z^N_t$ has played a crucial role, especially in the proof of Lemma \ref{le:3}. However, if we define for each $R>0$
\[ T^{N,R}_{\text{Ext}} = T^{N}_{\text{Ext}}\wedge\sigma^N_R,\]
where $\sigma_N^R=\inf\{t>0,\, I^N_t+S^N_t\ge R\}$, it is clear that we have reduced our situation to a bounded state space, and the exact same proofs leading to Proposition \ref{pro:exitT} and Theorem \ref{th:exitT}, which easily adapted to this new situation. Moreover, we have the
\begin{lemma}
As $R\to\infty$, $V_R:=\inf_{z=x+y\ge R}V(z^*,z)\to\infty$.
\end{lemma}
\bpf We use the Pontryagin maximum principle and refer to the notations in Section \ref{se:Pont}. Here $d=2$ and $k=5$, $B=\begin{pmatrix}1 & -1 & -1 & 0 & 0\\ -1 & 0 & 0 & 1 & -1 \end{pmatrix}$, $\beta_1(x,y)=\lambda xy$,
$\beta_2(x,y)=\gamma x$, $\beta_3(x,y)=\mu x$, $\beta_4(x,y)=\mu$, $\beta_5(x,y)=\mu y$. The forward-backward ODE system reads
\begin{align*}
\dot{x}_t&=\lambda x_ty_te^{p_t-q_t}-(\gamma+\mu)x_te^{-p_t},\quad x_0=\frac{\mu}{\gamma+\mu}-\frac{\mu}{\lambda}\\
\dot{y}_t&=-\lambda x_ty_te^{p_t-q_t}+\mu e^{q_t}-\mu y_te^{-q_t},\quad y_0=\frac{\gamma+\mu}{\lambda}\\
\dot{p}_t&=\lambda y_t+\gamma+\mu-\lambda y_ye^{p_t-q_t}-\gamma e^{-p_t}-\mu e^{-p_t},\\
\dot{q}_t&=\lambda x_t+\mu-x_te^{p_t-q_t}-\mu e^{-q_t},\quad  p_{\hat{T}}=q_{\hat{T}}.
\end{align*}
Condition \eqref{eq:Ha} at time $\hat{T}$ together with the condition $p_{\hat{T}}=q_{\hat{T}}$ allows us to conclude that
\[ p_{\hat{T}}=q_{\hat{T}}=\log\left(R+\frac{\gamma}{\mu}x\right).\]
It is clear that $\dot{p}_{\hat{T}}>\dot{q}_{\hat{T}}$. In fact it is not hard to show that, as long as $p_t\ge0$,
$p_t<q_t$. However, $\dot{p}_t\le \lambda y_t+\gamma+\mu\le \lambda R+\gamma+\mu$.
Let $a=\frac{1}{2}\frac{\log R}{\lambda R+\gamma+\mu}$. For any $\hat{T}-a\le t\le \hat{T}$,
$p_t\ge\log(R+\frac{\gamma}{\mu}x_{\hat{T}})-\frac{1}{2}\log R\ge \frac{1}{2}\log R>0$. Next we notice that
$\dot{x}_t+\dot{y}_t\le\mu e^{q_t}$. As long as $\hat{T}-a\le t\le \hat{T}$, we both have that $p_t\le q_t$ and $q_t\ge0$, hence $\dot{q}_t\ge0$, and $0<q_t\le q_{\hat{T}}=\log\left(R+\frac{\gamma}{\mu}x_{\hat{T}}\right)$.
Consequently $\dot{x}_t+\dot{y}_t\le\mu \left(R+\frac{\gamma}{\mu}x_{\hat{T}}\right)\le (\gamma+\mu)R$.
Finally, for $\hat{T}-a\le t\le \hat{T}$, $x_t+y_t\ge R-\frac{(\gamma+\mu)R\log R}{2(\lambda R+\gamma+\mu)}\ge
\frac{1}{2}R$, for $R$ large enough.

We can now lowed bound $V_R$.
We use the expression on the left of \eqref{eq:IC} for the instantaneous cost. We have
\begin{align*}
V_R&=\int_0^{\hat{T}}\Big[\lambda x_ty_t(1-e^{p_t-q_t}+(p_t-q_t)e^{p_t-q_t})
+(\mu+\gamma)x_t(1-e^{-p_t}-p_te^{-p_t})\\
&\quad\quad\quad\quad+\mu(1-e^{q_t}+q_te^{q_t})+\mu y_t(1-e^{-q_t}-q_te^{-q_t})\Big] dt \\
&\ge\int_{\hat{T}-a}^{\hat{T}}\mu(x_t+y_t)\inf\{1-e^{-p_t}-p_te^{-p_t},1-e^{-q_t}-q_te^{-q_t}\}dt\\
&\ge \frac{\mu}{8}\frac{R\log R}{\lambda R+\mu+\gamma}\\
&\to+\infty,
\end{align*}
as $R\to\infty$.
\epf

It follows from Corollary \ref{exitLoc} that as soon as $V_R>\hat{V}$, the probability that $Z^N$ exits the truncated domain through the ``extinction boundary'' $\{z^1=0\}$ goes to $1$ as $N\to\infty$. Also, for fixed $N$, $\P(T^{N}_{\text{Ext}}<\sigma^N_R)\to1$, as $R\to\infty$.

\begin{theorem} Let $T^{N,z}_{\text{Ext}}$ denote the extinction time in the $N$--SIR model with demography starting from $z_N=\frac{[zN]}{N}$.
Given $\eta>0$, for all $z\in \R^2_+$ with $z_1>0$,
\begin{equation*}
\lim_{N\to\infty}\mathbb{P}\big(\exp\{N(\overline{V}-\eta)\}<T^{N,z}_{\text{Ext}}<\exp\{N(\overline{V}+\eta)\}\big)=1.
\end{equation*}
Moreover, for all $\eta>0$, $z\in \R^2_+$ with $z_1>0$ and $N$ large enough,
\begin{equation*}
\exp\{N(\overline{V}-\eta)\}\leq\mathbb{E}(T^{N,z}_{\text{Ext}})\leq\exp\{N(\overline{V}+\eta)\}.
\end{equation*}
\end{theorem}

\setcounter{chapter}{1}
\setcounter{section}{0}
\renewcommand{\thechapter}{\Alph{chapter}}
\chapter*{Appendix}
\addcontentsline{toc}{chapter}{Appendix}
\pagestyle{appendixI}

This Appendix presents several mathematical notions, mostly from the theory of stochastic processes, as well as a couple of notions related to continuity of real-valued functions, which are used in the previous chapters.
Most proofs are given. Otherwise we refer to existing monographs.

\section{Branching processes}\label{TB-EP_sec_Br-proc}
We present the basic facts about branching processes, which are useful in these Notes. We give most of the proofs. Those which are missing can be found in classical monographs on branching processes,
see e.g.\  Athreya and Ney \cite{ANI} or Jagers \cite{PJI}, unless we give a precise reference in the text.

\subsection{Discrete time branching processes}\label{Dtbp}
Consider an ancestor (at generation 0) who has $\xi_0$ children, such that
\[
\P(\xi_0=k)=q_k,\ k\ge0\quad\text{and }\sum_{k\ge0}q_k=1.\]
Define $m=\E[\xi_0]=\sum_{k\ge1}k\ q_k$ and $g(s)=\E\left[s^{\xi_0}\right]$.

Each child of the ancestor belongs to generation 1. The $i$-th of those children has himself $\xi_{1,i}$ children,
where the random variables
$\{\xi_{k,i},\ k\ge0, i\ge1\}$ are i.i.d., all having the same law as  $\xi_0$. If we define $X_n$ as the number of individuals in generation $n$, we have
\[ X_{n+1}=\sum_{i=1}^{X_{n}}\xi_{n,i}.\]
We have $g(0)=q_0$, $g(1)=1$, $g'(1)=m$, $g'(s)>0$, $g''(s)>0$, for all $0\le s\le1$ (we assume that $q_0>0$
and $q_0+q_1<1$).
Let us compute the generating function of $X_n$: $g_n(s)=\E[s^{X_n}]$.
\begin{align*}
g_n(s)&=\E\left[s^{\sum_{i=1}^{X_{n-1}}\xi_{n-1,i}}\right]\\
&=\E\left[\E\left[s^{\sum_{i=1}^{X_{n-1}}\xi_{n-1,i}}\Big| X_{n-1}\right]\right]\\
&= \E\left[g(s)^{X_{n-1}}\right]\\
&=g_{n-1}\circ g(s).
\end{align*}
If we iterate this argument, we obtain
\[ g_n(s)=g\circ\cdots\circ g(s),\]
and also
\begin{align*}
\P(X_n=0)&=g^{\circ n}(0)\\
&=g\left[g^{\circ (n-1)}(0)\right].
\end{align*}
Hence if $z_n=\P(X_n=0)$, $z_n=g(z_{n-1})$, and $z_1=q_0$.
\begin{figure}
   \begin{minipage}[c]{.46\linewidth}
   \scalebox{0.45}{
      \includegraphics{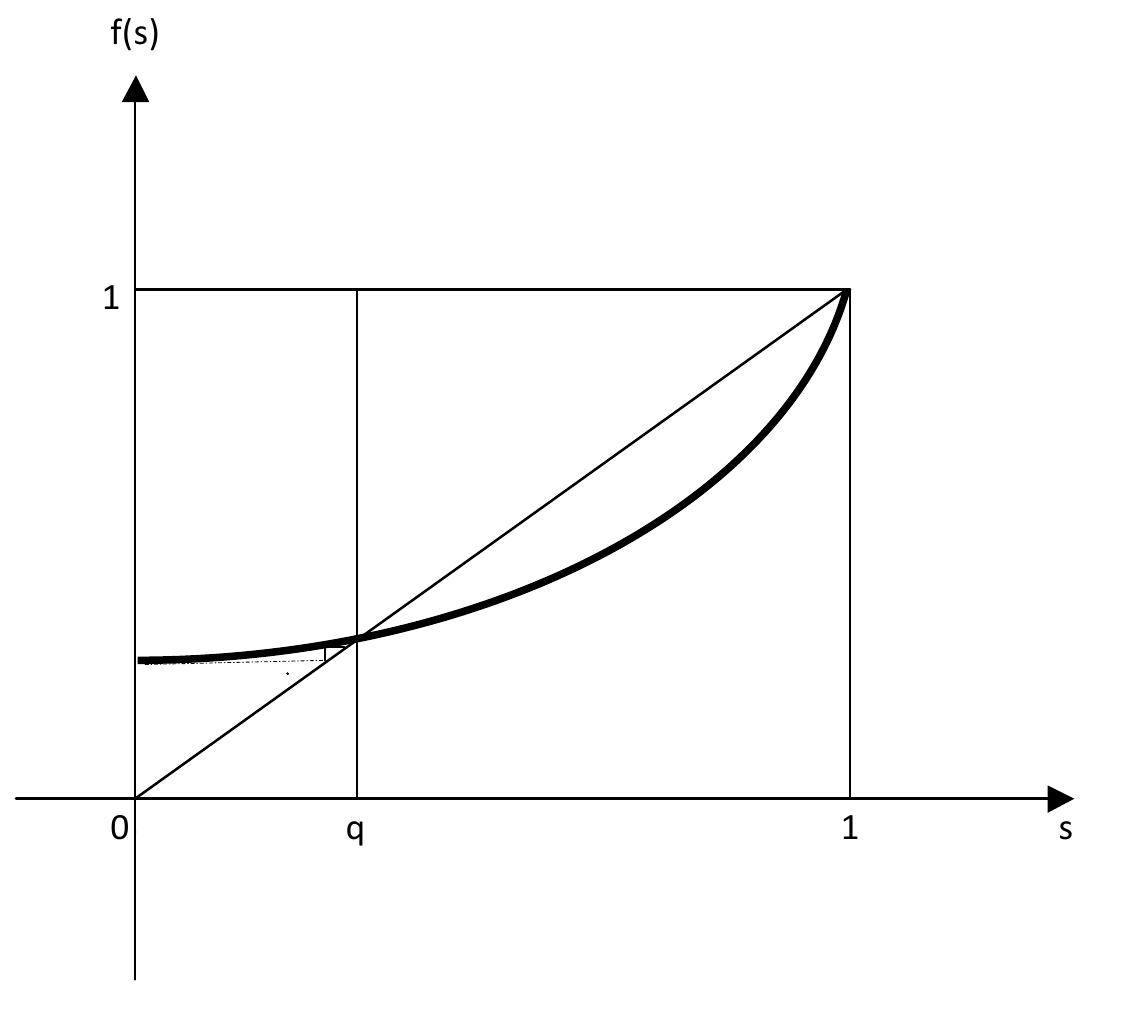}
      }
   \end{minipage} \hfill
   \begin{minipage}[c]{.46\linewidth}
   \scalebox{0.45}{
      \includegraphics{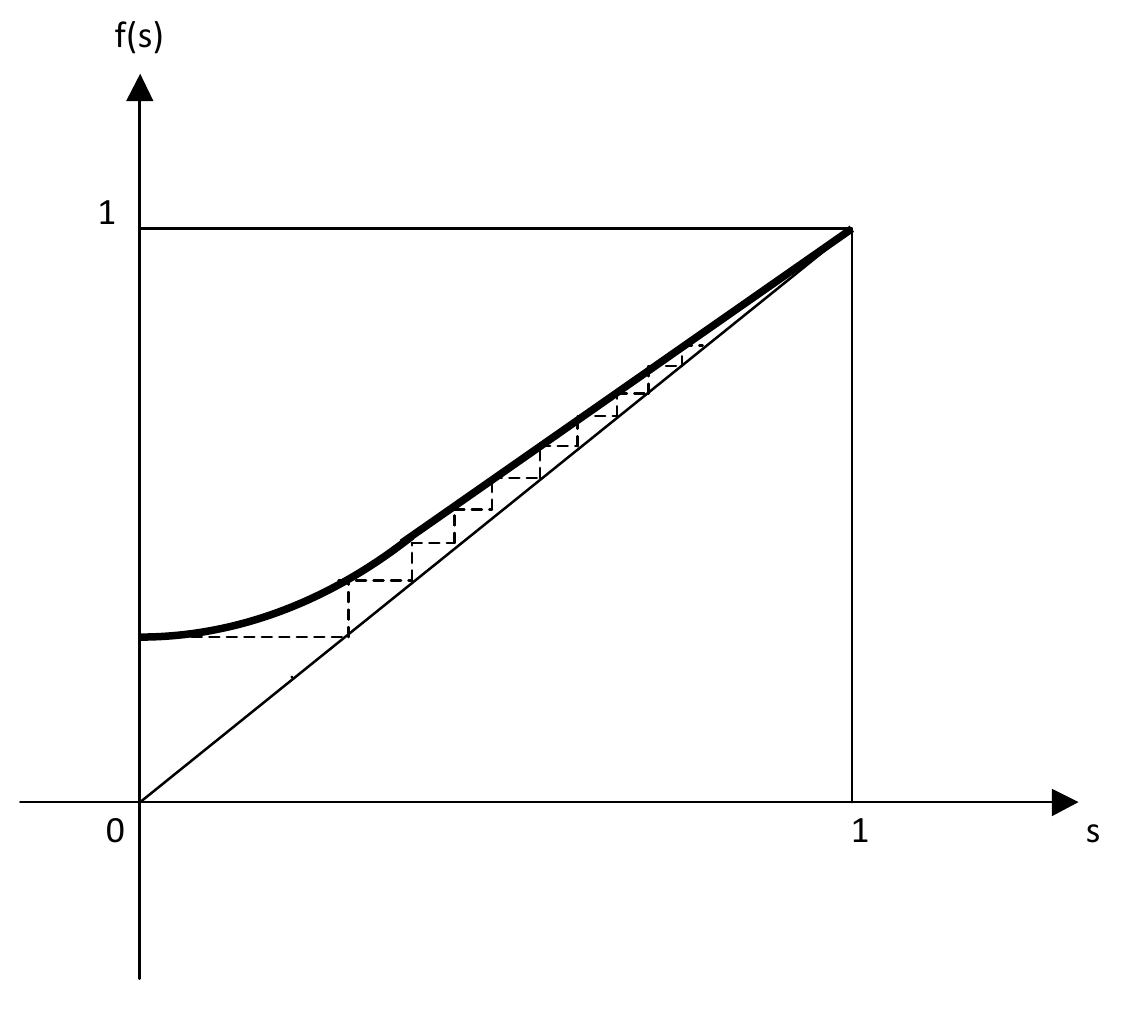}
      }
         \end{minipage}
          \caption {\label{fig0} Graphs of $g$ in case $m>1$ (left) and in case $m\leq1$ (right). The successive heights of the dashed line are the successive values of $\P(X_n=0)$.}
\end{figure}
We have  $z_n\uparrow z_\infty$, where
$z_\infty=\P(X_n=0 \text{ from some $n$})$.
The proof of the following Proposition is essentially
clear from Figure \ref{fig0}.
\begin{proposition}\label{BGW}
If $m\le 1$, then $\P(X_n=0)\to1$ as $n\to\infty$, and $z_\infty=1$.

If $m>1$, $\P(X_n=0)\to z_\infty=q$ as $n\to\infty$, where $q$ is the smallest solution
of the equation $z=g(z)$.
\end{proposition}
Note that on the event $\cup_{n=0}^\infty\{X_n=0\}$, which has probability one in the first case, the population goes extinct after a finite number of generations, and the total progeny is finite.

In the second case, with probability $1-z_\infty$, the branching process does not go extinct.

Let us show that $W_n=m^{-n}\ X_n$ is a martingale.
\begin{align*}
\E(W_{n+1}|X_n)&=m^{-n}\E\left(m^{-1}\sum_1^{X_n}\xi_{n,i}|X_n\right)\\
&=m^{-n}X_n\\
&=W_n.
\end{align*}
One can show that $W_n\to W$ a.s. as $n\to\infty$, and moreover, provided
$\sum_{j\ge1}q_j j\log j<\infty$,
\[ \E[W]=1,\ \text{and }\P(W>0)=\P(\{\text{the branching process does not go extinct}\}). \]
In the case  $\sum_{j\ge1}q_j j\log j=\infty$, then $\P(W=0)=1$.

\subsection{Continuous time branching processes}\label{Sec_Cont_Br_Pr}
We shall consider only binary continuous time branching processes, i.e.\  where at most one child is born at a given time.
This process starts with a single ancestor born at time $t=0$. This ancestor is characterized by a pair $(L_0,\{N_0(t),\; t\ge0\})$,
where  $L_0$ is the life length of the ancestor, and $N_0(t)$ is the number of children of the ancestor born on the time interval $[0,t]$. We assume that $N_0(\infty)=N_0(L_0)$, that is the ancestor does not give birth to offspring after his death. We now assume that the individuals are numbered in the order of their birth. To the individual $i$ is attached
a pair $(L_i,\{N_i(t)\})$, such that the sequence of pairs $\{(L_i,\{N_i(t)\})\}_{i\ge0}$ is i.i.d. If the individual $i$ is born at time
$B_i$, the offspring of individual $i$ are born at the jump times of the process $\{N_i(t-B_i),\ B_i\le t\le B_i+L_i\}$. Note that since $B_i$ depends only upon the pairs $\{(L_j,\{N_j(t)\})\}_{0\le j<i}$, $B_i$ and $(L_i,\{N_i(t)\})$ are independent.

  Let $X_t$ denote the number of individuals in the population alive at time $t$. This process is Markovian if and only if the law of the
  pair $(L_i,N_i(t))$ is such that $L_i$ and $\{N_i(t),\ t\ge0\}$ are independent, $L_i$ is an exponential random variable with parameter $d$, and
  $N_i(t)$ is a rate $b$ Poisson process. We first assume that we are in this situation. We shall denote by $X^k_t$ the number of descendants at time $t$ of $k$ ancestors at time $0$. The branching property implies that $\{X^k_t,\ t\ge0\}$ is the sum of $k$ independent copies of $\{X_t,\ t\ge0\}$. We have the following result.
  \pagebreak

\begin{proposition}\label{contBGW}
The generating function of the process $X$ is given by
$$\E\left(s^{X^k_t}\right)=\psi_t(s)^k,\ s\in[0,1],\ k\ge1,$$
where
$$\frac{\partial\psi_t(s)}{\partial t}=\Phi(\psi_t(s)),\quad \psi_0(s)=s,$$
and the function $\Phi$ is defined by
\begin{align*}
\Phi(s)&=d(1-s)+b(s^2-s)\\
&=(b+d)(h(s)-s),\ s\in[0,1],
\end{align*}
where $h$ is the generating function of the probability measure $\frac{d}{b+d}\delta_0+\frac{b}{d+b}\delta_2$.
\end{proposition}
\bpf   The process $X_t$ is a continuous time $\Z_+$-valued jump Markov process. Denote by $Q$ its infinitesimal generator. The non-zero elements of the $n$-th row of $Q$ are given by
$$Q_{n,m}=\begin{cases} nd,&\text{if $m= n-1$},\\
                   -n(b+d),&\text{if $m=n$};\\
                   nb,&\text{if $m=n+1$}. \end{cases}$$
Define $f:\N\to[0,1]$ by $f(k)=s^k$, $s\in[0,1]$. Then $\psi_t(s)=P_tf(1):=\E[f(X^1_t)]$ (we use the unusual notation $X^1=X$ to stress the fact that the process starts from $X_0=1$).   It follows from the backward Kolmogorov equation for the process $X$ (see e.g.\ Theorem 3.2, Chapter 7 in Pardoux \cite{EPI}) that
\begin{align*}
\frac{dP_tf(1)}{dt}&=(QP_tf)(1)\\
\frac{\partial \psi_t(s)}{\partial t}&= Q_{1,0}+Q_{1,1}\psi_t(s)+Q_{1,2}\psi_t(s)^2\\
  &= d-(b+d)\psi_t(s)+b\psi_t(s)^2\\
  &=\Phi(\psi_t(s)).
  \end{align*}
  \epf

  \begin{corollary}
  We have
  \[ \E[X^k_t]=k e^{r t},\quad\text{where }r=b-d.\]
  \end{corollary}
  \bpf Differentiating with respect to $s$ the above equation for $\psi_t(s)$ yields
  \begin{align*}
  \frac{\partial}{\partial t}\left(\frac{\partial}{\partial s}\psi_t(s)\right)&=\Phi'(\psi_t(s))\frac{\partial}{\partial s}\psi_t(s)\\
  &=(b+d)(h'(s)-1)\frac{\partial}{\partial s}\psi_t(s).
  \end{align*}
   The last equation at $s=1$ yields
  \begin{align*}
  \frac{d}{dt}\E[X_t]=r\E[X_t],
  \end{align*}
  where $X_t=X^1_t$. The result follows for $k=1$, and then the general case, since the mean number of offspring of $k$ ancestors equals
  $k$ times the mean number of offspring of one ancestor.
  \epf

 The quantity $r$ is often referred to as the Malthusian parameter. It is the mean number of births minus the mean number of death per unit time. Another important quantity is the mean number of offspring of each individual, which is equal to $m=b/d$. The process $X^k_t$ is said to be subcritical if $m<1$, i.e.\  $r<0$. In that case $X^k_t\to0$ in $L^1(\Omega)$, and it is easy to show that $X^k_t=0$ for $t$ large enough. This last conclusion holds in the critical case ($m=1$, i.e.\  $r=0$) as well. In those two cases, the total progeny is finite a.s. We now study the large time behaviour of $X^k_t$ in the supercritical case.
  In the next proposition, we again write $X_t$ for $X^1_t$.
  \begin{proposition}\label{contBGW_large_t}
  If $m>1$, or equivalently $r>0$, there exists a non-negative random variable $W$ such that $X_t\sim We^{r t}$ almost surely, as $t\to\infty$.
Moreover $\{W=0\}=\{\exists t>0 \text{ such that }X_t=0\}$ and
\[\P(W=0)=\P(\{\exists t>0 \text{ such that }X_t=0\})=\frac{d}{b}.\]
\end{proposition}
\bpf The first part of the result follows readily from the fact that $e^{-r t}X_t$ is a positive martingale, which converges a.s. to a limit $W$ as $t\to\infty$. Moreover it is not hard to show that $\sup_{t>0}\E[e^{-2rt}(X_t)^2]<\infty$, hence the convergence holds in $L^1(\Omega)$, so $\E[W]=1$.
Now clearly $\{\exists t>0 \text{ s.t. }X_t=0\}\subset\{W=0\}$. If we start with $k$ ancestors, the limiting $W$ is clearly the sum of
$k$ i.i.d.\ copies of $W$ when starting with one ancestor, and $\P_k(W=0)=(\P_1(W=0))^k$. It is now easy to deduce that
\[ \P_1(W=0| X_t)=\P_1(W=0))^{X_t}.\]
Taking the expectation in this identity and writing $q=\P_1(W=0)$, we obtain $q=\E[q^{X_t}]$. Differentiating that identity at $t=0$
and taking advantage of Proposition \ref{contBGW}, we deduce that $q$ solves $bq^2-(b+d)q+d=0$. Moreover since $\E(W)=1$,
$q<1$, hence $q=d/b$.
Finally
$\P(X_t=0)=\psi_t(0)$ is the solution of the ODE $\dot{x}(t)=bx(t)^2-(b+d)x(t)+d$, $x(0)=0$. It is clear that as $t\to\infty$,
$\psi_t(0)$ increases to the smallest solution of the equation $bs^2-(b+d)s+d=0$, again $d/b$.
\epf

We now consider non-Markovian continuous time binary branching processes. The non-Markovian continuous time branching processes which we have described at the beginning of this section are called Crump--Mode--Jagers processes.
Now the law of the pairs $(L_i,\{N_i(t)\})$ can be quite general. For the application to epidemics models, we can consider the case
where $L_i$ and $\{N_i(t)\}$ are independent, $N_i$ being a Poisson process, but the law of $L_i$ is no longer exponential.
We denote again by $m=\E[N_0(L_0)]$ the mean number of offspring of each individual. Of course, the process is subcritical, critical, or supercritical according as $m<1$, $m=1$ or $m>1$. We denote again by $X_t$ the number of individuals alive at time $t$. We define $F(t)=\E[N(t)]$ and $G(t)=\P(L\le t)$. We assume that $F$ is non-lattice, and $F(0^+)<1$. Doney \cite{RADI} showed the following two results.
\begin{proposition}\label{pro:malthus}
If $1<m<\infty$, then there exists a unique $r>0$ such that
\[\int_0^\infty e^{-r t} F(dt)=1\]
and $\E[X_t]\sim a e^{r t}$, where
\begin{equation*}
0 < a=\frac{\int_0^\infty(1-G(t))e^{-r t}dt}{\int_0^\infty te^{-r t}F(dt)}<\infty.
\end{equation*}
\end{proposition}
 Again $r$ is called the Malthusian parameter. In the next statement, we use the notation
\[ Y=\int_0^\infty e^{-rt}N(dt).\] It is clear that
\[ \E[Y]=\int_0^\infty e^{-rt}F(dt)=1.\]
\begin{theorem}\label{th:supercrit}
Suppose that $1<m<\infty$. Then, as $t\to\infty$
\[ \frac{X_t}{\E[X_t]}\to W\quad\text{in law}.\]
$W$ is not identically $0$ if and only if $\E[Y\log(Y)]<\infty$, in which case $\E[W]=1$ and $\P(W=0)=\P(\{\exists t>0 \text{ s.t. }X_t=0\})$.
Moreover, the law of $W$ has an atom at $0$ and is absolutely continuous on $(0,\infty)$.
\end{theorem}

\section[The Poisson process and Poisson point process]{The Poisson process and Poisson point process}\label{sec-Po-proc}

The Poisson process is central in this whole volume.
Let $\lambda>0$ be given. A rate $\lambda$ Poisson (counting) process is defined as
\[ P_t=\sup\{k\ge1,\ T_k\le t\},\]
where $0=T_0<T_1<T_2<\cdots<T_k<\cdots<\infty$, the random variables $\{T_k-T_{k-1},\ k\ge1\}$ being independent and identically distributed, each following the law $\mathrm{Exp}(\lambda)$.
We have
\begin{proposition}\label{prop:Poisson}
For all $n\ge1$, $0<t_1<t_2<\cdots<t_n$, the random variables $P_{t_1}, P_{t_2}-P_{t_1},\ldots,P_{t_n}-P_{t_{n-1}}$
are independent, and
for all $1\le k\le n$, $P_{t_k}-P_{t_{k-1}}\sim\text{Poi}[\lambda(t_k-t_{k-1})]$.
\end{proposition}
\bpf Let us first prove that for all $t,s>0$,
\[\P(P_{t+s}-P_t=0|P_t=k, T_1,T_2,\ldots,T_k)=\exp(-\lambda s).\]
Indeed
\begin{align*}
\P(P_{t+s}-P_t=0&|P_t=k, T_1,T_2,\ldots,T_k)\\&=\P(T_{k+1}>t+s|P_t=k,T_k)\\
&=\P(T_{k+1}-T_k>t+s-T_k|T_{k+1}-T_k>t-T_k>0)\\
&=\P(T_{k+1}-T_k>s)\\
&=\exp(-\lambda s).
\end{align*}
Let now $n\ge1$. For $1\le i\le n$, we define $X_{n,i}={\bf1}_{\{P_{t+is/n}-P_{t+(i-1)s/n}\ge1\}}$, and finally
$S_n=X_{n,1}+X_{n,2}+\cdots+X_{n,n}$. It follows from the first part of the proof that conditionally upon
$\sigma\{P_r,\ 0\le r\le t\}$, the random variables $X_{n,1}, X_{n,2}, \ldots, X_{n,n}$ are i.i.d., each Bernoulli with parameter
$1-e^{-\lambda s/n}$. Then conditionally upon $\sigma\{P_r,\ 0\le r\le t\}$, $S_n$ is binomial with parameters
$(n, 1-e^{-\lambda s/n})$. But $S_n\to P_{t+s}-P_t$ a.s. as $n\to\infty$, while its conditional law given
$\sigma\{P_r,\ 0\le r\le t\}$ converges towards the Poisson distribution with parameter $\lambda s$,
according to the following lemma. The proposition follows.  \epf

We have used the following well-known result. Recall the notation $\text{Bin}(n,p)$ for the binomial law with parameters $n$ and $p$, where $n\ge1$ and $0<p<1$.
\begin{lemma}
For all $n\ge1$, let $U_n$ be a $\text{Bin}(n,p_n)$ random variable. If $np_n\to\lambda$ as $n\to\infty$, with $\lambda>0$,
then $U_n$ converges in law towards Poi$(\lambda)$.
\end{lemma}

A Poisson process will be called standard if its rate is 1. If $P$ is a standard Poisson process, then
$\{P(\lambda t), \  t\ge0\}$ is a rate $\lambda$ Poisson process.

We will also use the following
\begin{exercise}\label{ex:select}
Let $\{P_t, \  t\ge0\}$ be a rate $\lambda$ Poisson process, and $\{T_k,\ k\ge1\}$ the random points of this Poisson process, i.e.\  for all $t>0$, $P_t=\sup\{k\ge1,\ T_k\le t\}$. Let $0<p<1$. Suppose that each $T_k$ is selected with probability $p$, not selected with probability $1-p$, independently from the others.  Let $P'_t$ denote the number of selected points on the interval $[0,t]$. Then $\{P'_t, \  t\ge0\}$ is a rate $\lambda p$ Poisson process.
\end{exercise}

A rate $\lambda$  Poisson process ($\lambda>0$) is a counting process $\{R_t,\ t\ge0\}$ such that
  $R_t-\lambda t$ is a martingale. Let $\{ P(t),\ t\ge0\}$ be a standard Poisson process (i.e.\  with rate $1$).
 Then $P(\lambda t)-\lambda t$ is martingale, and it is not hard to show that $\{P(\lambda t),\ t\ge0\}$ is a rate
 $\lambda$ Poisson process. Let now $\{\lambda(t),\ t\ge0\}$ be a measurable and locally integrable $\R_+$-valued function.
 Then the process $\{R_t:=P\left(\int_0^t\lambda(s)ds\right),\ t\ge0\}$ is called a rate $\lambda(t)$ Poisson
 process.
Clearly $R_t-\int_0^t\lambda(s)ds$ is a martingale.

We now want to consider the case where $\lambda$ is random.  For that purpose, it is convenient to
give an alternative definition of  the above process $R_t$.

Consider a standard Poisson random measure $Q$ on $\R_2^+$, which is defined as follows.
$M$ is the counting process associated to a random cloud of points in $\R_+^2$. One way to construct that cloud of points is as follows. We can consider $\R_+^2=\cup_{i=1}^\infty A_i$, where the $A_i$'s are disjoint squares with Lebesgue measure $1$. Let $K_i$, $i\ge1$ be i.i.d.\ mean one Poisson random variables. Let
$\{X^i_j,\ j\ge1, i\ge1\}$ be independent random points of $\R_+^2$, which are such that for any $i\ge1$, the
$X^i_j$'s are uniformly distributed in $A_i$. Then
\[  Q(dx)=\sum_{i=1}^\infty\sum_{j=1}^{K_i}\delta_{X^i_j}(dx).\]
$\lambda(t)$ denoting a positive-valued measurable function, the above $\{R_t,\ t\ge0\}$ has the same law as
\[R_t=\int_0^t\int_0^{\lambda(s)}Q(ds,du).\]

 Now let $\{\lambda(t),\ t\ge0\}$ be an $\R_+$-valued stochastic process, which
is assumed to be predictable, in the following sense. Let for $t\ge0$,
\[\F_t=\sigma\{Q(A),\ A \text{ Borel subset of }[0,t]\times\R_+\},\]
and consider the $\sigma$-algebra of subset
of $[0,\infty)\times\Omega$ generated by the subsets of the form ${\bf1}_{(s,t]}{\bf1}_F$, where $0\le s<t$
and $F\in\F_s$, which is called the predictable $\sigma$-algebra. Note that if $X_t$ is $\F_t$-progressively measurable and left-continuous, then it is predictable. If $X_t$ is progressively measurable and  right-continuous, then $X_{t-}$ is predictable.

We assume moreover that $\E\int_0^t\lambda(s)ds<\infty$ for all $t>0$. We now define the process $R_t$ as above:
\[ R_t=\int_0^t\int_0^{\lambda(s)}Q(ds,du).\]
We have (see the next subsection for the definition of a martingale)
\begin{lemma}\label{le:PoisIntMart}
$R_t-\int_0^t\lambda(s)ds$ is a martingale.
\end{lemma}
\bpf
  For any $\delta>0$, let
  \[R^\delta_t=\int_0^t\int_0^{\lambda(s-\delta)}Q(ds,du),\]
 where $\lambda(s)=0$ for $s<0$. It is not hard to show  that
 $R^\delta_t-\int_0^t\lambda(s-\delta)ds$ is a martingale which converges in $L^1(\Omega)$
 to $R_t-\int_0^t\lambda(s)ds$. Indeed, it suffices to show that if $0<s<t$ with $t-s\le\delta$,
 the restriction of the random measure $M$ to $(s,t]\times(0,+\infty)$ is independent of
 $\{\lambda(r-\delta),\, s<r\le t\}$, which is $\F_s$ measurable hence
 \[ \E^{\F_s}(R^\delta_t-R^\delta_s)=\E^{\F_s}\int_s^t\lambda(r-\delta)dr.\]
 The result follows.
 \epf

 The process $R_t$ is sometimes called ``a doubly stochastic Poisson process'' or a Cox process. Of course the increments of
 $R_t$ are not Poisson distributed.
 If we let $\sigma(t)=\inf\{r>0,\ \int_0^r\lambda(s)ds>t\}$, we have that $P(t):=R_{\sigma(t)}$ is
 a standard Poisson process, and it is clear that
 $R_t=P\left(\int_0^t\lambda(s)ds\right)$.

  In particular, the process which counts the new infections, which appears
   in Section \ref{TB-EP_sec_LLN}, takes the form
 \[
 P\left(\frac{\lambda}{N}\int_0^t I(r)S(r)dr\right)=\int_0^t\int_0^{\infty}{\bf1}_{u\le\frac{\lambda}{N}I(r-)S(r-)}Q(ds,du).\]
 If we let $\overline{Q}(ds,du)=Q(ds,du)-ds\times du$ and $M(t):=P(t)-t$, it is clear that, as a consequence of the above Lemma, we have
 \begin{corollary}\label{co:PoisIntMart}
 Define $M(\cdot)$ by
 \begin{align*}
 M\left(\frac{\lambda}{N}\int_0^t I(r)S(r)dr\right)&=\int_0^t\int_0^{\infty}{\bf1}_{u\le\frac{\lambda}{N}I(r-)S(r-)}\overline{Q}(ds,du)\\
 &=\int_0^t\int_0^{\infty}{\bf1}_{u\le\frac{\lambda}{N}I(r-)S(r-)}Q(ds,du)-\frac{\lambda}{N}\int_0^t I(r)S(r)dr.
 \end{align*}
 Then $M(t)$ is a martingale (see Definition \ref{de:martcont} below).
 \end{corollary}
 Note that $\int_0^t I(r-)S(r-)dr=\int_0^t I(r)S(r)dr$ since the two integrands coincide
 $dr$ a.e. since they differ on each interval $[0,t]$ at most at finitely many points. We use
 the second formulation, since it is simpler.

\section{Cram\'er's theorem for Poisson random variables}\label{CramerPoisson}
In order to explain what Large Deviations is about, let us first establish Cram\'er's Theorem,
in the particular case of Poisson random variables.
Let $X_1,X_2,\ldots,X_n,\ldots$ be mutually independent Poi$(\mu)$ random variables. The Law of Large Numbers tells us that
\[\frac{1}{N}\sum_{i=1}^NX_i\to\mu\quad\text{  a.s. as }N\to\infty.\]
Let us first define, for $X\sim\text{Poi}(\mu)$ the logarithm of its Laplace transform
\[\Lambda(\lambda)=\log\E[\exp(\lambda X)]=\mu(e^\lambda-1),\]
and the Fenchel--Legendre transform of the latter
\[\Lambda^\ast(x)=\sup_{\lambda\in\R}\left\{\lambda x-\Lambda(\lambda)\right\}=
x\log\left(\frac{x}{\mu}\right)-x+\mu.\]
Note that the minimum of $\Lambda^\ast$ is achieved at $x=\mu$, and $\Lambda^\ast$
is zero at that point.

Let $\nu_N$ denote the law of the random variable $\frac{1}{N}\sum_{i=1}^NX_i$.
We can now state Cram\'er's theorem.
\begin{theorem}\label{th:cramer}
Let $F\subset\R$ be a closed set.
\begin{align*}
\text{For any }N\ge1,\ \nu_N(F)&\le \exp\left(-N\inf_{x\in F}\Lambda^\ast(x)\right).\\
\text{Hence }\ \limsup_{N\to\infty}\frac{1}{N}\log\nu_N(F)&\le-\inf_{x\in F}\Lambda^\ast(x).
\end{align*}
Let $G\subset\R$ be an open set.
\begin{align*}
\text{For any }N\ge1,\ \nu_N(G)&\ge \exp\left(-N\inf_{x\in G}\Lambda^\ast(x)\right).\\
\text{Hence }\ \liminf_{N\to\infty}\frac{1}{N}\log\nu_N(G)&\ge-\inf_{x\in G}\Lambda^\ast(x).
\end{align*}
\end{theorem}
\bpf
\noindent{\sc First step. Proof of the upper bound}
Let $X_1,X_2,\ldots,X_n,\ldots$ be
mutually independent Poi$(\mu)$ random variables. For $\sigma>\mu$, we want to estimate 
\[\P\left(\frac{1}{N}\sum_{i=1}^NX_i\ge\sigma\right),\]
which is the probability of a Large Deviation from the LLN, since we know that for large $N$, $\frac{1}{N}\sum_{i=1}^NX_i\simeq\mu$.

For any $\lambda>0$, using Chebycheff's inequality,
\begin{align*}
\P\left(\frac{1}{N}\sum_{i=1}^NX_i\ge\sigma\right)
&=\P\left(\exp\left\{\lambda\left(\sum_{i=1}^NX_i-N\sigma\right)\right\}\ge1\right)\\
&\le\E\exp\left\{\lambda\left(\sum_{i=1}^NX_i-N\sigma\right)\right\}\\
&=\exp\left[-N(\lambda \sigma-\Lambda(\lambda)\right].
\end{align*}
The best possible upper bound is then (since with $\sigma>\mu$, $\Lambda^\ast(\sigma)$ is obtained by taking the supremum over $\lambda>0$)
\begin{align*}
\P\left(\frac{1}{N}\sum_{i=1}^NX_i\ge\sigma\right)&\le e^{-N\Lambda^\ast(\sigma)}\\
&=\exp\left[-N\left(\sigma\log\left(\frac{\sigma}{\mu}\right)-\sigma+\mu\right)\right].
\end{align*}

Similarly, if $\sigma<\mu$, for any $\lambda<0$,
\begin{align*}
\P\left(\frac{1}{N}\sum_{i=1}^NX_i\le\sigma\right)&\le\E\exp\left\{\lambda\left(\sum_{i=1}^NX_i-N\sigma\right)\right\}\\
&=\exp\left[-N(\lambda \sigma-\Lambda(\lambda)\right].
\end{align*}
Since with $\sigma<\mu$, $\Lambda^\ast(\sigma)$ is obtained by taking the supremum over $\lambda<0$,
the above computation leads again to
\[\P\left(\frac{1}{N}\sum_{i=1}^NX_i\le\sigma\right)\le
\exp\left[-N\left(\sigma\log\left(\frac{\sigma}{\mu}\right)-\sigma+\mu\right)\right].\]
It is not hard to see that the upper bound follows from the two above estimates.
\smallskip

\noindent{\sc Second step. Proof of the lower bound}
For any $\delta>0$,
\begin{align*}
\nu_N((-\delta,\delta))&\ge\nu_N(\{0\})=e^{-N\mu}, \quad\text{hence }
\frac{1}{N}\log\nu_N((-\delta,\delta))\ge-\mu=-\Lambda^\ast(0).
\end{align*}
Since transforming $X$ into $Y=X-x$ results in $\Lambda$ and $\Lambda^\ast$ being transformed into
$\Lambda_Y(\lambda)=\Lambda(\lambda)-\lambda x$ and $\Lambda_Y^\ast(\cdot)=\Lambda^\ast(\cdot+x)$, the above yields that for all $x>0$,
\[\frac{1}{N}\log\nu_N((x-\delta,x+\delta))\ge-\Lambda^\ast(x).\]
The lower bound follows readily.  \epf

\section{Martingales}\label{TB-EP_sec_martingales}

\subsection{Martingales in discrete time}
$(\Omega,\F,\P)$ being our standing probability space, let be given an increasing
sequence $\{\F_n,\ n\ge0\}$ of sub-$\sigma$-algebras of $\F$.
\begin{definition}
A sequence $\{X_n,\ n\ge0\}$ of random variables is called a martingale if
\begin{enumerate}
\item For all $n\ge0$, $X_n$ is $\F_n$-measurable and integrable,
\item For all $n\ge0$,  $\E(X_{n+1}|\F_n)=X_n$ a. s.
\end{enumerate}
A sub-martingale is a sequence which satisfies the first condition and \\ $\E(X_{n+1}|\F_n)\ge X_n$.
A super-martingale is a sequence which satisfies the first condition and
$\E(X_{n+1}|\F_n)\le X_n$.
\end{definition}

It follows readily from Jensen's inequality for conditional expectations the
\begin{proposition}\label{Jensen}
If $\{X_n,\ n\ge0\}$ is a martingale, and $\varphi:\R\to\R$ is a convex function such that
$\varphi(X_n)$ is integrable for all $n\ge0$, then $\{\varphi(X_n),\ n\ge0\}$ is a
 sub-martingale.
 \end{proposition}
 We shall need the notion of stopping time
 \begin{definition} A stopping time $\tau$ is an $\Z_+\cup\{+\infty\}$-valued random variable
which satisfies $\{\tau=n\}\in\F_n$, for all $n\ge0$.
\end{definition}
We also let
\[ \F_\tau=\{B\in\F,\ B\cap\{\tau=n\}\in\F_n,\ \forall n\in\Z_+\}.\]
We have Doob's optional sampling theorem:
 \begin{theorem} If $\{X_n,\ n\ge0\}$ is a martingale (resp.\ a sub-martingale), and
 $\tau_1$, $\tau_2$ two stopping times s.t. $\tau_1\le\tau_2\le N$ a.s., then
 $X_{\tau_i}$ is $\F_{\tau_i}$ measurable and integrable, $i=1,2$, and moreover
 \begin{align*}
 \E(X_{\tau_2}|\F_{\tau_1})&=X_{\tau_1}\\
 (\text{resp.\ }\
  \E(X_{\tau_2}|\F_{\tau_1})&\ge X_{\tau_1}).
  \end{align*}
  \end{theorem}
  \bpf  For all $A\in\BB$, $n\ge0$,
  $$\{X_{\tau_i}\in A\}\cap\{\tau_i=n\}=\{X_n\in A\}\cap\{\tau_i=n\}\in\F_n,$$
  and moreover
  $$|X_{\tau_i}|\le\sum_{k=1}^N |X_k|,$$ which establishes the first part of the statement.

 Let $A\in\F_{\tau_1}$. Then
 $$A\cap\{\tau_1<k\le\tau_2\}=A\cap\{\tau_1\le k-1\}\cap\{\tau_2\le k-1\}^c\in\F_{k-1}.$$
 Indeed, we have
 \begin{align*}
 A\cap\{\tau_1\le k-1\}=\cup_{j=1}^{k-1} A\cap\{\tau_1=j\}\ \in\F_{k-1}, \
 \text{and } \{\tau_2\le k-1\}^c\in\F_{k-1}.
 \end{align*}
 Let $\Delta_k=X_k-X_{k-1}$. We have, with $A\in\F_{\tau_1}$,
 \begin{align*}
 \int_A(X_{\tau_2}-X_{\tau_1})d\P&=\int_A\sum_{k=1}^n{\bf1}_{\{\tau_1<k\le\tau_2\}}\Delta_kd\P\\
 &=\sum_{k=1}^n\int_{A\cap\{\tau_1<k\le\tau_2\}}\Delta_kd\P\\
 &=0
 \end{align*}
 or else $\ge0$, depending upon whether $\{X_n,\ n\ge0\}$ is a martingale or a sub-martingale.
 \epf

 We have a first Doob's inequality
  \begin{proposition}\label{Doob1}
 If $X_1,\ldots,X_n$ is a sub-martingale, then for all $\alpha>0$,
 $$\P\left(\max_{1\le i\le n}X_i\ge\alpha\right)\le\frac{1}{\alpha}\E(X_n^+).$$
 \end{proposition}
 \bpf  Define the stopping time $\tau=\inf\{0\le k\le n,\ X_k\ge\alpha\}$ and let
 $M_k=\max_{1\le i\le k}X_i$. We have
 $$\{M_n\ge\alpha\}\cap\{\tau\le k\}=\{M_k\ge\alpha\}\in\F_k.$$
Hence $\{M_n\ge\alpha\}\in\F_{\tau}$. From the optional sampling Theorem,
\begin{align*}
\alpha\P(M_n\ge\alpha)&\le\int_{\{M_n\ge\alpha\}}X_\tau d\P\\
&\le \int_{\{M_n\ge\alpha\}}X_nd\P\\
&\le \int_{\{M_n\ge\alpha\}}X_n^+d\P\\
&\le\E(X_n^+).
\end{align*}
\epf

We have finally a second Doob's inequality
\begin{proposition}\label{Doob2}
If $M_1,\ldots,M_n$ is a martingale, then
\[\E\left[\sup_{0\le k\le n}|M_k|^2\right]\le 4 \E\left[|M_n|^2\right].\]
\end{proposition}
\bpf
Let $X_k=|M_k|$. From Proposition \ref{Jensen}, $X_1,\ldots,X_n$ is a sub-martingale.
It follows from the proof of Proposition \ref{Doob1} that, with the notation
$X^\ast_k=\sup_{0\le k\le n} X_k$,
\[  \P(X^\ast_n>\lambda)\le \frac{1}{\lambda}\E\left(X_n{\bf1}_{X^\ast_n>\lambda}\right).\]
Consequently
\begin{align*}
\int_0^\infty \lambda\P(X^\ast_n>\lambda)d\lambda&\le\int_0^\infty
\E\left(X_n{\bf1}_{X^\ast_n>\lambda}\right) d\lambda\\
\E\left(\int_0^{X^\ast_n}\lambda d\lambda\right)&\le \E\left(X_n\int_0^{X^\ast_n}d\lambda\right)\\
\frac{1}{2}\E\left[|X^\ast_n|^2\right]&\le \E( X_n X^\ast_n)\\
&\le \sqrt{E(|X_n|^2)}\sqrt{E(|X_n^\ast|^2)},
\end{align*}
from which the result follows. \epf

\subsection{Martingales in continuous time}\label{sec:contmart}
We are now given an increasing collection $\{\F_t,\ t\ge0\}$ of sub-$\sigma$-algebras in continuous time.
\begin{definition}\label{de:martcont}
A process $\{X_t,\ t\ge0\}$  is called a martingale if
\begin{enumerate}
\item for all $t\ge0$, $X_t$ is $\F_t$-measurable and integrable;
\item for all $0\le s<t$,  $\E(X_t|\F_s)=X_s$ a. s.
\end{enumerate}
A sub-martingale is a process which satisfies the first condition and $\E(X_t|\F_s)$ $\ge X_s$.
A super-martingale is a process which satisfies the first condition and
$\E(X_t|\F_s)\le X_s$.
\end{definition}

Suppose $\{M_t,\ t\ge0\}$ is a right-continuous martingale. For any $n\ge1$,
$0=t_0<t_1<\cdots<t_n$, $(M_{t_0}, M_{t_1},\ldots,M_{t_n})$ is a discrete time martingale, to which
Proposition \ref{Doob2} applies. Since
\[ \sup_{0\le s\le t}|M_s|=\sup_{\text{Partitions of }[0,t]}\sup_{1\le k\le n}|M_{t_k}|,\]
 Proposition \ref{Doob2} implies readily
\begin{proposition}\label{Doob2c}
 If $\{M_t,\ t\ge0\}$ is a right-continuous martingale,
  \[\E\left[\sup_{0\le s\le t}|M_s|^2\right]\le 4 \E\left[|M_t|^2\right].\]
\end{proposition}

We now establish a particular (essentially obvious) instance of It\^o's formula.
Recall that an $\R$-valued function of $t$ has locally bounded variations if and only if it is the difference of an increasing and a decreasing function. This class of functions excludes all non-zero continuous martingales,
e.g.\ Brownian motion. But all processes considered in these Notes, except for the limit in the functional central limit theorem, are locally of bounded variations.
Given such a locally bounded variation
right-conti\-nu\-ous $1$-dimensional process $X_t$, we define the  bracket
$[X,X]_t=\sum_{0\le s\le t}|\Delta X_s|^2$, where $\Delta X_s=X_s-X_{s-}$ is the jump of $X$ at time $s$. It follows from the fact that $X$ has bounded variation on any compact interval that the set
$\{s\ge0,\, \Delta X_s\not=0\}$ is at most countable, hence the above sum makes sense. If $X$ and $Y$ are two processes of the above type, then
\[[X,Y]_t=\sum_{0\le s\le t}\Delta X_s\Delta Y_s=\frac{1}{2}\left([X+Y,X+Y]_t-[X,X]_t-[Y,Y]_t\right).\]
Now we have what we call It\^o's formula. If $X_t$ and $Y_t$ are right-continuous and have left limits at any $t$, have bounded variations on any compact interval, then for any $t>0$,
\begin{equation}\label{ito}
X_tY_t=X_0Y_0+\int_0^tX_{s-}dY_s+\int_0^tY_{s-}dX_s+[X,Y]_t.
\end{equation}
In case all jumps of $X$ and $Y$ are isolated, which is the only situation treated in these Notes, the result follows clearly by analyzing the evolution of both sides of the identity between the jumps, and at the jump times. The result in the more general situation is easily deduced by approximation.

If $M_t$ is  a right-continuous $\R$-valued martingale with locally bounded variation, we define as above its quadratic variation as
\[ [M,M]_t=\sum_{0\le s\le t}|\Delta M_s|^2,\]
and $\langle M,M\rangle_t$ as the unique increasing predictable process such that $[M,M]_t-\langle M,M\rangle_t$ is a martingale.
Note that both $M^2_t-[M,M]_t$ and $M^2_t-\langle M,M\rangle_t$ are martingales. Consequently, we have in particular
\begin{proposition}\label{SecondMomentMart}
Let $M_t$ be  a square--integrable right-continuous $\R$-valued martingale with finite variation such that $M_0=0$.
Then for all $t>0$,
\[ \E\left(M_t^2\right)=\E\left(\sum_{0\le s\le t}|\Delta M_s|^2\right).\]
\end{proposition}

\section{Tightness and weak convergence in path space}\label{tight}
In these Notes we consider continuous time processes with values in $\R^d$. Most of our processes
are discontinuous. Their trajectories belong to the set $D([0,+\infty);\R^d)$ of functions which are
right continuous and have left limits at any point $t\in[0,+\infty)$. It is not very convenient to use
the topology of locally uniform convergence on this set, since we would like for instance the two Heaviside type functions ${\bf1}_{[1,+\infty)}(t)$ and ${\bf1}_{[1+\eps,+\infty)}(t)$ to be close for
$\eps$ small. The Skorokhod topology essentially says that two functions are close if after a time change which is close to the identity, they are (at least locally) close in the supremum topology.
The only weak convergence (i.e.\  convergence in law) results we consider in these Notes are convergence results towards a continuous process. In this case, convergence in the sense of the Skorokhod topology is equivalent to locally uniform convergence.

Note also that weak convergence of a sequence of processes $X^n$ towards $X$ is equivalent to the two following facts:
\begin{enumerate}
\item The sequence $\{X^n\}_{n\ge1}$ is tight, as a sequence of random elements of\linebreak
$D([0,+\infty);\R^d)$ equipped with the Skorokhod topology.
\item For any $k\ge1$, $0\le t_1<t_2<\cdots<t_k$, $(X^n_{t_1},\ldots,X^n_{t_k})\Rightarrow
(X_{t_1},\ldots,X_{t_k})$, in the sense of weak convergence in $\R^{d\times k}$.
\end{enumerate}
If only 2 is satisfied, then one has convergence in the sense of finite-dimensional distributions.

What do we mean by tightness? A sequence $\{X_n\}_{n\ge1}$ of random variables with values in a topological space $S$ is said to be tight
if for any $\eps>0$, there exists a compact set $K\subset S$ such that $\P(X_n\in K)\ge1-\eps$ for all $n\ge1$.

Consider the product $X_n Y_n$, where $X_n$ and $Y_n$ are real-valued. If one of the two sequences is tight and the other tends to $0$ in probability, then $X_nY_n\to0$ in probability. This easy result is used in the proof of Theorem \ref{th:CLT}.

In the proof of Lemma \ref{TCL-PoisProc}, we use the following argument: a sequence of continuous time martingales $M^n_t$ satisfying $M^n_0=0$ is
tight as soon as the associated sequence of predictable increasing processes $\langle M^n,M^n\rangle_t$
is $C$--tight, in the sense that both it is tight, and any weak limit of a converging sub--sequence is continuous, see e.g.\ Theorem VI.4.13 in Jacod and Shiryaev \cite{JSI}. In the
situation of Lemma \ref{TCL-PoisProc},
$\langle M^n,M^n\rangle_t=t$ which is $C$--tight, since it does not depend upon $n$ and is continuous.

\section{Pontryagin's maximum principle}\label{se:Pont}
In this section, we present the Pontryagin maximum principle in optimal control, which is useful in order to compute or give some estimates for the exponent
in the asymptotic evaluation of the time to extinction derived from large deviation theory. We refer the reader for a more general presentation, proofs and references to Tr\'elat \cite{TreI} and Pontryagin et al. \cite{PBGMI}.

The quantity of interest, denoted by $\overline{V}$ in Section \ref{sec:textSIRS} and the following pages, is the value function of an optimal control problem which is of the following type. $x\in C([0,\infty);\R^d)$ solves the controlled ODE
\[ \dot{x}_t=Bu_t,\quad x_0=x^\ast,\]
where $B$ is a $d\times k$ matrix, and $u\in L^1([0,\infty);\R^k_+)$ is to be chosen together with the final time $T$ such as to minimize a cost functional
\[ C(u)=\sum_{j=1}^k\int_0^T g(u_j(t),\beta_j(x_t)) dt,\]
while the following constraint must be satisfied: $x_T\in M_1$, where $M_1$ is some affine subspace of $\R^d$. The function $g$ is the one which appears in Section \ref{sec:rate}, namely $g(a,b)=a\log(a/b)-a+b$, while the $\beta_j$'s are some mappings from $\R^d$ into $\R_+$, which, like the matrix $B$, depend upon the particular model we consider. Note that in our case all entries of $B$ are either $1$, $0$, or $-1$.

We associate to this optimal control problem a Hamiltonian which takes the form
\[ H(x,p,u)=\langle p, Bu\rangle -\sum_{j=1}^k g(u_j,\beta_j(x)),\]
where $p\in C([0,T];\R^d)$ is the adjoint state. The next statement constitutes Pontryagin maximum principle, applied to our particular situation.

\begin{theorem}\label{th:pontryag}
If $(\hat{u},\hat{T})$ is an optimal pair, then there exists an adjoint state, such that the following is satisfied
\begin{align*}
\dot{x}_t&= B \hat{u}_t,\quad x_0=x^\ast,\ x_{\hat{T}}\in M_1, \\
\dot{p}_t&=\sum_{j=1}^k\left[\nabla\beta_j(x_t)-\hat{u}_j(t)\frac{\nabla\beta_j(x_t)}{\beta_j(x_t)}\right],\quad p_{\hat{T}}\perp M_1,\\
H(x_t,p_t,\hat{u}_t)&=\max_{v\in\R^k_+}H(x_t,p_t,v)=0,\quad0\le t\le \hat{T}.
\end{align*}
\end{theorem}
Of course, the first equation could be of the more general form $\dot{x}=f(x,u)$. The general form of the adjoint equation reads
$\dot{p}=-\nabla_xH$. The Hamiltonian is zero at time $\hat{T}$
since the final time is not fixed and there is no final cost. The Hamiltonian is constant along the optimal trajectory because none of the coefficients depends upon $t$.
\bigskip

Since $u\to (B^\ast p)_ju-g(u,\beta_j(x))$ is concave, the maximum is the zero of its derivative if it is non-negative. Hence
\[ \hat{u}_j=e^{(B^\ast p)_j}\beta_j(x),\]
and the two above equations can be written as
\begin{align*}
\dot{x}_t=\sum_{j=1}^k e^{(B^\ast p_t)_j}\beta_j(x_t)h_{j},\quad
\dot{p}_t=\sum_{j=1}^k (1-e^{(B^\ast p_t)_j})\nabla\beta_j(x_t),
\end{align*}
and the Hamiltonian along the optimal trajectory reads
\begin{equation}\label{eq:Ha}
 H(x_t,p_t,\hat{u}_t)=\sum_{j=1}^k\beta_j(x_t)(e^{(B^\ast p_t)_j}-1)=0.
 \end{equation}
 \bigskip

Finally the instantaneous cost takes the form
\begin{equation}\label{eq:IC}
\sum_{j=1}^k \left(1\!-\!e^{(B^\ast p_t)_j}\!+\!(B^\ast p_t)_je^{(B^\ast p_t)_j}\right)\beta_j(x_t)=\!\sum_{j=1}^k(B^\ast p_t)_je^{(B^\ast p_t)_j}\beta_j(x_t),
\end{equation}
where this identity follows from \eqref{eq:Ha}.

\section{Semi- and equicontinuity}
Let $\mathcal{X}$ be a metric space, equipped with a distance $d$, and $f$ be a mapping from $\mathcal{X}$ into $\R\cup\{-\infty,\infty\}$.

\begin{definition}\label{de:sc}
$f$ is said to be lower (resp.\ upper) semi-continuous if for any $x_0\in\mathcal{X}$,
\[\liminf_{x\to x_0}f(x)\ge f(x_0)\quad (\text{resp.\ } \limsup_{x\to x_0}f(x)\le f(x_0)).\]
\end{definition}
Clearly $f$ is continuous if and only if it is both lower and upper semi-continuous.

A lower (resp.\ upper) semi-continuous $(-\infty,\infty]$-valued (resp.\ $[-\infty,\infty)$-valued) function achieves its minimum (resp.\ maximum) on a compact subset of $\mathcal{X}$.

The pointwise supremum (resp.\ infimum) of a collection of continuous functions is lower (resp.\ upper) semi-continuous.
\bigskip

Let now $\{f_i,\ i\in I\}$ be a collection of elements of $C(\mathcal{X})$ (i.e.\  of continuous functions from $\mathcal{X}$ into $\R$), where $I$ is an arbitrary index set.

\begin{definition}\label{de:equicont}
The collection $\{f_n,\ n\ge1\}$ is said to be equicontinuous if for any $x_0\in\mathcal{X}$,
$\sup_{i\in I}|f_i(x)-f_i(x_0)|\to0$, as $x\to x_0$. The same collection is said to be uniformly
equicontinuous if $\sup_{i\in I}\sup_{d(x,y)\le\delta}|f_i(x)-f_i(y)|\to0$, as $\delta\to0$.
\end{definition}
Note that when $\mathcal{X}$ is compact, equicontinuity and uniform equicontinuity are equivalent.

\section{Solutions to selected exercises}\label{Sec_Solutions}

\textbf{Solution to Exercise \ref{x2.1}}. $R_0=\lambda E(I)=\lambda/\gamma=1.8$. The escape probability from a given under infected individual equals $E(e^{-\lambda I/N})=\gamma/( \gamma + \lambda/N)$, since $\psi_I(-\lambda/N)=\gamma/(\gamma + \lambda/N)$ when $I\sim \mathrm{Exp}(\gamma)$. For $\lambda=1.8,\ \gamma = 1,\ N=100$ we get 0.9823.

\vskip0.5cm\noindent \textbf{Solution to Exercise \ref{x2.2}}.  For the Reed--Frost epidemic we hence have the same $R_0=\lambda E(I)=1.8$. As for the escape probability we get
\[P(\text{avoid infection from an infective})=e^{-\lambda \iota/N}=0.9822.\]
The escape probabilities are not identical, but very similar for the two models.

\vskip0.5cm\noindent \textbf{Solution to Exercise \ref{xEarly.1}}.
If $I\equiv 1$, then $X\sim \mathrm{Poi}(R_0)$. If $I\sim \mathrm{Exp}(1/\iota )$, then $X\sim \mathrm{MixPoi}(\lambda I)$. So
\[P(X=k)=\int_0^\infty P(X=k| I=s)e^{-s/\iota}/\iota ds = (R_0/(R_0+1))^k (1/(R_0+1)),\]
so $X\sim \mathrm{Geo}(p=1/(R_0+1))$.

\vskip0.5cm\noindent \textbf{Solution to Exercise \ref{xEarly.2}}. The probability of a minor outbreak corresponds to the probability of extinction in the approximating branching process. This probability $q$ was derived in Section \ref{TB-EP_sec_Br-proc}  by conditioning on the number $k$ infected in the first generation, the offspring distribution: if $k$ get infected these all start new independent branching processes so the probability that all go extinct equals $q^k$. The general equation is hence
$$
q=\sum_{k=0}^\infty q^kP(X=k).
$$
The offspring distribution $X$ depends on the infectious period distribution $I$. Given that $I=s$, $X$ has a Poisson distribution with mean $\lambda s$, so $X\sim \mathrm{MixPoi}(\lambda I)$. In situation 2 (cont-time R-F) $I\equiv 1$ so $X\sim \mathrm{Poi}(\lambda=1.5)$. This gives the following equation
$$
q=\sum_{k=0}^\infty q^k\frac{\lambda^k e^{-\lambda}}{k!}= ... =e^{-R_0(1-q)}=e^{-1.5(1-q)}.
$$
If this equation is solved numerically it gives the result that $q=1-0.583=0.417$. So for the Reed--Frost case the probability of a major outbreak, equals 0.583.

As for the Markovian SIR, where $I\sim \mathrm{Exp}(1)$ we get
\begin{align*}
P(X=k)&=\int_0^\infty P(X=k|I=s)f_I(s)ds\\
 &= \int_0^\infty \frac{(\lambda s)^k e^{-\lambda s}}{ k!}e^{-s}ds=\cdots= \frac{1}{1+\lambda}\left( \frac{\lambda }{1+\lambda}\right)^k
\end{align*}
i.e.\ the geometric distribution, which should not come as a surprise (each time, the event is either infection or recovery, and the latter has probability $1/(\lambda + 1)$). We then get
$$
q=\sum_{k=0}^\infty q^kP(X=k) =q^k\left( \frac{\lambda}{\lambda + 1} \right)^k\frac{1}{\lambda + 1}=\frac{1}{1+(1-q)\lambda }.
$$
As a consequence, the probability of a minor outbreak for the Markovian SIR hence equals $q=1/\lambda=1/R_0=1/1.5=0.67$. The probability of a major outbreak is hence only 0.33. The randomness of the infectious period hence reduces the risk for a major outbreak. It can actually be proven that having a constant infectious period maximizes the outbreak probability among all distributions of the infectious period.

\vskip0.5cm\noindent \textbf{Solution to Exercise \ref{xEarly.3}}.  The exponential growth rate (or decay rate if $R<1$) $r$ is the solution to Equation (\ref{Malthus-eq}), where $h (s)$ is the average rate of infectious contacts $s$ units after infection: $h(s)=\lambda P(L\le s\le L+I)$. For the Markovian SIR (for which $L\equiv 0$ and $I\sim \mathrm{Exp}(\gamma=1/\iota)$) we hence have $h (s)=\lambda e^{-s/\iota}=1.5e^{-s}$, and the solution equals $r=\lambda - 1/\iota$ For $R_0=1.5$ and $\gamma=\iota=1$ this gives the exponential growth rate $r=0.5$.

For the continuous time Reed--Frost model we have $h (s)=\lambda 1_{(s<\iota )}$. The equation then becomes $\int_0^\iota e^{-rs}\lambda ds =\frac{\lambda}{r}\left(1-e^{-r\iota}\right) =1$. The equation is hence $r/\lambda =1-e^{-(r/\lambda)R_0}$.  When $R_0=1.5$ we numerically get $r/\lambda =0.583$, so $r=0.874$ for the continuous time Reed--Frost model. This epidemic hence grows quicker than the Markovian SIR epidemic with the same parameters. The main reason for this is that even if the two infectious periods have equal mean $\iota =1$, the average time of the infectious contacts are not the same. For the Reed--Frost the mean time to a randomly selected infectious contact (the mean of the generation time distribution) is of course 0.5 (the generation time distribution is uniform on $[0,\ 1]$, whereas for the Markovian SIR it equals 1 (the generation time distribution is $\mathrm{Exp}(1)$).

For the third case, with exponentially distributed latency and infectious periods,  we have $h(s)=P(L<s<L+I) = \frac{\lambda\nu}{\gamma-\nu}\left( e^{-\nu s}-e^{-\gamma s}\right)$. Solving $\int_0^\infty e^{-rs}h(s)ds =1$ gives the solution
$$
r=\sqrt{\nu (\lambda-\gamma) + \left(\frac{\gamma+\nu}{2}\right)^2 } - \frac{\gamma+\nu}{2} \approx 0.2247.
$$
Of course, adding a latency period before the infectious period will reduce the growth rate $r$ of the epidemic.

\vskip0.5cm\noindent \textbf{Solution to Exercise \ref{xvacc-numeric}}.
$v_c=1-1/R_0=0.5$. When $v=0.33$, $z_v$ solves the equation $1-z_v=e^{-(1-v)R_0z_v}$, and the numerical solution equals $z_v=0.4544$. The over-all fraction infected is hence $(1-v)z_v=0.3029$. As for the probability of a major outbreak we have that for the Markovian SIR $P(\text{major outbreak}) = 1-1/R_v=0.25$, since $R_v=(1-v)R_0=0.67\cdot 2=1.33$.

\vskip0.5cm\noindent \textbf{Solution to Exercise \ref{xnon-perf-vacc}}.
The new rate at which an infectious individual makes infectious contacts when $v=33\%$ are vaccinated is $\lambda' =\lambda p v + \lambda (1-v)$ where $p=0.2$ (this is true irrespective of whether the infector was vaccinated or not). Since the average infectious period equals $E(I)=1$ we have $R_v =\lambda'E(I)=1.467$ (instead of $R_0=2$ when no one is vaccinated).

\vskip0.5cm\noindent \textbf{Solution to Exercise \ref{x4.2}}. $R_0=1.5$: 0.583,  $R_0=3$: 0.940, $R_0=15$: 1.000 (of course not exactly, but to this precision).

\vskip0.5cm\noindent \textbf{Solution to Exercise \ref{x3.1}}.
\begin{align*}
p_3^{(3)} &=p^3 + \binom{3}{2}p^2(1-p)*(1-(1-p)^2)
\\
&\hskip2cm + \binom{3}{1}p(1-p)^2 *p^2 + \binom{3}{1}p(1-p)^2 * \binom{2}{1}p(1-p) * p.
\end{align*}

\vskip0.5cm\noindent \textbf{Solution to Exercise \ref{xCLT1}}
The limiting mean equals $Nz$ where $z$ solves $1-z=e^{-R_0z}$ so with $R_0=\lambda\iota=1.5$ we get $z=0.583$ and the limiting mean equals 583 for both scenarios. The limiting variance of $Z^N$ equals $N\frac{z(1-z)(1+r^2(1-z)R_0^2)}{(1-(1-z)R_0)^2}$, where $r$ is the coefficient of variation of the infectious period. For the Reed--Frost case with non-random infectious period we have $r=0$ implying that the limiting variance equals 1737, so the standard deviation equals 41.7, so one can expect that the final size will be somewhere in the interval $583\pm 80$ with about 95\% probability. The Markovian SIR has exponential infectious period which has $r=1$ giving a variance of 3367 and standard deviation 58.0. So, the fact that the infectious period is exponential as compared to fixed makes the standard deviation of the final size increase by close to 50\%.

\vskip0.5cm\noindent \textbf{Solution to Exercise \ref{xSubcrit}}.
The numerical values are: the final size equals $z=0.583$ and $R_0(1-z)=0.626 < 1$.

\vskip0.5cm\noindent \textbf{Solution to Exercise \ref{xDur.1}}.
Computing the two leading terms is equivalent to computing $r$ and $r^*$. For the Markovian SIR we have $r=0.5$ and $r^*=-0.3742$, for the continuous time Reed--Frost we get $r=0.8742$ and $r^*=-0.8741$, and for the Markovian SEIR we have $r=0.2247$ and $r^*=-0.2089$.

\vskip0.5cm\noindent \textbf{Solution to Exercise \ref{Nasell}}.
Denoting by $U(t)$ the vector of the Gaussian fluctuations around
$\begin{pmatrix}s(t)\\ i(t)\end{pmatrix}$, deduce from Theorem \ref{th:CLT}
that this vector solves the linear SDE
\[ U(t)=\int_0^t A(r)U(r)dr+\int_0^t C(r)dB_r,\]
where $B(t)$ is a standard five-dimensional Brownian motion and
\begin{align*}
 A(t)&=\mu \begin{pmatrix} -1-\frac{R_0}{\varepsilon}i(t) & -\frac{R_0}{\varepsilon}s(t)\\ \frac{R_0}{\varepsilon}i(t)& \varepsilon^{-1}(R_0s(t)-1)\end{pmatrix},\\
C(t)&=\mu \begin{pmatrix}\sqrt{\mu}&-\sqrt{\frac{\mu R_0}{\varepsilon}s(t)i(t)}&-\sqrt{\mu s(t)}&0\\
0&\sqrt{\frac{\mu R_0}{\varepsilon}s(t)i(t)}&0&-\sqrt{\frac{\mu}{\varepsilon}i(t)}\end{pmatrix}.
\end{align*}
Show that, as $t\to\infty$,
\begin{align*}
A(t)\to \mu\begin{pmatrix} -R_0&-1/\varepsilon\\ R_0-1&0\end{pmatrix}\!, \
C(t)C^\ast(t)\to\frac{\mu}{R_0}\begin{pmatrix}2R_0&-(R_0-1)\\-(R_0-1)&2(R_0-1)\end{pmatrix}.
\end{align*}
Show that the eigenvalues of $A=\lim_{t\to\infty}A(t)$ are complex, as soon as $\varepsilon<4/R_0$, and that the real parts of those eigenvalues are negative.
Conclude from a combination of Exercise  \ref{exerOU} and Lemma \ref{le:cov_inv} that the covariance matrix of the stationary distribution of $U(t)$ reads
\[\begin{pmatrix} \frac{1}{R_0}+\frac{1}{\varepsilon R_0^2}&-\frac{1}{R_0}\\
-\frac{1}{R_0}&\frac{1}{R_0}-\frac{1}{R_0^2}+\varepsilon\end{pmatrix}.\]
Conclude by taking into account that we expect to have $\varepsilon<<R_0^{-1}$.

\vskip0.5cm\noindent \textbf{Solution to Exercise \ref{xopen-end-level1}}. The relative length of the infectious period $\varepsilon$ affects the critical community size $N_c$  much more than $R_0$ does, since it is squared in the approximation of $N_c$. As an illustration, if the infectious period is doubled (with half infectivity per unit of time thus keeping $R_0$ fixed)  $N_c$ will decrease by a factor 4, whereas if the basic reproduction number is doubled (keeping everything else fixed) only decreases $N_c$ by a factor close to 2.

\vskip0.5cm\noindent \textbf{Solution to Exercise \ref{xopen-end-level2}}. There are two effects of this vaccination strategy. The first is that vaccinated individuals can be ignored, so the relevant population (of unvaccinated people) is now $N^{(unvacc)}=N(1-v)$. Secondly, since infected individuals have contact with both types of individuals, the rate of having contact with the population of interest is reduced to $\lambda (1-v)$ implying that the reproduction number is changed to $R_v=R_0(1-v)$. The critical population size of unvaccinated people $N_c^{(unvacc)}$ is then simply obtained in the same way, but for these new parameters, so
$$
N_c^{(unvacc)}= \frac{9}{ \varepsilon^2( 1-\frac{1}{R_v} )^2 R_v} = \frac{9}{ \varepsilon^2( 1-\frac{1}{(1-v)R_0} )^2 (1-v)R_0}.
$$
However, a more interesting quantity is the critical community size counting all individuals, hence also vaccinated. Since $N=N^{(unvacc)}(/1-v)$, the critical community size for a population in which a fraction $v$ of the new-born are continuously being vaccinated is given by
$$
N_c^{(v)}=  \frac{9}{ (1-v)^2\varepsilon^2( 1-\frac{1}{(1-v)R_0} )^2R_0}.
$$
By numerical studies it is easily shown that the critical community size grows very big with $v$, also agreeing with empirical evidence since e.g.\ measles is no longer endemic in England (or anywhere else in the world having high vaccination coverage).

\newpage

\pagestyle{referencesI}

\renewcommand{\bibsection}{\chapter*{References for Part I}}

\renewcommand{\thechapter}{\arabic{chapter}}





\backmatter


\end{document}